\documentclass[shortAfour,sageh,times]{sagej}
\usepackage{hyperref}
\usepackage{amsmath}
\usepackage{amssymb}
\usepackage{mathtools}
\usepackage[algoruled,lined,longend,nofillcomment]{algorithm2e}
\usepackage{multirow}
\usepackage{booktabs}
\newcommand{\assgn}{\ensuremath\mathrel{\mathop:}=}
\newcommand{\diag}{\ensuremath\mathop{\mathrm{diag}}}
\newcommand{\sign}{\ensuremath\mathop{\mathrm{sign}}}
\newcommand{\Real}{\ensuremath\mathop{\mathrm{Re}}}
\newcommand{\Imag}{\ensuremath\mathop{\mathrm{Im}}}
\def\volumeyear{2020}
\newlength\fbw
\begin{document}
\runninghead{Novakovi\'{c} and Singer}
\title{Implicit Hari--Zimmermann algorithm for the generalized SVD on the GPUs}
\author{Vedran Novakovi\'{c}\affilnum{1} and Sanja Singer\affilnum{2}}
\affiliation{\affilnum{1}Completed a major part of this research while
  being affiliated to Universidad Jaime I, Av.~Vicent Sos Baynat,
  12071 Castell\'{o}n de la Plana, Spain\\
  \affilnum{2}University of Zagreb, Faculty of Mechanical Engineering
  and Naval Architecture, Ivana Lu\v{c}i\'{c}a 5, 10000 Zagreb,
  Croatia}
\corrauth{Sanja Singer, University of Zagreb, Faculty of Mechanical
  Engineering and Naval Architecture, Ivana Lu\v{c}i\'{c}a 5, 10000
  Zagreb, Croatia}
\email{ssinger@fsb.hr}
\begin{abstract}
  A parallel, blocked, one-sided Hari--Zimmermann algorithm for the
  generalized singular value decomposition (GSVD) of a real or a
  complex matrix pair $(F,G)$ is here proposed, where $F$ and $G$ have
  the same number of columns, and are both of the full column rank.
  The algorithm targets either a single graphics processing unit
  (GPU), or a cluster of those, performs all non-trivial computation
  exclusively on the GPUs, requires the minimal amount of memory to be
  reasonably expected, scales acceptably with the increase of the
  number of GPUs available, and guarantees the reproducible, bitwise
  identical output of the runs repeated over the same input and with
  the same number of GPUs.
\end{abstract}
\maketitle
%
%
\section{Introduction}\label{s:1}
%
%
The two-sided Hari--Zimmermann
algorithm~\cite{Hari-84,Hari-2018,Hari-2019,Zimmermann-69} is a
Jacobi-type method for computing the generalized eigenvalue
decomposition (GEVD) of a matrix pair $(A,B)$, where both matrices are
Hermitian of the same order and $B$ is positive definite.

If $A$ and $B$ are instead given \emph{implicitly\/} by their factors
$F$ and $G$ (not necessarily square nor with the same number of rows),
respectively, such that $(A,B) = (F^{\ast} F, G^{\ast} G)$, then the
GEVD of $(A,B)$ can be computed implicitly, i.e., without assembling
$A$ and $B$ in entirety from the factors, by a modification of the
Hari--Zimmermann algorithm~\cite{Novakovic-Singer-Singer-2015}.
However, pivot submatrices of $A$ and $B$ of a certain, usually small
order are formed explicitly throughout the computation.

The modified algorithm is a method that jointly orthogonalizes the
pairs of columns of $F$ and $G$ by a sequence of transformations that
are applied from the right side of the factors only.  Such a one-sided
algorithm computes $U$, $\Sigma_F^{}$, $V$, $\Sigma_G^{}$, and $Z$,
where $F Z = U \Sigma_F^{}$, $G Z = V \Sigma_G^{}$, and
$U^{\ast} U^{} = V^{\ast} V = I$.  The matrix $Z$ is square and
nonsingular, while $\Sigma_F^{}$ and $\Sigma_G^{}$ are non-negative,
diagonal, and scaled such that $\Sigma_F^2 + \Sigma_G^2 = I$.  The
method thus implicitly computes the GEVD of $(A,B)$, but
\emph{explicitly\/} the generalized singular value decomposition
(GSVD; see, e.g., \cite{Paige-Saunders-81,VanLoan-76}) of $(F,G)$,
with the generalized singular values forming the diagonal of
$\Sigma \assgn \Sigma_G^{-1} \Sigma_F^{}$ (all of them finite, since
$\Sigma_G^{}$ has a positive diagonal).  Furthermore, the generalized
singular values can be considered to be sorted descendingly by a
symmetric permutation, i.e., $\Sigma = P_0^T \Sigma' P_0^{}$, and thus
$U = U' P_0^{}$, $V = V' P_0^{}$, and $Z = Z' P_0^{}$, where
$F Z' = U' \Sigma_F'$, $G Z' = V' \Sigma_G'$, and
$\Sigma' = \Sigma_G^{\prime -1} \Sigma_F'$ constitute a decomposition
of $(F,G)$ possessing all other aforementioned properties.

\looseness=-1
The GEVD of $(A,B)$, if required, can be recovered by letting
$\Lambda \assgn \Sigma^2$ and noting that $A Z = B Z \Lambda$, i.e.,
the columns of $Z$ are the generalized eigenvectors, and the diagonal
of $\Lambda$ contains the generalized eigenvalues of $(A,B)$.
However, the converse is \emph{not\/} numerically sound, i.e., the
GEVD should not, in general, be used for computing the GSVD\@.  For a
further clarification, see Appendix~\ref{s:G}.

The right generalized singular vectors $X \assgn Z^{-1}$, if needed,
can either be computed from $Z$, or can be obtained simultaneously
with $Z$ by accumulating the inverses of the transformations that have
been multiplied to form
$Z$~\cite{Singer-DiNapoli-Novakovic-Caklovic-2020}.  With
$\widetilde{\Theta}$ from subsection~\ref{ss:2.4}, if
\begin{displaymath}
  Z = Z_0^{} \widetilde{Z} \widetilde{\Theta} = Z_0^{} \cdot \widetilde{Z}_0^{} \cdot \widetilde{Z}_1^{} \cdots \widetilde{Z}_N^{} \cdot \widetilde{\Theta},
\end{displaymath}
when $N+1$ transformations have been applied, then
\begin{displaymath}
  X = \widetilde{\Theta}^{-1} \cdot \widetilde{Z}_N^{-1} \cdot \widetilde{Z}_{N-1}^{-1} \cdots \widetilde{Z}_0^{-1} \cdot Z_0^{-1}.
\end{displaymath}

The recent work~\cite{Novakovic-Singer-Singer-2015} has shown that
such method can be blocked and parallelized for the shared memory
nodes and for the clusters of those, albeit only the real matrix pairs
were considered therein.  Even the sequential but blocked version
outperformed the GSVD algorithm in LAPACK~\cite{Anderson-et-al-99},
and the parallel ones exhibited a decent scalability.

On the other hand, an efficient, parallel and blocked one-sided
Jacobi-type algorithm for the ``ordinary'' and the hyperbolic
SVD~\cite{Novakovic-2015,Novakovic-2017} of a single real matrix has
been developed for the GPUs, that utilizes the GPUs almost fully, with
the CPU serving only the controlling purpose in the single-GPU case.

This work aims to merge the experience of those two approaches, and
present a parallel and blocked one-sided (also called ``implicit'')
Hari--Zimmermann algorithm for the GSVD on the GPU(s) as an extension
of the latter, but for the complex matrix pairs as well as for the
real ones.

Even though the research in parallelization of the GSVD has a long
history~\cite{Bai-94,Luk-85}, three novel and major differences from
the earlier, Kogbetliantz-based procedures aim to ensure both the high
performance and the high relative accuracy of this one: using the
implicit Hari--Zimmermann algorithm as the basic method, that is
blocked to exploit the GPU memory hierarchy, and the massive
parallelism of the GPUs that suits the algorithm (and vice versa)
perfectly.

In the last twenty years, many applications of GSVD have been found in
science and technology.  To mention just a few applications, the GSVD
is used for dimension reduction for clustered text
data~\cite{Howland-Jeon-Park-2003} and for face recognition
algorithms~\cite{Howland-Wang-Park-2006}, where in both cases the
matrix pair is naturally given implicitly, i.e., in a factored form.

In~\cite{Alter-Brown-Botstein-2003} the GSVD serves for comparison of
two different organisms to find their biological similarities based on
a genome-scale expression data sets.  Also, the GSVD can be used in
beamforming~\cite{Senaratne-Tellambura-2013} and separation of
partially overlapping data packets~\cite{Zhou-vanderVeen-2017} in
communication systems, machine condition monitoring when looking for
symptoms of wear~\cite{Cempel-2009}, and filtering of brain activities
while preforming two different tasks~\cite{Zhao-et-al-2010}.  In the
last case, matrices could be very large.

This paper continues with section~\ref{s:2}, where the complex and the
real one-sided Hari--Zimmermann algorithms are introduced, together
with the general, architecturally agnostic principles of their
blocking and parallelization.  In section~\ref{s:3} the single-GPU
implementation are described in detail, while in section~\ref{s:4} the
most straightforward multi-GPU implementation approach is suggested.
The numerical testing results are summarized in section~\ref{s:5}, and
the paper concludes with some directions for future research in
section~\ref{s:6}.  In Appendix~\ref{s:A} a non-essential method for
enhancing the accuracy of the real and the complex dot-products on the
GPUs is proposed.
%
%
\section{The complex and the real one-sided Hari--Zimmermann algorithms}\label{s:2}
%
%
In this section the complex and the real one-sided Hari--Zimmermann
algorithms are briefly described.  Please
see~\cite{Hari-84,Hari-2018,Hari-2019} for a more thorough overview of
the two-sided algorithms, and~\cite{Novakovic-Singer-Singer-2015} for
a detailed explanation of the real implicit Hari--Zimmermann
algorithm.  In this paper the terminology and the
implementation decisions
of~\cite{Singer-DiNapoli-Novakovic-Caklovic-2020}, where the complex
generalized hyperbolic SVD based on the implicit Hari--Zimmermann
approach has been introduced, are closely followed, but without the
hyperbolic scalar products (i.e., the signature matrix $J$ is taken to
be identity here) and without forming the right generalized singular
vectors $X$ from $Z$.

Let the matrices $F$ and $G$ be of size $m_F^{} \times n$ and
$m_G^{} \times n$, respectively, with $\min\{m_F^{}, m_G^{}\}\ge n$.
Then, $Z$ is square of order $n$, and assume that $n \ge 2$.
Otherwise, for $n = 1$, the GSVD of $(F,G)$ is obtained by taking
\begin{displaymath}
  \begin{aligned}
    U & \assgn \|F\|_F^{-1} F,\\
    V & \assgn \|G\|_F^{-1} G,\\
    Z & \assgn \frac{1}{\sqrt{\|F\|_F^2 + \|G\|_F^2}},
  \end{aligned}
  \qquad
  \begin{aligned}
    \Sigma_F^{} & \assgn \frac{\|F\|_F^{}}{\sqrt{\|F\|_F^2 + \|G\|_F^2}},\\
    \Sigma_G^{} & \assgn \frac{\|G\|_F^{}}{\sqrt{\|F\|_F^2 + \|G\|_F^2}}.
  \end{aligned}
\end{displaymath}

Even though the algorithm works on the rectangular matrices, it might
be beneficial performance-wise to avoid transforming very tall and
skinny (block)columns by working on the square matrices instead.  To
shorten $F$ and $G$, the problem is transformed by computing the QR
factorization of $F$ with the column pivoting,
$F P_1^{} = Q_F^{} R_F^{}$, and then $G$, with its columns prepermuted
by $P_1^{}$, is shortened by the column-pivoted QR factorization,
$(G P_1^{}) P_2^{} = Q_G^{} R_G^{}$.  The square matrices
$F'' \assgn R_F^{} P_2^{}$ and $G'' \assgn R_G^{}$, both of order $n$,
take the place of $F$ and $G$ in the algorithm, respectively.  With
$\Sigma = \Sigma''$ in the decompositions of $(F,G)$ and of
$(F'',G'')$, the matrix $Z$ from the former, sought-for decomposition
can be recovered by using $P'' \assgn P_1^{} P_2^{}$ and the computed
$Z''$ from the latter as $Z \assgn P'' Z''$.

It is assumed that $\diag(B) = I$, i.e., the column norms of $G$ are
unity.  Should it not be the case, $F$ and $G$ are then prescaled by
a nonsingular, diagonal matrix $Z_0^{}$, where
$(Z_0^{})_{\null\!jj} \assgn 1/\!\sqrt{\|g_{\null\!j}^{}\|_F^{}}$,
$g_{\null\!j}^{}$ is the $j$th column of $G$ and $1 \le j \le n$;
otherwise, $Z_0^{} \assgn I$.  The iterative transformation phase
starts with the matrix pair $(F_0^{}, G_0^{})$, where
$F_0^{} \assgn F Z_0^{}$, and $G_0^{} \assgn G Z_0^{}$.
Implicitly, $A$ and $B$ have been transformed by a congruence with
$Z_0^{}$ as $A_0^{} \assgn F_0^{\ast} F_0^{}$ and
$B_0^{} \assgn G_0^{\ast} G_0^{}$.
%
%
\subsection{Simultaneous diagonalization of a pair of pivot matrices}\label{ss:2.1}
%
%
An iteration (or ``step'') $k \ge 0$ of the sequential non-blocked
Hari--Zimmermann algorithm consists of selecting a pair of indices
$(i_k^{},j_k^{})$, $1 \le i_k^{} < j_k^{} \le n$, and thus two
$2 \times 2$ pivot submatrices, one of
$A_k^{} \assgn F_k^{\ast} F_k^{}$,
\begin{displaymath}
  \widehat{A}_k^{} \assgn
  \begin{bmatrix}
    a_{i_k^{}i_k^{};k}^{} & a_{i_k^{}j_k^{};k}^{}\\
    \bar{a}_{i_k^{}j_k^{};k}^{} & a_{j_k^{}j_k^{};k}^{}
  \end{bmatrix} =
  \begin{bmatrix}
    f_{i_k^{};k}^{\ast} f_{i_k^{};k}^{} & f_{i_k^{};k}^{\ast} f_{j_k^{};k}^{}\\
    f_{j_k^{};k}^{\ast} f_{i_k^{};k}^{} & f_{j_k^{};k}^{\ast} f_{j_k^{};k}^{}
  \end{bmatrix},
\end{displaymath}
and one of $B_k^{} \assgn G_k^{\ast} G_k^{}$,
\begin{displaymath}
  \widehat{B}_k^{} \assgn
  \begin{bmatrix}
    1 & b_{i_k^{}j_k^{};k}^{} \\
    \bar{b}_{i_k^{}j_k^{};k}^{} & 1
  \end{bmatrix} =
  \begin{bmatrix}
    1 & g_{i_k^{};k}^{\ast} g_{j_k^{};k}^{}\\
    g_{j_k^{};k}^{\ast} g_{i_k^{};k}^{} & 1
  \end{bmatrix},
\end{displaymath}
which are then jointly diagonalized by a congruence transformation
with a nonsingular matrix $\widehat{Z}_k^{}$, to be defined in
subsections~\ref{sss:2.1.1} and \ref{sss:2.1.2}, as
\begin{displaymath}
  \begin{aligned}
    \widehat{A}_{k+1}^{}&\assgn
    \widehat{Z}_k^{\ast} \widehat{A}_k^{} \widehat{Z}_k^{} =
    \begin{bmatrix}
      a_{i_k^{}i_k^{};k+1} & 0\\
      0 & a_{j_k^{}j_k^{};k+1}
    \end{bmatrix},\\
    \widehat{B}_{k+1}^{}&\assgn
    \widehat{Z}_k^{\ast} \widehat{B}_k^{} \widehat{Z}_k^{} = I_2^{}.
  \end{aligned}
\end{displaymath}

If $\widehat{Z}_k^{}$ is embedded into an $n \times n$ matrix
$\widetilde{Z}_k^{}$ such that
$\widetilde{Z}_{i_k^{}i_k^{};k}^{} \assgn \widehat{Z}_{11;k}^{}$,
$\widetilde{Z}_{i_k^{}j_k^{};k}^{} \assgn \widehat{Z}_{12;k}^{}$,
$\widetilde{Z}_{j_k^{}i_k^{};k}^{} \assgn \widehat{Z}_{21;k}^{}$,
$\widetilde{Z}_{j_k^{}j_k^{};k}^{} \assgn \widehat{Z}_{22;k}^{}$,
while letting $\widetilde{Z}_k^{}$ be the identity matrix elsewhere,
then looking two-sidedly the congruence with $\widetilde{Z}_k^{}$
transforms the pair $(A_k^{},B_k^{})$ into a pair
$(A_{k+1}^{},B_{k+1}^{})$, where
$A_{k+1}^{} \assgn \widetilde{Z}_k^{\ast} A_k^{} \widetilde{Z}_k^{}$
and 
$B_{k+1}^{} \assgn \widetilde{Z}_k^{\ast} B_k^{} \widetilde{Z}_k^{}$.
One-sidedly, the transformation by $\widetilde{Z}_k^{}$ orthogonalizes
the $i_k^{}$th and the $j_k^{}$th pivot columns of $F_k^{}$ and
$G_k^{}$ to obtain $F_{k+1}^{} \assgn F_k^{} \widetilde{Z}_k^{}$ and
$G_{k+1}^{} \assgn G_k^{} \widetilde{Z}_k^{}$.  Also,
$\widetilde{Z}_k^{}$ is accumulated into the product
$Z_{k+1}^{} \assgn Z_k^{} \widetilde{Z}_k^{}$.  In a one-sided
sequential step only the $i_k^{}$th and the $j_k^{}$th columns of
$F_k^{}$, $G_k^{}$, and $Z_k^{}$ are effectively transformed, in-place
(i.e., overwritten), postmultiplying them by the $2 \times 2$ matrix
$\widehat{Z}_k^{}$, while the other columns of these matrices remain
intact:
\begin{displaymath}
  \begin{aligned}
    \begin{bmatrix}
      f_{i_k^{};k+1}^{} & f_{j_k^{};k+1}^{}
    \end{bmatrix} & =
    \begin{bmatrix}
      f_{i_k^{};k}^{} & f_{j_k^{};k}^{}
    \end{bmatrix}
    \cdot \widehat{Z}_k^{},\\
    \begin{bmatrix}
      g_{i_k^{};k+1}^{} & g_{j_k^{};k+1}^{}
    \end{bmatrix} & =
    \begin{bmatrix}
      g_{i_k^{};k}^{} & g_{j_k^{};k}^{}
    \end{bmatrix}
    \cdot \widehat{Z}_k^{},\\
    \begin{bmatrix}
      z_{i_k^{};k+1}^{} & z_{j_k^{};k+1}^{}
    \end{bmatrix} & =
    \begin{bmatrix}
      z_{i_k^{};k}^{} & z_{j_k^{};k}^{}
    \end{bmatrix}
    \cdot \widehat{Z}_k^{}.
  \end{aligned}
\end{displaymath}

As $\diag(\widehat{B}_{k+1}^{}) = \diag(\widehat{B}_k^{}) = I_2^{}$,
it follows that $\diag(B_{k+1}^{}) = \diag(B_k^{}) = I_n^{}$.
However, due to the floating-point rounding errors, these equations
might not hold.  To prevent $\diag(\widehat{B}_k^{})$ to drift too far
away from $\diag(I_2^{})$ as the algorithm progresses, the squared
Frobenius norms of $g_{i_k^{};k}^{}$ and $g_{j_k^{};k}^{}$ could be
recomputed for each $k$ as
$b_{i_k^{}i_k^{};k}^{} = g_{i_k^{};k}^{\ast} g_{i_k^{};k}^{}$ and
$b_{j_k^{}j_k^{};k}^{} = g_{j_k^{};k}^{\ast} g_{j_k^{};k}^{}$.  Then,
a rescaling of $\widehat{A}_k^{}$ and $\widehat{B}_k^{}$ as
$\widehat{A}_k' \assgn \widehat{D}_k^{\ast} \widehat{A}_k^{} \widehat{D}_k^{}$
and
$\widehat{B}_k' \assgn \widehat{D}_k^{\ast} \widehat{B}_k^{} \widehat{D}_k^{}$,
by a diagonal matrix $\widehat{D}_k^{}$ such that
$\widehat{D}_{11;k}^{} = 1/\!\sqrt{b_{i_k^{}i_k^{};k}^{}}$ and
$\widehat{D}_{22;k}^{} = 1/\!\sqrt{b_{j_k^{}j_k^{};k}^{}}$, should
bring back $\diag(\widehat{B}_k^{})$ close to $\diag(I_2^{})$.  From
$\widehat{A}_k'$ and $\widehat{B}_k'$ it is then possible to compute
$\widehat{Z}_k'$, with the final
$\widehat{Z}_k^{} \assgn \widehat{D}_k \widehat{Z}_k'$.  In this
version of the algorithm it is not necessary to rescale the columns of
$F$ and $G$ by $\widetilde{Z}_0^{}$ at the start, since such rescaling
happens at each step, so $\widetilde{Z}_0^{} \assgn I$.  If
$\widehat{D}_k^{} = I_2^{}$, this version is equivalent to the
standard (previously described) one, for which it can be formally set
$\widehat{A}_k' \assgn \widehat{A}_k^{}$ and
$\widehat{B}_k' \assgn \widehat{B}_k^{}$.

Suppose that $\widehat{Z}_k''$ has been computed (by either version)
such that it diagonalizes $\widehat{A}_k^{}$ and $\widehat{B}_k^{}$,
but $a_{i_k^{}i_k^{};k+1}^{} < a_{j_k^{}j_k^{};k+1}^{}$.  To keep
$\diag(\widehat{A}_k^{})$ sorted descendingly, swap the columns of
$\widehat{Z}_k''$ by a permutation
$\widehat{P}_k^{} \assgn \begin{bsmallmatrix} 0 & 1 \\ 1 & 0 \end{bsmallmatrix}$
to obtain $\widehat{Z}_k^{} \assgn \widehat{Z}_k'' \widehat{P}_k^{}$.
Such $\widehat{Z}_k^{}$ will swap the $i_k^{}$th and the $j_k^{}$th
columns of $F_k^{}$ and $G_k^{}$ as it orthogonalizes them.  Sorting
in each step is a heuristic that speeds up the algorithm notably in
practice (see section~\ref{s:5}), but it makes reasoning about the
convergence harder and is not strictly necessary.

Computing $\widehat{Z}_k^{}$ from $\widehat{A}_k^{}$ and
$\widehat{B}_k^{}$ is more involved in the complex case than
in the real one.  However, in both cases, first it is established
whether the $i_k^{}$th and the $j_k^{}$th columns of $F_k^{}$ and
$G_k^{}$ are numerically relatively orthogonal,
\begin{displaymath}
  \begin{aligned}
    |a_{i_k^{}j_k^{};k}'|&<\sqrt{a_{i_k^{}i_k^{};k}'}\cdot\sqrt{a_{j_k^{}j_k^{};k}'}\cdot\varepsilon\cdot\sqrt{n},\\
    |b_{i_k^{}j_k^{};k}'|&<\varepsilon\cdot\sqrt{n},
  \end{aligned}
\end{displaymath}
where $\varepsilon$ is the precision of the chosen floating-point
datatype.  The relation relies on the expected (as opposed to the
worst case) rounding error for the dot-products~\cite{Drmac-97} that
form the elements of $\widehat{A}_k'$ and $\widehat{B}_k'$, and while
sensible in the real case, it is probably too tight in the complex
case, where a more careful analysis of the complex dot-products might
be employed in the future work and a handful of transformations
subsequently might be skipped.  If the aforesaid columns are
relatively orthogonal, no non-trivial transformation is to take place,
and $\widehat{Z}_k^{} \assgn \widehat{P}_k^{}$, since still the column
swap may be warranted.  Rescaling by $\widehat{D}_k^{}$ is thus not
performed even for $\widehat{D}_k^{} \ne I_2^{}$, since it might
perturb the columns sufficiently enough for them to cease to be
numerically orthogonal.
%
%
\subsubsection{The complex case}\label{sss:2.1.1}
%
%
The transformation matrix $\widehat{Z}_k'$ is sought in a
form~\cite{Hari-84,Singer-DiNapoli-Novakovic-Caklovic-2020}
\begin{displaymath}
  \widehat{Z}_k' \assgn \frac{1}{t_k^{}}
  \begin{bmatrix}
    \hphantom{-e^{-i\beta_k^{}}}\cos\varphi_k^{} & e^{i\alpha_k^{}}\sin\varphi_k^{}\\
    -e^{-i\beta_k^{}}\sin\psi_k^{} & \hphantom{e^{i\alpha_k^{}}}\cos\psi_k^{}
  \end{bmatrix}.
\end{displaymath}
To that end, let $\displaystyle x_k^{} \assgn |b_{i_k^{}j_k^{};k}'|$,
$\displaystyle \zeta_k^{} \assgn \arg(b_{i_k^{}j_k^{};k}')$, or
$\zeta_k^{} \assgn 0$ if $b_{i_k^{}j_k^{};k}' = 0$,
$\displaystyle z_k^{} \assgn e^{-i\zeta_k^{}}a_{i_k^{}j_k^{};k}'$, and
define $\sign(a,b)$ to be $|a|$ with the sign of $b$ for $a$ and $b$
real.  Then, let $t_k^{} \assgn \sqrt{1 - x_k^2}$, set
\begin{displaymath}
  \begin{aligned}
    u_k^{} & \assgn \Real(z_k^{}),\\
    v_k^{} & \assgn \Imag(z_k^{}),
  \end{aligned}
  \qquad
  \begin{aligned}
    h_k^{} & \assgn a_{j_k^{}j_k^{};k}' - a_{i_k^{}i_k^{};k}',\\
    \tau_k^{} & \assgn \sign(1,h_k^{}),
  \end{aligned}
\end{displaymath}
and, noting that $t_k^{} > 0$ since $\widehat{B}_k'$ is positive
definite, with these quantities compute
\begin{displaymath}
  \begin{aligned}
    \tan(2\vartheta_k^{})&\assgn \tau_k^{}\frac{2u_k^{} - (a_{i_k^{}i_k^{};k}' + a_{j_k^{}j_k^{};k}')x_k^{}}{t_k^{}\sqrt{h_k^2 + 4v_k^2}},\\
    \tan\gamma_k^{}&\assgn 2\frac{v_k^{}}{h_k^{}},
  \end{aligned}
\end{displaymath}
where $-\pi/4 < \vartheta_k^{} \le \pi/4$ and
$-\pi/2 < \gamma_k^{} \le \pi/2$.  In these ranges of the angles, for
$\theta\in\{2\vartheta_k^{},\gamma_k^{}\}$ the trigonometric
identities $\cos\theta = 1 / (1 + \tan^2\theta)$ and
$\sin\theta = \tan\theta \cos\theta$ hold when $\theta < \pi/2$.
Otherwise, $\tan\theta = \infty$, $\cos\theta = 0$, and
$\sin\theta = 1$.  Then, compute
$c_{2\vartheta}^{}\assgn\cos(2\vartheta_k^{})$,
$s_{2\vartheta}^{}\assgn\sin(2\vartheta_k^{})$,
$c_{\gamma}^{}\assgn\cos\gamma_k^{}$, and
$s_{\gamma}^{}\assgn\sin\gamma_k^{}$, and with them finally obtain
\begin{displaymath}
  \begin{aligned}
    \cos\varphi_k^{}&\assgn \frac{1}{\sqrt{2}}\sqrt{1+x_k^{}s_{2\vartheta}^{}+t_k^{}c_{\gamma}^{}c_{2\vartheta}^{}},\\
    \cos\psi_k^{}&\assgn \frac{1}{\sqrt{2}}\sqrt{1-x_k^{}s_{2\vartheta}^{}+t_k^{}c_{\gamma}^{}c_{2\vartheta}^{}},\\
    e^{i\alpha_k^{}}\sin\varphi_k^{}&\assgn e^{i\zeta_k^{}}\frac{(s_{2\vartheta}^{}-x_k^{})+it_k^{}s_{\gamma}^{}c_{2\vartheta}^{}}{2\cos\psi_k^{}},\\
    e^{-i\beta_k^{}}\sin\psi_k^{}&\assgn e^{-i\zeta_k^{}}\frac{(s_{2\vartheta}^{}+x_k^{})-it_k^{}s_{\gamma}^{}c_{2\vartheta}^{}}{2\cos\varphi_k^{}},
  \end{aligned}
\end{displaymath}
where $0 \leq \varphi_k^{} < \pi/2$ and $0 \leq \psi_k^{} < \pi/2$.

\paragraph{An exception}%
If $v_k^{}=h_k^{}=0$, i.e., if
$\arg(b_{i_k^{}j_k^{};k}') = \arg(a_{i_k^{}j_k^{};k}')$
and $a_{i_k^{}i_k^{};k}' = a_{j_k^{}j_k^{};k}'$, then
$\tan\gamma_k^{}$ is undefined, and $\tan(2\vartheta_k^{})$ might also
be.  In that case, it can be shown that $\widehat{A}_k'$ and
$\widehat{B}_k'$ are diagonalized by
\begin{displaymath}
  \widehat{Z}_k'\assgn\frac{1}{\sqrt{2}}
  \begin{bmatrix}
    \frac{1}{\sqrt{1 + x}} & \frac{-e^{i\zeta_k^{}}}{\sqrt{1 - x}}\\[6pt]
    \frac{e^{-i\zeta_k^{}}}{\sqrt{1 + x}} & \frac{1}{\sqrt{1 - x}}
  \end{bmatrix}.
\end{displaymath}
%
%
\subsubsection{The real case}\label{sss:2.1.2}
%
%
The transformation matrix $\widehat{Z}_k'$ is sought in a
form~\cite{Hari-84,Novakovic-Singer-Singer-2015}
\begin{displaymath}
  \widehat{Z}_k' \assgn \frac{1}{t_k^{}}
  \begin{bmatrix}
    \hphantom{-}\cos\varphi_k^{} & \sin\varphi_k^{}\\
    -\sin\psi_k^{} & \cos\psi_k^{}
  \end{bmatrix}.
\end{displaymath}
To that end, let $\displaystyle x_k^{} \assgn b_{i_k^{}j_k^{};k}'$ and
$\displaystyle t_k^{} \assgn \sqrt{1 - x_k^2} > 0$.  Then, set
\begin{displaymath}
  \begin{aligned}
  \xi_k^{} & \assgn \frac{x_k^{}}{\sqrt{1 + x_k^{}} + \sqrt{1 - x_k^{}}},\\
  \eta_k^{} & \assgn \frac{x_k^{}}{\left(1 + \sqrt{1 + x_k^{}}\right)\left(1 + \sqrt{1 - x_k^{}}\right)},
  \end{aligned}
\end{displaymath}
and compute
\begin{displaymath}
  \cot(2\vartheta_k^{})\assgn
  \frac{t_k^{}(a_{j_k^{}j_k^{};k}'-a_{i_k^{}i_k^{};k}')}{2a_{i_k^{}j_k^{};k}' - (a_{i_k^{}i_k^{};k}' + a_{j_k^{}j_k^{};k}')x_k^{}},
\end{displaymath}
where $-\pi/4 < \vartheta_k^{} \le \pi/4$.

Note that $\cot(2\vartheta_k^{})$ and $\cot\vartheta_k^{}$ (and the
corresponding tangents) have the same sign in the range of
$\vartheta_k^{}$.  Assuming that the floating-point arithmetic unit
does not trap on $\pm 1/0$ and $1/\infty$, obtain $\tan\vartheta_k^{}$
as
\begin{displaymath}
  \tan\vartheta_k^{} \assgn \frac{\sign(1, \cot(2\vartheta_k^{}))}{|\cot(2\vartheta_k^{})| + \sqrt{1 + \cot^2(2\vartheta_k^{})}},
\end{displaymath}
and from it $\cos\vartheta_k^{}$ and $\sin\vartheta_k^{}$ using the
same trigonometric identities  as in the complex case.  Finally,
compute
\begin{displaymath}
  \begin{aligned}
    \cos\varphi_k^{} & \assgn \cos\vartheta_k^{} +
    \xi_k^{}(\sin\vartheta_k^{} - \eta_k^{}\cos\vartheta_k^{}),\\
    \cos\psi_k^{} & \assgn \cos\vartheta_k^{} -
    \xi_k^{}(\sin\vartheta_k^{} + \eta_k^{}\cos\vartheta_k^{}),\\
    \sin\varphi_k^{} & \assgn \sin\vartheta_k^{} -
    \xi_k^{}(\cos\vartheta_k^{} + \eta_k^{}\sin\vartheta_k^{}),\\
    \sin\psi_k^{} & \assgn \sin\vartheta_k^{} +
    \xi_k^{}(\cos\vartheta_k^{} - \eta_k^{}\sin\vartheta_k^{}),
  \end{aligned}
\end{displaymath}
where $0\leq\varphi_k^{}<\pi/2$ and $0\leq\psi_k^{}<\pi/2$.

\paragraph{An exception}%
Since the real case is in fact a simplification of the complex case,
when $\cot(2\vartheta_k^{})$ is undefined, being $0/0$, i.e., when
$a_{i_k^{}i_k^{};k}' = a_{j_k^{}j_k^{};k}'$ and
$a_{i_k^{}j_k^{};k}' = a_{i_k^{}i_k^{};k}' b_{i_k^{}j_k^{};k}'$ (or,
in other words, when $\widehat{A}_k'$ and $\widehat{B}_k'$ are
proportional), define
\begin{displaymath}
  \widehat{Z}_k' \assgn
  \frac{1}{\sqrt{2}}
  \begin{bmatrix}
    \frac{1}{\sqrt{1 + |x|}} & \frac{-1}{\sqrt{1 - |x|}}\\[6pt]
    \frac{1}{\sqrt{1 + |x|}} & \frac{\hphantom{-}1}{\sqrt{1 - |x|}}
  \end{bmatrix}.
\end{displaymath}

\begin{figure*}[h!btp]
  \centering
  \includegraphics[keepaspectratio,width=\textwidth]{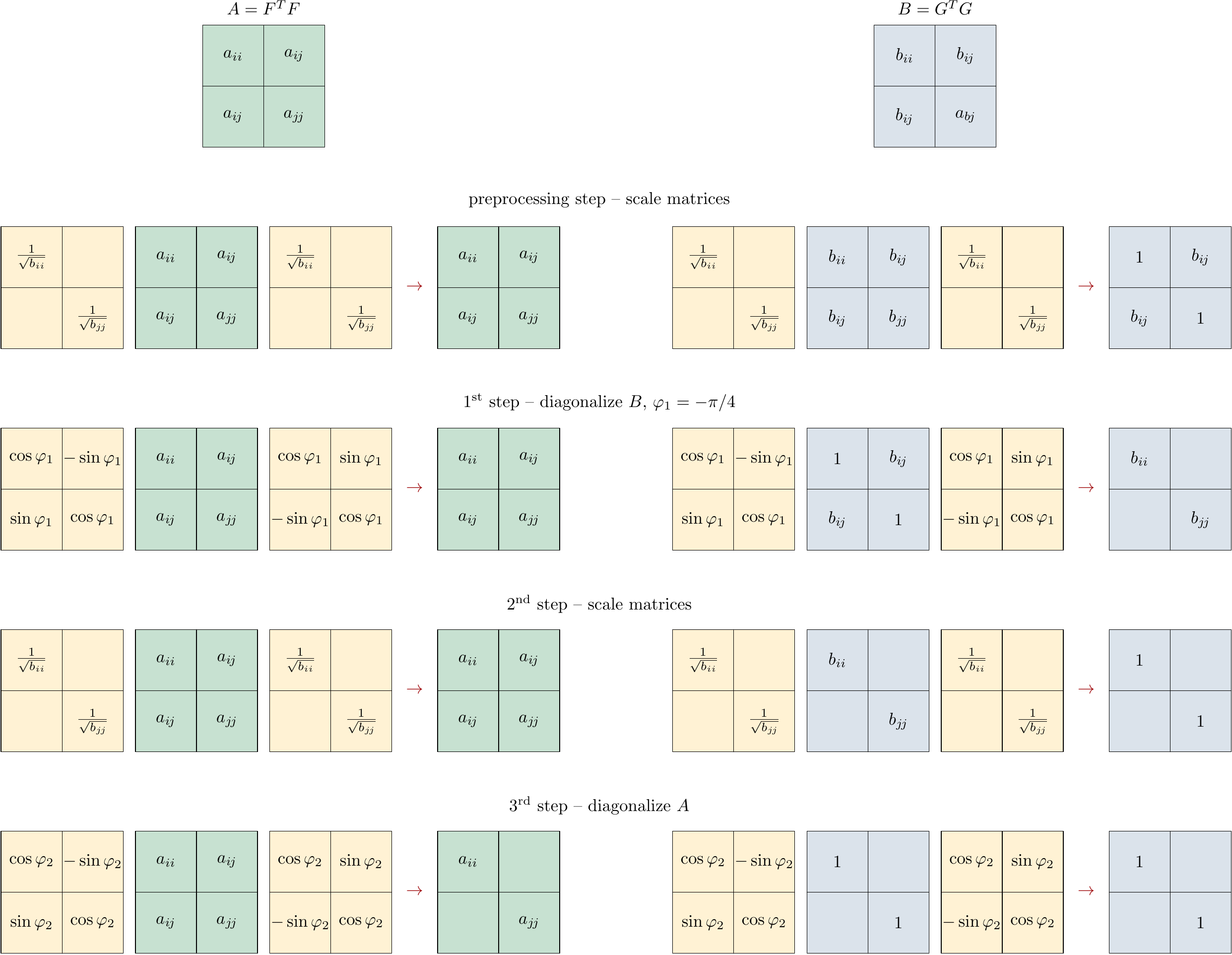}
  \caption{\looseness=-1
    A decomposition of the two-sided Hari--Zimmermann joint
    diagonalization of a pair $(A,B)$ of $2\times 2$ symmetric
    matrices into four simple transformations.  The second and the
    last are the orthogonal Jacobi rotations, but the first and the
    third (the diagonal scalings) are not orthogonal in general and
    also demonstrate why the positive definiteness of $B$ is essential
    for the method.}
  \label{fig:cartoon1}
\end{figure*}

\looseness=-1
Figure~\ref{fig:cartoon1} shows a schematic derivation of the
two-sided Hari--Zimmermann transformations.  Starting with a pair
$(A,B)$ of $2\times 2$ symmetric matrices, where $B$ is positive
definite, they are jointly transformed in four steps, i.e., twice by a
diagonal scaling followed by a Jacobi rotation, after which both
matrices become diagonal.  The above formulas for $\widehat{Z}_k'$
follow by combining the last three steps into a convenient computation
without the intermediate matrices.
%
%
\subsection{Parallelization of the one-sided algorithm}\label{ss:2.2}
%
%
The sequential one-sided algorithm in each step chooses a single pivot
index pair, according to some criterion that is called a sequential
Jacobi strategy.  However, at most $\lfloor n/2\rfloor$ pivot column
pairs of each matrix can be transformed concurrently if the indices in
all index pairs are distinct.

In a parallel step $k\ge 0$ a sequence $(i_k^{(\ell)},j_k^{(\ell)})$
of pivot index pairs, where $1\le\ell\le\lfloor n/2\rfloor$, such that
each index in the range from 1 to $n$ appears at most (and for even
$n$, exactly) once in the sequence, addresses $\lfloor n/2\rfloor$
pivot column pairs of $A_k^{}$ and $B_k^{}$ to be transformed---each
pair by a separate, concurrent task.  All permutations of a given
$(i_k^{(\ell)},j_k^{(\ell)})$ are equivalent from the numerical point
of view, since the resulting $A_{k+1}^{}$ and $B_{k+1}^{}$ are the
same for every reordering of the sequence, and therefore any
reordering represents the entire equivalence class.

For simplicity, a barrier is assumed between the successive parallel
steps, i.e., all tasks of a step have to be completed before those of
the following step are started.

A criterion to choose a pivot index pair sequence for each parallel
step is called a parallel Jacobi strategy.  Among the strategies that
are simplest to compute are the ones that prescribe a pivot sequence
for each step, until all $n(n-1)/2$ index pairs $(i,j)$ are selected
at least once.  The choice of the steps is then periodically repeated.
Let $s$ be the shortest such period.  The first $s$ steps constitute
the first sweep, the following $s$ steps the second sweep, and so on.

If in any sweep exactly $n(n-1)/2$ different index pairs are chosen,
such a strategy is called cyclic; otherwise, some index pairs are
repeated in a sweep, and the strategy is called quasi-cyclic.  For
even $n$, $s \ge n-1$, and the equality holds if and only if the
strategy is cyclic.

A strategy is defined for a fixed $n$; however, by a slight abuse of
the usual terminology, a single principle by which the particular
strategies are generated for some given matrix orders will simply be
called a strategy kind, or even a strategy for short.

Based on the previous experience with the one-sided Jacobi-like
algorithms, two parallel Jacobi strategy kinds have been selected for
testing: the modified modulus (\textsc{mm}; see,
e.g.,~\cite{Novakovic-SingerSanja-2011,Novakovic-Singer-Singer-2015}),
quasi-cyclic with $s=n$, and the generalized Mantharam--Eberlein
(\textsc{me}; see~\cite{Mantharam-Eberlein-93,Novakovic-2015}) cyclic
one.  Please see Figures~1 and 2 in the supplementary material, where
a sweep of \textsc{me} and of \textsc{mm}, respectively, is shown
two-sidedly on a matrix of order 32.
%
%
\subsection{Blocking of the one-sided algorithm}\label{ss:2.3}
%
%
Parallelization alone is not sufficient for achieving a decent
performance of the algorithm on the modern architectures with multiple
levels of the memory hierarchy.

The pointwise algorithm just described is therefore modified to work
on the block columns of the matrices, instead of the columns proper.
Each block column comprises an arbitrary but fixed number \texttt{w},
$1 < \mathtt{w} < \lfloor n/2 \rfloor$, of consecutive matrix
columns.  Instead of $2\times 2$ pivot submatrices of $A_k^{}$ and
$B_k^{}$, in the blocked algorithm $2\mathtt{w}\times 2\mathtt{w}$
pivot submatrices $\widehat{\mathsf{A}}_k^{(\ell)}$ and
$\widehat{\mathsf{B}}_k^{(\ell)}$ are formed in the $k$th (parallel or
sequential) step by matrix multiplications,
\begin{displaymath}
  \begin{aligned}
    \widehat{\mathsf{A}}_k^{(\ell)}\assgn&
    \begin{bmatrix}
      A_{i_k^{(\ell)}i_k^{(\ell)};k}^{} & A_{i_k^{(\ell)}j_k^{(\ell)};k}^{}\\
      A_{i_k^{(\ell)}j_k^{(\ell)};k}^{\ast} & A_{j_k^{(\ell)}j_k^{(\ell)};k}^{}
    \end{bmatrix}\\=&
    \begin{bmatrix}
      F_{i_k^{(\ell)};k}^{\ast} F_{i_k^{(\ell)};k}^{} & F_{i_k^{(\ell)};k}^{\ast} F_{j_k^{(\ell)};k}^{}\\
      F_{j_k^{(\ell)};k}^{\ast} F_{i_k^{(\ell)};k}^{} & F_{j_k^{(\ell)};k}^{\ast} F_{j_k^{(\ell)};k}^{}
    \end{bmatrix},\\
    \widehat{\mathsf{B}}_k^{(\ell)}\assgn&
    \begin{bmatrix}
      B_{i_k^{(\ell)}i_k^{(\ell)};k}^{} & B_{i_k^{(\ell)}j_k^{(\ell)};k}^{}\\
      B_{i_k^{(\ell)}j_k^{(\ell)};k}^{\ast} & B_{j_k^{(\ell)}j_k^{(\ell)};k}^{}
    \end{bmatrix}\\=&
    \begin{bmatrix}
      G_{i_k^{(\ell)};k}^{\ast} G_{i_k^{(\ell)};k}^{} & G_{i_k^{(\ell)};k}^{\ast} G_{j_k^{(\ell)};k}^{}\\
      G_{j_k^{(\ell)};k}^{\ast} G_{i_k^{(\ell)};k}^{} & G_{j_k^{(\ell)};k}^{\ast} G_{j_k^{(\ell)};k}^{}
    \end{bmatrix},
  \end{aligned}
\end{displaymath}
where $F_{i_k^{(\ell)};k}^{}$, $F_{j_k^{(\ell)};k}^{}$,
$G_{i_k^{(\ell)};k}^{}$, $G_{j_k^{(\ell)};k}^{}$,
$Z_{i_k^{(\ell)};k}^{}$, and $Z_{j_k^{(\ell)};k}^{}$ are the
$i_k^{(\ell)}$th and $j_k^{(\ell)}$th block columns of $F_k^{}$,
$G_k^{}$, and $Z_k^{}$ of width \texttt{w}.

Now, $\widehat{\mathsf{A}}_k^{(\ell)}$ and
$\widehat{\mathsf{B}}_k^{(\ell)}$ can either be jointly diagonalized
by a matrix $\widehat{\mathsf{Z}}_k^{(\ell)}$, which leads to the full
block (\textsc{fb}) algorithm~\cite{Hari-SingerSanja-SingerSasa-2014},
as called in the context of the Jacobi methods, or their off-diagonal
norms can be reduced by a sequence of congruences accumulated into
$\widehat{\mathsf{Z}}_k^{(\ell)}$, which is called the block-oriented
(\textsc{bo}) algorithm~\cite{Hari-SingerSanja-SingerSasa-2010}.  The
idea behind blocking is that $\widehat{\mathsf{A}}_k^{(\ell)}$,
$\widehat{\mathsf{B}}_k^{(\ell)}$, and
$\widehat{\mathsf{Z}}_k^{(\ell)}$ fit, by choosing \texttt{w}, into
the small but fast cache memory (e.g., the shared memory of a GPU), to
speed up the computation with them, as well as employing BLAS~3
(matrix multiplies) operations for the block column updates by
$\widehat{\mathsf{Z}}_k^{(\ell)}$ afterwards:
\begin{displaymath}
  \begin{aligned}
    \begin{bmatrix}
      F_{i_k^{(\ell)};k+1}^{} & F_{j_k^{(\ell)};k+1}^{}
    \end{bmatrix} & =
    \begin{bmatrix}
      F_{i_k^{(\ell)};k}^{} & F_{j_k^{(\ell)};k}^{}
    \end{bmatrix}
    \cdot \widehat{\mathsf{Z}}_k^{(\ell)},\\
    \begin{bmatrix}
      G_{i_k^{(\ell)};k+1}^{} & G_{j_k^{(\ell)};k+1}^{}
    \end{bmatrix} & =
    \begin{bmatrix}
      G_{i_k^{(\ell)};k}^{} & G_{j_k^{(\ell)};k}^{}
    \end{bmatrix}
    \cdot \widehat{\mathsf{Z}}_k^{(\ell)},\\
    \begin{bmatrix}
      Z_{i_k^{(\ell)};k+1}^{} & Z_{j_k^{(\ell)};k+1}^{}
    \end{bmatrix} & =
    \begin{bmatrix}
      Z_{i_k^{(\ell)};k}^{} & Z_{j_k^{(\ell)};k}^{}
    \end{bmatrix}
    \cdot \widehat{\mathsf{Z}}_k^{(\ell)}.
  \end{aligned}
\end{displaymath}

The computation of $\widehat{\mathsf{Z}}_k^{(\ell)}$ in either
\textsc{fb} or \textsc{bo} can be done by any convergent method; a
two-sided method can be applied straightforwardly, but for the
one-sided approach $\widehat{\mathsf{A}}_k^{(\ell)}$ and
$\widehat{\mathsf{B}}_k^{(\ell)}$ have to be factorized first by,
e.g., the Cholesky factorization
\begin{displaymath}
  \widehat{\mathsf{A}}_k^{(\ell)} =
  \widehat{\mathsf{F}}_k^{(\ell)\ast}\widehat{\mathsf{F}}_k^{(\ell)},\quad
  \widehat{\mathsf{B}}_k^{(\ell)} =
  \widehat{\mathsf{G}}_k^{(\ell)\ast}\widehat{\mathsf{G}}_k^{(\ell)},
\end{displaymath}
and then the same implicit
Hari--Zimmermann method, pointwise or blocked, and in both cases,
either parallel or sequential, can be recursively applied to
$\widehat{\mathsf{F}}_k^{(\ell)}$ and
$\widehat{\mathsf{G}}_k^{(\ell)}$.

In the single-GPU algorithm, there is only one level of such a
recursion, i.e., one level of blocking.  The block, outer level of the
algorithm and the pointwise, inner level do not need to employ the same
strategy kind.  Both levels, however, are parallel.  The sweeps of the
outer level are called the block (or outer) sweeps, and those of the
inner level are called the pointwise (or inner) sweeps, which for
\textsc{fb} are limited to 30 ($\widehat{\mathsf{A}}_k^{(\ell)}$ and
$\widehat{\mathsf{B}}_k^{(\ell)}$ are usually fully diagonalized in
less than that number of sweeps), and for \textsc{bo} are limited to
only one inner sweep.   Apart from that, there is no other substantial
difference between \textsc{fb} and \textsc{bo}.

\looseness=-1
The Cholesky factorization is not the only way to form
$\widehat{\mathsf{F}}_k^{(\ell)}$ and
$\widehat{\mathsf{G}}_k^{(\ell)}$.  One numerical stability
improvement would be to use a diagonally pivoted version of the
factorization
instead~\cite{SingerSanja-SingerSasa-Novakovic-Uscumlic-Dunjko-2012},
\begin{displaymath}
  \widehat{\mathsf{A}}_k^{(\ell)} =
  P_{F;k}^{(\ell)} \widetilde{\mathsf{F}}_k^{(\ell)\ast} \widetilde{\mathsf{F}}_k^{(\ell)} P_{F;k}^{(\ell)T},\ \ \
  \widehat{\mathsf{B}}_k^{(\ell)} =
  P_{G;k}^{(\ell)} \widetilde{\mathsf{G}}_k^{(\ell)\ast} \widetilde{\mathsf{G}}_k^{(\ell)} P_{G;k}^{(\ell)T}.
\end{displaymath}
Another one would be to skip forming $\widehat{\mathsf{A}}_k^{(\ell)}$
and $\widehat{\mathsf{B}}_k^{(\ell)}$ explicitly by shortening the
pivot block columns by the column-pivoted QR factorization
directly~\cite{Singer-DiNapoli-Novakovic-Caklovic-2020},
\begin{displaymath}
  \begin{aligned}
    \widetilde{\mathsf{F}}_k^{(\ell)}&\assgn
    \begin{bmatrix}
      \widetilde{F}_{i_k^{(\ell)};k}^{} & \widetilde{F}_{j_k^{(\ell)};k}^{}
    \end{bmatrix}\\&=
    Q_{F;k}^{(\ell)\ast}\cdot
    \begin{bmatrix}
      F_{i_k^{(\ell)};k}^{} & F_{j_k^{(\ell)};k}^{}
    \end{bmatrix}
    \cdot P_{F;k}^{(\ell)},\\
    \widetilde{\mathsf{G}}_k^{(\ell)}&\assgn
    \begin{bmatrix}
      \widetilde{G}_{i_k^{(\ell)};k+1}^{} & \widetilde{G}_{j_k^{(\ell)};k}^{}
    \end{bmatrix}\\&=
    Q_{G;k}^{(\ell)\ast} \cdot
    \begin{bmatrix}
      G_{i_k^{(\ell)};k}^{} & G_{j_k^{(\ell)};k}^{}
    \end{bmatrix}
    \cdot P_{G;k}^{(\ell)}.
  \end{aligned}
\end{displaymath}
In both cases, let
\begin{displaymath}
  \widehat{\mathsf{F}}_k^{(\ell)} \assgn \widetilde{\mathsf{F}}_k^{(\ell)} P_{F;k}^{(\ell)T},\quad
  \widehat{\mathsf{G}}_k^{(\ell)} \assgn \widetilde{\mathsf{G}}_k^{(\ell)} P_{G;k}^{(\ell)T},
\end{displaymath}
where $P_{F;k}^{(\ell)}$ and $P_{G;k}^{(\ell)}$ are permutation
matrices, while $Q_{F;k}^{(\ell)}$ and $Q_{G;k}^{(\ell)}$ are unitary
and are not required to be stored, implicitly or explicitly, for any
further computation.

\looseness=-1
However, the QR factorization (even without the column pivoting) of
a pair of the tall and skinny block columns comes with a significant
performance penalty on a GPU compared to the Cholesky factorization of
a small, square pivot submatrix~\cite{Novakovic-2015}, and the pivoted
Cholesky factorization does not avoid a possibility of getting a
severely ill-conditioned $\widehat{\mathsf{A}}_k^{(\ell)}$ or
$\widehat{\mathsf{B}}_k^{(\ell)}$ by multiplying an ill-conditioned
pair of block columns by itself.  Both of these enhancements are
therefore only mentioned here, with a performance comparison of the
in-kernel QR factorizations versus the formation of the Grammian
matrices and their Cholesky factorizations available in
Appendix~\ref{s:E}.  If the batched tall-and-skinny QR factorizations
prove indispensable for a particularly ill-conditioned problem,
\texttt{cublasXgetrfBatched} routine (with
$\text{\texttt{X}}\in\{\text{\texttt{D}},\text{\texttt{Z}}\}$) and
\cite{Boukaram-et-al-2018} could also be considered.

\begin{figure*}[h!btp]
  \centering
  \includegraphics[keepaspectratio,width=\textwidth]{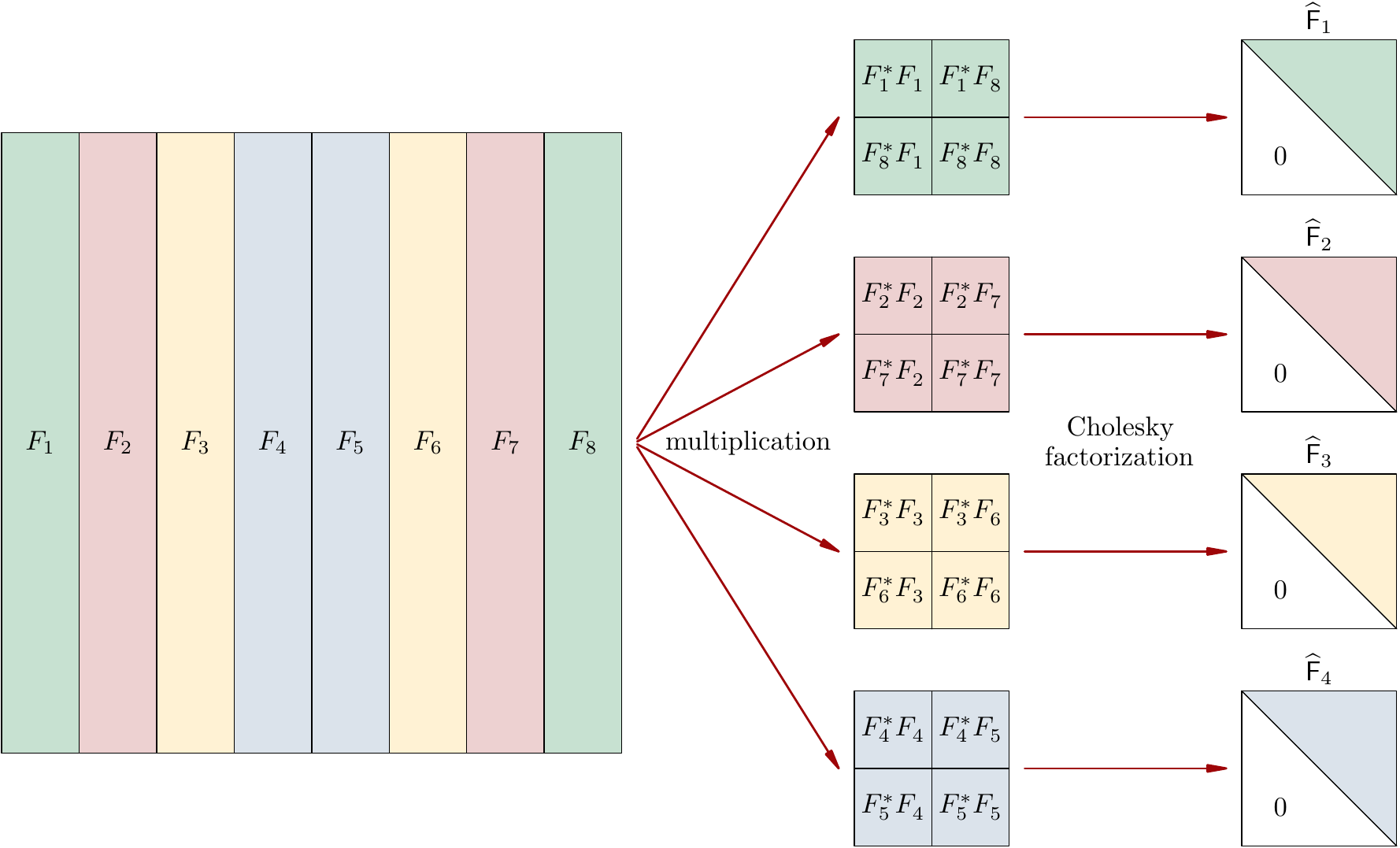}
  \caption{\looseness=-1
    Formation of the Grammian matrices from the pairs of block columns
    of $F$ and their subsequent Cholesky factorizations.  Each pair is
    indicated by a different hue, and varies with a block step.  The
    same process is repeated for $G$, to obtain the factors
    $\widehat{\mathsf{G}}_k^{}$.}
  \label{fig:cartoon2}
\end{figure*}

In the following, the blocked algorithm is assumed to form the pivot
submatrices as
\begin{displaymath}
  \begin{aligned}
    \widehat{\mathsf{A}}_k^{(\ell)} & \assgn
    \begin{bmatrix}
      F_{i_k^{(\ell)};k}^{} & F_{j_k^{(\ell)};k}^{}
    \end{bmatrix}^{\ast}\cdot
    \begin{bmatrix}
      F_{i_k^{(\ell)};k}^{} & F_{j_k^{(\ell)};k}^{}
    \end{bmatrix},\\
    \widehat{\mathsf{B}}_k^{(\ell)} & \assgn
    \begin{bmatrix}
      G_{i_k^{(\ell)};k}^{} & G_{j_k^{(\ell)};k}^{}
    \end{bmatrix}^{\ast}\cdot
    \begin{bmatrix}
      G_{i_k^{(\ell)};k}^{} & G_{j_k^{(\ell)};k}^{}
    \end{bmatrix},
  \end{aligned}
\end{displaymath}
i.e., each one by a \texttt{ZHERK} (\texttt{DSYRK} in the real case)
like operation in the BLAS terminology, and the non-pivoted Cholesky
factorization is then used to obtain $\widehat{\mathsf{F}}_k^{(\ell)}$
and $\widehat{\mathsf{G}}_k^{(\ell)}$, as demonstrated in
Figure~\ref{fig:cartoon2}, where eight block columns of $F$ are
depicted.  The same illustration holds if $F$ is replaced by $G$. The
block columns of the same hue are paired together, according to the
first step of the \textsc{me} strategy, giving four square blocks to
be formed and factorized.

\begin{figure*}[h!btp]
  \centering
  \includegraphics[keepaspectratio,width=\textwidth]{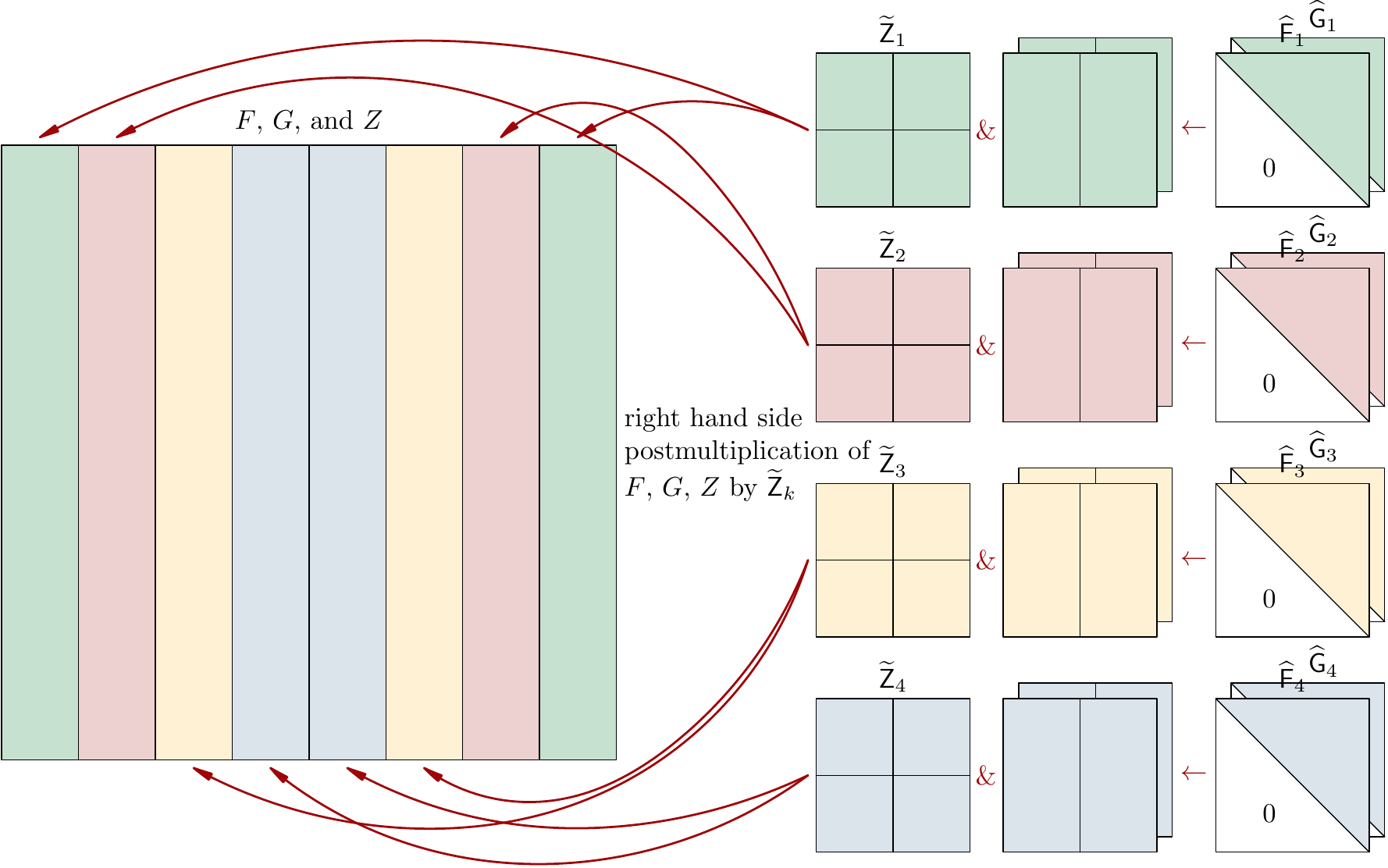}
  \caption{The pointwise Hari--Zimmermann GSVD of four
    $(\widehat{\mathsf{F}}_k^{},\widehat{\mathsf{G}}_k^{})$ pairs
    results in two unused scaled left generalized singular vector
    matrices per pair, and a single accumulated and rescaled
    transformation matrix $\widetilde{\mathsf{Z}}_k^{}$.  Each of the
    four original block column pairs of $F$, $G$, and $Z$ is then
    updated by multiplying it from the right by the corresponding
    $\widetilde{\mathsf{Z}}_k^{}$ (indicated by a pair of arrows).}
  \label{fig:cartoon3}
\end{figure*}

\looseness=-1
Figure~\ref{fig:cartoon3} shows how each pair of the factors
$\widehat{\mathsf{F}}_k^{}$ and $\widehat{\mathsf{G}}_k^{}$ is
processed by the pointwise Hari--Zimmermann GSVD, leaving two matrices
of the scaled left generalized singular vectors that are not used
further, and a single matrix (rescaled, as noted in the following
subsection~\ref{ss:2.4}) $\widetilde{\mathsf{Z}}_k^{}$ of the
accumulated transformations.  The block column pairs of $F$, $G$, and
$Z$, with the physically disjoint but logically contiguous block
columns, are then postmultiplied, each from the right by the
corresponding $\widetilde{\mathsf{Z}}_k^{}$, and replaced by the
result.
%
%
\subsection{Rescalings}\label{ss:2.4}
%
%
Observe that $\widehat{\mathsf{Z}}_k^{(\ell)}$ is a product of several
non-unitary matrices, elements of which can be larger than 1 by
magnitude, so the norm of $\widehat{\mathsf{Z}}_k^{(\ell)}$ can build
up significantly by such accumulation of the transformations.  Also,
if $\widehat{\mathsf{Z}}_k^{(\ell)}$ diagonalizes
$\widehat{\mathsf{A}}_k^{(\ell)}$ and
$\widehat{\mathsf{B}}_k^{(\ell)}$, or reduces their off-diagonal
norms, so does any matrix
$\widehat{\mathsf{Z}}_k^{(\ell)}\widehat{\mathsf{\Theta}}_k^{(\ell)}$,
where $\widehat{\mathsf{\Theta}}_k^{(\ell)}$ is a real, diagonal
matrix with its diagonal elements positive and smaller than 1.

Let $\widetilde{\mathsf{\Theta}}_k^{(\ell)}$ be such a matrix that
reduces the norm of $\widehat{\mathsf{Z}}_k^{(\ell)}$,
\begin{displaymath}
  \left(\widetilde{\mathsf{\Theta}}_k^{(\ell)}\right)_{\null\!jj}^{}\assgn
  \left(\left\|\left(\widehat{\mathsf{F}}_k^{(\ell)}\right)_{\null\!j}'\right\|_F^2+\left\|\left(\widehat{\mathsf{G}}_k^{(\ell)}\right)_{\null\!j}'\right\|_F^2\right)^{-1/2},
\end{displaymath}
where $\left(\widehat{\mathsf{F}}_k^{(\ell)}\right)_{\null\!j}'$ and
$\left(\widehat{\mathsf{G}}_k^{(\ell)}\right)_{\null\!j}'$ stand for
the $j$th column of the final, transformed
$\widehat{\mathsf{F}}_k^{(\ell)}$ and
$\widehat{\mathsf{G}}_k^{(\ell)}$, respectively, of which the latter
has unit norm, and thus
$\max_{\null\!j}^{}\left(\widetilde{\mathsf{\Theta}}_k^{(\ell)}\right)_{\null\!jj}^{} < 1$.

This is exactly the same scaling as it would be performed in the last,
post-iterative phase of the algorithm,
\begin{displaymath}
  (\widetilde{\Sigma}_F^{})_{\null\!jj}^{}\assgn \|(F_N^{})_{\null\!j}^{}\|_F^{},\quad
  (\widetilde{\Sigma}_G^{})_{\null\!jj}^{}\assgn \|(G_N^{})_{\null\!j}^{}\|_F^{},
\end{displaymath}
\begin{displaymath}
  \widetilde{\Theta}_{\null\!jj}^{}\assgn ((\widetilde{\Sigma}_F^{})_{\null\!jj}^2 + (\widetilde{\Sigma}_G^{})_{\null\!jj}^2)^{-1/2},
\end{displaymath}
\begin{displaymath}
  \begin{aligned}
    \Sigma_F &\assgn \widetilde{\Sigma}_F^{} \widetilde{\Theta},\\
    U &\assgn F_N \widetilde{\Sigma}_F^{-1},
  \end{aligned}\quad
  \begin{aligned}
    \Sigma_G &\assgn \widetilde{\Sigma}_G^{} \widetilde{\Theta},\\
    V &\assgn G_N \widetilde{\Sigma}_G^{-1},
  \end{aligned}\quad
  \begin{aligned}
    \Sigma &\assgn \Sigma_G^{-1} \Sigma_F^{},\\
    Z &\assgn Z_N \widetilde{\Theta},
  \end{aligned}
\end{displaymath}
except that $\widehat{\mathsf{F}}_k^{(\ell)\prime}$ and
$\widehat{\mathsf{G}}_k^{(\ell)\prime}$ do not have to be rescaled and
the norms of their columns do not have to be kept as they are all
discarded immediately after $\widetilde{\mathsf{\Theta}}_k^{(\ell)}$
has been computed.

Then,
$\widetilde{\mathsf{Z}}_k^{(\ell)} \assgn \widehat{\mathsf{Z}}_k^{(\ell)} \widetilde{\mathsf{\Theta}}_k^{(\ell)}$,
is applied to the pivot block column pair of $F_k^{}$, $G_k^{}$, and
$Z_k^{}$ instead of $\widehat{\mathsf{Z}}_k^{(\ell)}$, and is
considered embedded into $\widetilde{Z}_k^{}$ in a similar way as
$\widehat{Z}_k^{(\ell)}$ would be in the pointwise case, i.e.,
starting from $\widetilde{Z}_k^{}$ being $I_n^{}$, for each $\ell$ let
\begin{displaymath}
  \begin{aligned}
    \widetilde{Z}_k^{}&((i_k^{(\ell)}-1)\cdot\mathtt{w}+1:i_k^{(\ell)}\cdot\mathtt{w},\,(i_k^{(\ell)}-1)\cdot\mathtt{w}+1:i_k^{(\ell)}\cdot\mathtt{w})\\
                 \null&\assgn\widetilde{\mathsf{Z}}_k^{(\ell)}(1:\mathtt{w},\,1:\mathtt{w})\\
    \widetilde{Z}_k^{}&((i_k^{(\ell)}-1)\cdot\mathtt{w}+1:i_k^{(\ell)}\cdot\mathtt{w},\,(j_k^{(\ell)}-1)\cdot\mathtt{w}+1:j_k^{(\ell)}\cdot\mathtt{w})\\
                 \null&\assgn\widetilde{\mathsf{Z}}_k^{(\ell)}(1:\mathtt{w},\,\mathtt{w}+1:2\cdot\mathtt{w})\\
    \widetilde{Z}_k^{}&((j_k^{(\ell)}-1)\cdot\mathtt{w}+1:j_k^{(\ell)}\cdot\mathtt{w},\,(i_k^{(\ell)}-1)\cdot\mathtt{w}+1:i_k^{(\ell)}\cdot\mathtt{w})\\
                 \null&\assgn\widetilde{\mathsf{Z}}_k^{(\ell)}(\mathtt{w}+1:2\cdot\mathtt{w},\,1:\mathtt{w})\\
    \widetilde{Z}_k^{}&((j_k^{(\ell)}-1)\cdot\mathtt{w}+1:j_k^{(\ell)}\cdot\mathtt{w},\,(j_k^{(\ell)}-1)\cdot\mathtt{w}+1:j_k^{(\ell)}\cdot\mathtt{w})\\
                 \null&\assgn\widetilde{\mathsf{Z}}_k^{(\ell)}(\mathtt{w}+1:2\cdot\mathtt{w},\,\mathtt{w}+1:2\cdot\mathtt{w}),
  \end{aligned}
\end{displaymath}
where the subscripting is to be interpreted as in Fortran.

To reduce the norm of the entire $Z_k^{}$, a similar rescaling can be
applied on $Z_k^{}$, using the column norms of $F_k^{}$ and $G_k^{}$,
after each but the last block sweep.  After the last block sweep, a
rescaling of all three matrices ($F_N^{}$, $G_N^{}$, and $Z_N^{}$) is
performed to obtain $U$, $V$, and $Z$, with the extraction of
$\Sigma_F^{}$, $\Sigma_G^{}$, and $\Sigma$.
%
%
\subsection{Convergence}\label{ss:2.5}
%
%
The inner level of the algorithm stops when there were no
transformations, apart from the sorting permutations, applied in a
sweep, or when the prescribed maximal number of sweeps has been
reached.  Then, the pivot block column pairs of $F_k^{}$, $G_k^{}$,
and $Z_k^{}$ are updated concurrently for all $\ell$ by
$\widetilde{\mathsf{Z}}_k^{(\ell)}$, which can be skipped for those
$\ell$ where $\widetilde{\mathsf{Z}}_k^{(\ell)} = I_{2\mathtt{w}}$.

The same criterion could be used for the outer level, where the count
of transformations applied in an outer sweep is a sum of all
transformations applied in the inner level in all steps of the outer
sweep.  However, this criterion has to be
relaxed~\cite{Novakovic-2015,Novakovic-Singer-Singer-2015}, since the
rounding errors in forming and factorizing the block pivot submatrices
could spoil the attained numerical orthogonality of the original
columns, and introduce a small number of unwarranted transformations
that prevent the algorithm from detecting convergence even if it has
in fact been reached.

Therefore, the transformations are divided in two classes: ``big''
and ``small''.  The latter are all $\widehat{Z}_k^{(\ell)}$ where
either:
\begin{enumerate}
\item[C1.] $(\cos\varphi_k^{})/t_k^{}=(\cos\psi_k^{})/t_k^{}=1$, or
\item[C2.] $\cos\varphi_k^{}=\cos\psi_k^{}=1$,
\end{enumerate}
i.e., where $\widehat{Z}_k^{(\ell)}$ is close to (a multiple of)
identity, and the former are all other transformations.  Note that
neither definition of the small transformations implies that
$\sin\varphi_k^{}$ or $\sin\psi_k^{}$ are \emph{numerically\/} equal
to zero (and are usually not).  Also, since $x_k^{}$ tends to zero and
therefore $t_k^{}$ to one in the last sweeps of the algorithm, the
first and the second definition should not differ significantly.

There are separate counters of the big transformations, and of all
transformations applied in the inner level of the algorithm.  The
inner level still halts when there were no transformations of any
class in a sweep, but the outer level stops when there were no big
transformations applied in an outer sweeps (i.e., small
transformations are allowed to occur but do not spoil the overall
convergence).  Such a heuristic criterion prevents in practice a long
sequence of outer sweeps at the end of the algorithm, with only a few
transformations close to identity in each.
%
%
\subsection{Variants of the algorithm}\label{ss:2.6}
%
%
To summarize the variants of the algorithm, see Table~\ref{tbl:var}.
\begin{table}[h!bt]
  \caption{Variants of the implicit Hari--Zimmermann algorithm.}
  \label{tbl:var}
  \centering
  \addtolength{\tabcolsep}{-1pt}
  \begin{tabular}{@{}cccc@{}}
    \toprule
    ID & convergence & transformations & dot-products\\
    \midrule
    \textbf{0} & criterion C1 & $\widehat{Z}_k^{}$ ($F_k^{},G_k^{}$ prescaled) & ordinary\\
    1 & criterion C1 & $\widehat{Z}_k^{}$ ($F_k^{},G_k^{}$ prescaled) & enhanced\\
    2 & criterion C1 & $\widehat{Z}_k'$ ($F_k^{},G_k^{}$ not scaled) & ordinary\\
    3 & criterion C1 & $\widehat{Z}_k'$ ($F_k^{},G_k^{}$ not scaled) & enhanced\\
    4 & criterion C2 & $\widehat{Z}_k^{}$ ($F_k^{},G_k^{}$ prescaled) & ordinary\\
    5 & criterion C2 & $\widehat{Z}_k^{}$ ($F_k^{},G_k^{}$ prescaled) & enhanced\\
    6 & criterion C2 & $\widehat{Z}_k'$ ($F_k^{},G_k^{}$ not scaled) & ordinary\\
    7 & criterion C2 & $\widehat{Z}_k'$ ($F_k^{},G_k^{}$ not scaled) & enhanced\\
    \bottomrule
  \end{tabular}
\end{table}
The first column, ID, sets a shorthand for the corresponding variant.
The second column specifies a convergence criterion used.  The third
column distinguished between assuming the column norms of the second
matrix to be unity, and rescaling of both matrices with each
transformation.  The fourth column relates to computing the
dot-products the usual way, or by an enhanced, possibly more accurate
procedure from Appendix~\ref{s:A}.  Unless specified otherwise, the
column sorting is employed in all cases.  Thus, e.g.,
DHZ3-(\textsc{mm}-\textsc{fb}-\textsc{me}) refers to the
double-precision real implicit blocked Hari--Zimmermann algorithm with
ID equal to 3, using \textsc{mm} at the outer and \textsc{me} at the
inner level of blocking of type \textsc{fb}.  Similarly,
ZHZ0-(\textsc{me}-\textsc{bo}-\textsc{me}) stands for the
double-precision complex implicit blocked Hari--Zimmerman algorithm
with ID equal to 0 (the ``standard'' variant), using \textsc{me} at
both levels of blocking of type \textsc{bo}.

From now on, when a numeric variant ID is mentioned in the text, it is
assumed that it should be looked up in Table~\ref{tbl:var}.
%
%
\section{The single-GPU implementation}\label{s:3}
%
%
In this section the single-GPU implementation of the complex and the
real one-sided Hari--Zimmermann algorithms are described.  The focus
is on the complex algorithm, and the real one is commented on when
substantially different.  The target framework is CUDA
C++~\cite{NVidia-2019} for the NVIDIA GPUs (Kepler series and newer),
but also another general-purpose GPU programming environment with the
analogous concepts and constructs could probably be used.
%
%
\subsection{Data layout and transfer}\label{ss:3.1}
%
%
Due to blocking employed by the algorithm, each matrix is viewed as
column-striped, with the block columns containing $\mathtt{w}=16$
consecutive columns each.  To simplify the implementation, assume that
$n$ is a multiple of 32, and let $\mathtt{n}\assgn n/\mathtt{w}$ (so
$\mathtt{n}$ is even).  If the assumption does not hold for the input,
the matrices should then be bordered by appending $32-(n\bmod 32)$
columns to the right, and as many rows to the bottom.  The elements
$(m_Y^{}+1, n+1)$, $(m_Y^{}+2, n+2)$, \ldots, in the columns newly
added to the matrix $Y \in \{F, G\}$ should be set to unity, to avoid
introducing zero columns, since it is essential for $Y$ to be of full
column rank.  Other bordering elements should be set to zero, to
prevent any transformation not implied by the original matrices from
happening (see a bordering example
in~\cite{Novakovic-SingerSanja-2011}).

Another assumption, to simplify the loop unrolling in various parts of
the code, is to have $m_F^{}$ and $m_G^{}$ as a multiple of 64.  If it
is not the case with the input, then, after a possible bordering as
described above, $64-(m_Y^{} \bmod 64)$ rows of zeros should be appended
to the bottom of the matrix $Y \in \{F, G\}$.
%
%
\subsubsection{The CPU and the GPU RAM layout and transfer}\label{sss:3.1.1}
%
%
Data is laid out in the GPU RAM (also called ``global memory'' in the
GPU context) in the following sequence:
\begin{displaymath}
  \Real(F), \Imag(F);\quad
  \Real(G), \Imag(G);\quad
  \Real(Z), \Imag(Z),
\end{displaymath}
after which follow the output-only vectors $\Sigma$, $\Sigma_F$,
$\Sigma_G$ (with double-precision elements, each of length $n$), and
$C$ (holding unsigned 8-byte integers, of length $\mathtt{n}$).  The
rest of data is used both for input and output, i.e., the six
double-precision matrices are constantly being read and overwritten
within the GPU as the algorithm progresses.  The matrices are loaded
to the GPU at the beginning of the algorithm's execution, if they are
not already in place as a result of another computation, and
optionally copied to the CPU at its end, as well as $\Sigma$,
$\Sigma_F^{}$, and $\Sigma_G^{}$.

In the pre- and post-processing stages on the CPU, input ($F$, $G$)
and output data ($U$ in place of $F$; $V$ in place of $G$; and $Z$),
respectively, is repacked from (or to) the standard representation of
complex matrices, in which the successive elements are complex numbers
$z=(\Real(z),\Imag(z))$.  Each double-precision matrix can therefore
be loaded to, or copied from, the GPU with a single CUDA call.

This decision to keep all data in real-typed variables by splitting
the real and the imaginary matrix parts and to perform the complex
arithmetic manually is a design choice, not a necessity, since an
implementation of the algorithm with the real and the imaginary parts
interleaved in the customary way is also possible.  There is no direct
support for the standard C (with \texttt{\_Complex} types) or C++
(with \texttt{std::complex} types) complex arithmetic in CUDA, so some
non-standard approach has to be used anyway; e.g., the datatypes and
the routines from the \texttt{cuComplex.h} header file, or those from
the \texttt{thrust} library, or a custom implementation---possibly
with a different memory layout---of complex numbers and the operations
with them.  The chosen, custom approach with the split data layout
makes reading or writing only one (real or imaginary) component of the
successive matrix elements straightforward, and such memory accesses
can be contiguous.

In the auxiliary vector $C$ there are two counters, $C_{\ell}^{(0)}$
and $C_{\ell}^{(1)}$, where $\ell$ is the index of a thread block in a
grid of the main computational kernel.  In $C_{\ell}^{(0)}$ the count
of the ``big'' transformations, and in $C_{\ell}^{(1)}$ the count of
all transformations applied in all kernel launches within a single
block sweep are accumulated.  At the beginning of each block sweep $C$
is zeroed out on the GPU, and is copied to a CPU vector
$\widetilde{C}$ at the end of the sweep.
%
%
\subsubsection{The shared memory layout}\label{sss:3.1.2}
%
%
For each thread block, the non-constant, non-register data (comprising
three complex matrices: $\mathsf{F}$, $\mathsf{G}$, and $\mathsf{Z}$)
for the main computational kernel is laid out in the shared memory as:
\begin{displaymath}
  \Real(\mathsf{F}), \Real(\mathsf{G}), \Real(\mathsf{Z});\quad
  \Imag(\mathsf{F}), \Imag(\mathsf{G}), \Imag(\mathsf{Z}).
\end{displaymath}
Each double-precision matrix is square, of order $32$, with the
elements stored in Fortran array order, for a total shared memory
requirement of
$(3\times 2)\times(32\times 32)\times 8\,\mathrm{B}=48\,\mathrm{kiB}$
for a thread block.  The shared memory is configured with 8-byte-wide
banks.  No other kernel requires any shared memory.

Let $\Real(\mathsf{F}_{64\times 32}^{})$ stand for the
contiguous memory space occupied by
$\Real(\mathsf{F})$ and
$\Real(\mathsf{G})$;
$\Imag(\mathsf{F}_{64\times 32}^{})$ for
$\Imag(\mathsf{F})$ and
$\Imag(\mathsf{G})$;
$\Real(\mathsf{G}_{64\times 32}^{})$ for
$\Real(\mathsf{G})$ and
$\Real(\mathsf{Z})$; and
$\Imag(\mathsf{G}_{64\times 32}^{})$ for
$\Imag(\mathsf{G})$ and
$\Imag(\mathsf{Z})$, as the real and the imaginary
parts of $\mathsf{F}_{64\times 32}^{}$ and
$\mathsf{G}_{64\times 32}^{}$ matrices that share the same storage
with $\mathsf{F}$, $\mathsf{G}$, and $\mathsf{Z}$.  Such overlapping
of data is necessary for the formation of $\mathsf{F}$ and
$\mathsf{G}$ from the block columns of $F$ and $G$, respectively, as
described below.  Also, let
$\Real(\mathsf{F}_{96\times 32}^{})$ stand for
$\Real(\mathsf{F})$,
$\Real(\mathsf{G})$,
$\Real(\mathsf{Z})$;
and $\Imag(\mathsf{F}_{96\times 32}^{})$ for
$\Imag(\mathsf{F})$,
$\Imag(\mathsf{G})$,
$\Imag(\mathsf{Z})$.

\paragraph{The real case}%
In the real Hari--Zimmermann algorithm $\Imag(\cdot)$ matrices do not
exist, so repacking of the input and the output data does not happen.
The other properties of two data layouts still hold.  The shared
memory requirements are half of those for the complex algorithm, i.e.,
$24\,\mathrm{kiB}$.
%
%
\subsubsection{The constant memory layout}\label{sss:3.1.3}
%
%
The constant memory on the GPU holds the pointers to the matrices and
the vectors described above, with their dimensions, to avoid sending
them as parameters in each kernel call.  The Jacobi strategy table for
the first, pointwise level of the algorithm is also stored in the
constant memory, since it does not depend on the actual input data.

The strategy table contains 31 (or 32) rows, same as the number of
steps of a chosen (quasi-)cyclic parallel strategy.  Each row is an
array of 16 index pairs $(\mathsf{p}, \mathsf{q})$, with
$\mathsf{p} < \mathsf{q}$, where no two indices in a row are the same.
A pair of such indices addresses a pair of columns of the matrices
$\mathsf{F}$ and $\mathsf{G}$ to be transformed concurrently with all
other column pairs in the step.
%
%
\subsubsection{Constants in the global memory}\label{sss:3.1.4}
%
%
The Jacobi strategy table for the second, block level of the algorithm
might not fit in the constant memory for the large $\mathtt{n}$, so it
has to be stored in the global memory in such a case.  It is similarly
formatted as the table for the pointwise level, but with
$\mathtt{n}-1$ (or $\mathtt{n}$) rows, each with $\mathtt{n}/2$ index
pairs.  Here, a pair $(\mathtt{p}, \mathtt{q})$, with
$\mathtt{p} < \mathtt{q}$, addresses a pair of block columns of the
matrices $F$ and $G$.  No two indices in a row are the same, i.e.,
every integer between 0 and $\mathtt{n}-1$ appears exactly once in a
row.  Each row encodes a step of the chosen block level (quasi-)cyclic
parallel strategy, which does not have to be of the same kind as the
one chosen for the pointwise level.

Both tables are precomputed on and preloaded from the
CPU~\cite{Novakovic-2017,Singer-DiNapoli-Novakovic-Caklovic-2020}
before any computation starts on the GPU.
%
%
\subsection{Arithmetic operations}\label{ss:3.2}
%
%
\looseness=-1
Since the data is held in the real-valued arrays only, the complex
arithmetic is performed manually, computing the real and the imaginary
parts of the result separately, rather than assembling the complex
operands in the CUDA format each time an operation has to be
performed, and disassembling the result when it has to be stored back
in memory.
%
%
\subsubsection{Complex arithmetic}\label{sss:3.2.1}
%
%
The arithmetic operations on complex numbers needed by the algorithm
are addition, subtraction, negation, complex conjugation,
multiplication by a complex or a real number (or an inverse of the
latter), and taking the absolute value.  Only $|z|$, $a \cdot b$, and
an FMA-like operation $a \cdot b + c$ (a complex multiplication and an
addition fused) require special attention, while the rest are trivial
to express by the real arithmetic directly in the code.

The absolute value is obtained as
$|z|\assgn\mathop{\mathtt{hypot}}(\Real(z),\Imag(z))$,
without undue overflow.  Still, it is possible that $|z|$ overflows
when at least one component of $z$ is close enough by magnitude to the
largest representable finite double-precision number, but such a
problem can be mitigated by a joint downscaling of two matrices under
transformation.  For example, a scaling by $1/2$ would suffice, and
would also keep the significand intact for all normalized (i.e.,
finite non-subnormal) numbers.  Such rescaling has not been
implemented, though it would not be overwhelmingly hard to apply the
rescaling and restart the computation if any thread detects that its
$|z|$ operation has overflowed, and makes that known to other threads
in a block by a subsequent \texttt{\_\_syncthreads\_count} CUDA
primitive invoked with a Boolean value indicating the presence of an
overflow.

For multiplication, an inlineable routine (\texttt{zmul}) computes
$z\assgn a \cdot b$ and returns the result via two output-only
arguments, referring to $\Real(z)$ and $\Imag(z)$.
With the CUDA FMA intrinsic \texttt{\_\_fma\_rn} it holds
\begin{displaymath}
  \Real(z)=\mathop{\mathtt{\_\_fma\_rn}}(\Real(a),\Real(b),-\!\Imag(a)\cdot\Imag(b)),
\end{displaymath}
computed in a way that requires three floating-point operations but
two roundings only.  Note that the operations are ordered arbitrarily,
thus \texttt{zmul} could also be realized by multiplying the real
parts of the factors first.  $\Imag(z)$ is obtained by
\begin{displaymath}
  \Imag(z)=\mathop{\mathtt{\_\_fma\_rn}}(\Real(a),\Imag(b),\Imag(a)\cdot\Real(b)),
\end{displaymath}
where only two instead of three floating-point operations are
required, with two roundings, and the choice of the real product
arguments is arbitrary.  In total, five operations (of which the
negation is trivial) instead of six are needed.

The FMA-like operation is modeled after the CUDA one in the
\texttt{cuComplex.h} header.  Let $z \assgn a \cdot b + c$.  Then,
\texttt{zfma} routine requires 3 operations with 2 roundings for
\begin{displaymath}
  \begin{aligned}
    d\assgn&\mathop{\mathtt{\_\_fma\_rn}}(-\!\Imag(a),\Imag(b),\Real(c)),\\
    \Real(z)=&\mathop{\mathtt{\_\_fma\_rn}}(\Real(a),\Real(b),d),
  \end{aligned}
\end{displaymath}
and 2 operations with 2 roundings for
\begin{displaymath}
  \begin{aligned}
    d\assgn&\mathop{\mathtt{\_\_fma\_rn}}(\Imag(a),\Real(b),\Imag(c)),\\
    \Imag(z)=&\mathop{\mathtt{\_\_fma\_rn}}(\Real(a),\Imag(b),d).
  \end{aligned}
\end{displaymath}
It holds $\mathop{\mathtt{zfma}}(a,b,0)=\mathop{\mathtt{zmul}}(a,b)$
for all $a$ and $b$.
%
%
\subsubsection{Real arithmetic}\label{sss:3.2.2}
%
%
The real arithmetic uses operations with the accuracy guarantees
mandated by the IEEE~754 standard for floating-point arithmetic in
rounding to nearest (ties to even) mode, except in the optional
enhanced dot-product computation, where rounding to $-\infty$ is also
employed, as described in Appendix~\ref{s:A}.

A correctly rounded (i.e., with the relative error of no more than
half ulp) double-precision
$\mathop{\mathtt{rsqrt}}(x) \assgn 1/\!\sqrt{x}$ device function,
provided by Norbert Juffa in private communication, that improves the
accuracy of the CUDA math library routine of the same name (let it be
referred to by $\mathtt{rsqrt\_rn}$ when a need arises to disambiguate
between the two, and by $\mathtt{rsqrt}$ when either is acceptable),
is called wherever such an expression has to be computed.
%
%
\subsubsection{Reproducibility}\label{sss:3.2.3}
%
%
In both the real and the complex code $\mathtt{rsqrt\_rn}$ function is
expected, but not extensively verified, to be correctly rounded and
thus reproducible.  Reproducibility of the results is guaranteed for
the complex code as long as it is for the $\mathtt{hypot}$ function in
all CUDA versions and on all GPUs under consideration.  All other
floating-point arithmetic operations with rounding (i.e., not
including the comparisons and the negations) are expressed in the
terms of the seven double precision CUDA intrinsics.
%
%
\subsubsection{Integer arithmetic}\label{sss:3.2.4}
%
%
To keep the memory requirements low, the pointwise level indices in
the strategy table are stored as unsigned 1-byte integers, while the
block level indices occupy 2 bytes each (i.e., $\mathtt{n}\le 65536$,
what is enough to exceed the RAM sizes of the present-day GPUs).

For dimensioning and indexing purposes the unsigned 4-byte integers
(after a possible promotion) are used, since their range allows for
addressing up to $32\,\mathrm{GiB}$ of double-precision floating-point
data, which is twice the quantity of GPU RAM available on the testing
hardware.  However, 8-byte integers should be used instead if the
future GPUs provide more memory than this limit.

Although Fortran array order is assumed throughout the paper and the
code, the indices on a GPU are zero-based.  The CUDA thread (block)
indices \texttt{blockIdx.x}, \texttt{threadIdx.x}, and
\texttt{threadIdx.y} are shortened as $\mathtt{b}_{\mathtt{x}}^{}$,
$\mathtt{t}_{\mathtt{x}}^{}$, and $\mathtt{t}_{\mathtt{y}}^{}$,
respectively.
%
%
\subsection{Initialization of $Z$ with optional rescaling of $F$ and $G$}\label{ss:3.3}
%
%
Here, \texttt{initFGZ}, the first of three computational kernels, is
described.  Its purpose is to initialize the matrix $Z$, having been
zeroed out after allocation, to $Z_0^{}$, a diagonal matrix such that
$(Z_0^{})_{\null\!jj}^{} \assgn 1/\|g_{\null\!j}^{}\|_F^{}$, and to
rescale $F$ and $G$ to $F_0^{}$ and $G_0^{}$, by multiplying the
elements of each column $j$ of the matrices by
$(Z_0^{})_{\null\!jj}^{}$ in the variants 0, 1, 4, and 5.  Else, in
other variants, $Z_0^{}=I_n^{}$.

The kernel is launched once, before the iterative phase of the
algorithm, with a one-dimensional grid of $n/2$ thread blocks, each of
which is also one-dimensional, with 64 threads (two warps of 32
consecutive-numbered threads).

A warp is in charge of one column of $F$, $G$, and $Z$, i.e., its
threads access only the elements $i$ of that column $j$, where
\begin{displaymath}
  j\assgn\mathtt{b}_{\mathtt{x}}^{}\cdot 2 + \lfloor\mathtt{t}_{\mathtt{x}}^{}/32\rfloor,\quad
  i\bmod 32 = \mathtt{t}_{\mathtt{x}}^{}\bmod 32.
\end{displaymath}

A warp reads 32 consecutive elements of $\Real(G)_{\null\!j}^{}$
and $\Imag(G)_{\null\!j}^{}$ at a time.  Each of its threads updates
its register-stored partial sums
\begin{displaymath}
  \hat{c}_{\mathrm{r}}'[\mathtt{t}_{\mathtt{x}}^{}]\assgn\hat{c}_{\mathrm{r}}^{}[\mathtt{t}_{\mathtt{x}}^{}]+\Real(G)_{ij}^2,\quad
  \hat{c}_{\mathrm{i}}'[\mathtt{t}_{\mathtt{x}}^{}]\assgn\hat{c}_{\mathrm{i}}^{}[\mathtt{t}_{\mathtt{x}}^{}]+\Imag(G)_{ij}^2,
\end{displaymath}
using one FMA operation for each update, and this is repeated by going
to rows $i\assgn i + 32$ until $i\ge m_G^{}$.  Initially,
$i = \mathtt{t}_{\mathtt{x}}^{}\bmod 32$ and
$\hat{c}_{\mathrm{r}}^{}[\mathtt{t}_{\mathtt{x}}^{}]=\hat{c}_{\mathrm{i}}^{}[\mathtt{t}_{\mathtt{x}}^{}]=0$.
After passing through the entire column, those partial sums are added
to obtain
$\hat{s}[\mathtt{t}_{\mathtt{x}}^{}]\assgn\hat{c}_{\mathrm{r}}^{}[\mathtt{t}_{\mathtt{x}}^{}]+\hat{c}_{\mathrm{i}}^{}[\mathtt{t}_{\mathtt{x}}^{}]$.
Then, $\hat{s}[\mathtt{t}_{\mathtt{x}}^{}]$ are summed and the result
is distributed across the warp by a warp-shuffling~\cite{NVidia-2019}
sum-reduction, described in Appendix~\ref{s:C}, yielding the sum of
squares of the magnitudes of the elements in the column, i.e.,
$\|g_{\null\!j}\|_F^2$.

Such a computation occurs in the variants 0 and 4, while in the
variants 1 and 5 the enhanced dot-product computation as in
Appendix~\ref{s:A} updates the per-thread, register-stored partial
sums
$c_{\mathrm{r}}^{}[\mathtt{t}_{\mathtt{x}}^{}]$,
$c_{\mathrm{i}}^{}[\mathtt{t}_{\mathtt{x}}^{}]$,
$d_{\mathrm{r}}^{}[\mathtt{t}_{\mathtt{x}}^{}]$,
$d_{\mathrm{i}}^{}[\mathtt{t}_{\mathtt{x}}^{}]$.
After a pass over the column completes,
$s[\mathtt{t}_{\mathtt{x}}^{}]$ are formed according to the rules
of Appendix~\ref{s:A} and summed as above.

Either way,
$z_{\null\!j}^{}[\mathtt{t}_{\mathtt{x}}^{}]\assgn 1/\!\sqrt{\|g_{\null\!j}\|_F^2}$
is then computed, and the $j$th columns of $F$ and $G$ are scaled by
$z_{\null\!j}^{}[\mathtt{t}_{\mathtt{x}}^{}]$ in a loop similar to the
one described above, i.e., for $i$ in steps of 32 while $i<m_F^{}$,
\begin{displaymath}
  \begin{aligned}
    \Real(F)_{ij}'&\assgn\Real(F)_{ij}\cdot z_{\null\!j}^{}[\mathtt{t}_{\mathtt{x}}^{}],\\
    \Imag(F)_{ij}'&\assgn\Imag(F)_{ij}\cdot z_{\null\!j}^{}[\mathtt{t}_{\mathtt{x}}^{}],
  \end{aligned}
\end{displaymath}
and then the same scaling is performed on $G$, with $i<m_G^{}$.

Finally, $z_{\null\!j}^{}[\mathtt{t}_{\mathtt{x}}^{}]$ is written to
$\Real(Z)_{\null\!lj}^{}$ by the lowest-numbered thread in a warp,
i.e., $\mathrm{t}_{\mathrm{x}}\equiv 0\pmod{32}$, where $l$ is an
index making a physical column $j$ treated as a logical column $l$.
In the single-GPU case, $l=j$.  In the variants 2, 3, 6, and 7,
$\Real(Z)_{\null\!lj}^{}$ is set to 1 and no other processing occurs.

This and any other computation of the Frobenius norm of a vector via
the sum of squares of its elements could overflow even if the result
itself would not.  See~\cite[Appendix~A]{Novakovic-2015} for one of
several possible remedies.
%
%
\subsection{Rescaling of $Z$ and extraction of $U$, $\Sigma_F^{}$, $V$, $\Sigma_G^{}$, and $\Sigma$}\label{ss:3.4}
%
%
After each block sweep, another kernel, \texttt{rescale}, is called,
with a Boolean flag \texttt{f} indicating whether it is the last
sweep.

If \texttt{f} is \texttt{false}, only $Z$ is rescaled according to the
rules of subsection~\ref{ss:2.4}, and otherwise the full results of
the GSVD computation ($U$, $\Sigma_F^{}$, $V$, $\Sigma_G^{}$, and
$\Sigma$) are produced.

The kernel's grid is identical, and the operation very similar to
\texttt{initFGZ}.  First, $\|f_{\null\!j}^{}\|_F^2$ is computed, and
if non-unity and \texttt{f}, $f_{\null\!j}^{}$ is scaled by
$1/\!\sqrt{\|f_{\null\!j}^{}\|_F^2}$.  If \texttt{f},
$\Sigma_{\null\!j}'[\mathtt{t}_{\mathtt{x}}^{}]\assgn\sqrt{\|f_{\null\!j}^{}\|_F^2}$.
Then, $\|g_{\null\!j}^{}\|_F^2$ is computed, and if non-unity and
\texttt{f}, $g_{\null\!j}^{}$ is scaled by
$1/\!\sqrt{\|g_{\null\!j}^{}\|_F^2}$, as well as
$\Sigma_{\null\!j}'[\mathtt{t}_{\mathtt{x}}^{}]$ to obtain
$\Sigma_{\null\!j}^{}[\mathtt{t}_{\mathtt{x}}^{}]$; else, if
\texttt{f},
$\Sigma_{\null\!j}^{}[\mathtt{t}_{\mathtt{x}}^{}]\assgn\Sigma_{\null\!j}'[\mathtt{t}_{\mathtt{x}}^{}]$.

Then
$\Sigma_{F;j}'[\mathtt{t}_{\mathtt{x}}^{}]\assgn\sqrt{\|f_{\null\!j}^{}\|_F^2}$,
$\Sigma_{G;j}'[\mathtt{t}_{\mathtt{x}}^{}]\assgn\sqrt{\|g_{\null\!j}^{}\|_F^2}$,
and
$\theta[\mathtt{t}_{\mathtt{x}}^{}]\assgn 1/\!\sqrt{\|f_{\null\!j}^{}\|_F^2+\|g_{\null\!j}^{}\|_F^2}$.
If $\theta[\mathtt{t}_{\mathtt{x}}^{}]\ne 1$, $z_{\null\!j}^{}$ is
scaled by $\theta[\mathtt{t}_{\mathtt{x}}^{}]$, as well as
$\Sigma_{F;j}'[\mathtt{t}_{\mathtt{x}}^{}]$ and
$\Sigma_{G;j}'[\mathtt{t}_{\mathtt{x}}^{}]$ to obtain
$\Sigma_{F;j}^{}[\mathtt{t}_{\mathtt{x}}^{}]$ and
$\Sigma_{G;j}^{}[\mathtt{t}_{\mathtt{x}}^{}]$; else,
$\Sigma_{F;j}^{}[\mathtt{t}_{\mathtt{x}}^{}]\assgn\Sigma_{F;j}'[\mathtt{t}_{\mathtt{x}}^{}]$
and
$\Sigma_{G;j}^{}[\mathtt{t}_{\mathtt{x}}^{}]\assgn\Sigma_{G;j}'[\mathtt{t}_{\mathtt{x}}^{}]$.

Finally, if \texttt{f}, $\Sigma_{\null\!j}^{}$, $\Sigma_{F;j}^{}$, and
$\Sigma_{G;j}^{}$ are written to the GPU RAM by a thread
$\mathrm{t}_{\mathrm{x}}\equiv 0\pmod{32}$.  All variables indexed by
$\mathtt{t}_{\mathtt{x}}^{}$ above are per-thread and register-stored,
unless a register spill occurs.
%
%
\subsection{The main computational kernel}\label{ss:3.5}
%
%
The main kernel comes in \texttt{bstep1s} and \texttt{bstep1n}
versions, where the former is the default one, with the column
sorting, while the latter is a non-sorting version.

The kernel is called once per a block step.  Each such call
constitutes the entire block step, and it cannot run concurrently with
any other GPU part of the algorithm since it can update almost the
whole allocated GPU memory.

The kernel's grid is one-dimensional, with $\mathtt{n}/2$
two-dimensional thread blocks, each of them having
$32\times\mathtt{w}=512$ threads.  A thread block
$\ell\assgn\mathtt{b}_{\mathtt{x}}^{}$ in the block step
$\mathtt{k}^{}\assgn k^{}\!\!\mod\mathtt{n}'$ is in charge of one
pivot block column pair,
$(\mathtt{p}_{\mathtt{k}}^{(\ell)},\mathtt{q}_{\mathtt{k}}^{(\ell)})$,
of $F$, $G$, and $Z$, where $\mathtt{n}'$ is $\mathtt{n}-1$ for the
\textsc{me} or $\mathtt{n}$ for the \textsc{mm} strategy kind.

The computational subphases of
\texttt{bstep1}(\texttt{s}/\texttt{n})($k$),
\begin{enumerate}
\item[\textbf{1}.] formation of $\widehat{\mathsf{A}}_k^{(\ell)}$ and
  $\widehat{\mathsf{B}}_k^{(\ell)}$ in the shared memory,
\item[\textbf{2}.] the Cholesky factorizations of
  $\widehat{\mathsf{A}}_k^{(\ell)}$ and
  $\widehat{\mathsf{B}}_k^{(\ell)}$ as
  $\widehat{\mathsf{F}}_k^{(\ell)\ast}\widehat{\mathsf{F}}_k^{(\ell)}$
  and
  $\widehat{\mathsf{G}}_k^{(\ell)\ast}\widehat{\mathsf{G}}_k^{(\ell)}$,
  respectively,
\item[\textbf{3}.] the pointwise implicit Hari--Zimmermann algorithm on
  the matrix pair
  $(\widehat{\mathsf{F}}_k^{(\ell)},\widehat{\mathsf{G}}_k^{(\ell)})$,
  yielding $\widetilde{\mathsf{Z}}_k^{(\ell)}$,
\item[\textbf{4}.] postmultiplication of the pair $\ell$ of pivot block
  columns of $F$, $G$, and $Z$ by $\widetilde{\mathsf{Z}}_k^{(\ell)}$,
\end{enumerate}
are all fused into a single kernel to effortlessly preserve the
contents of the shared memory between them.

All the required matrix algebra routines have been written as device
functions with the semantics similar to, but different from the
standard BLAS, due to the data distribution and the memory
constraints.  For example, a single call of the BLAS-compatible
\texttt{ZHERK} (or \texttt{DSYRK} in the real case) operation for the
subphase~\textbf{1} is not possible, since the two pivot block columns
do not have to be adjacent in the global memory.  The
subphase~\textbf{3} cannot use a single standard \texttt{ZGEMM} (or
\texttt{DGEMM}) call for the same reason, but also because the block
columns have to be overwritten in-place to avoid introducing any work
arrays.

Since no two pivot block index pairs share an index, all thread blocks
can be executed concurrently without any interdependencies or data
races.  Due to the shared memory requirement and a high thread count,
it is not possible that more than two (or, in the real case, four)
thread blocks could share a single GPU multiprocessor (an SM for
short, which cannot have more than 2048 threads resident at present).
On a Maxwell GPU, the profiler reports occupancy of 25\% for the real
and the complex \texttt{bstep1s}, i.e., at most one thread block is
active on an SM at any time.  That can be attributed to a huge
register pressure, since 128 registers per thread are used for the
main kernel (in the variant 0), with a significant amount of spillage,
thus completely exhausting the SM's register file.  Should more than
$2^{16}$ registers be available per SM, it might be possible to
achieve a higher occupancy.

Therefore, for the matrices large enough, only a fraction of all
thread blocks in the grid can execute at the same time on a GPU\@.  It
is a presumption (but not a requirement) that the CUDA runtime shall
schedule a thread block for execution at an early opportunity after a
running one terminates, thereby keeping the GPU busy despite of the
possible execution time variations (i.e., the number of the inner
sweeps and the transformations required) among the thread blocks,
especially in the \textsc{fb} case.

Note that $\mathtt{t}_{\mathtt{y}}^{}$ addresses a warp,
$0\le\mathtt{t}_{\mathtt{y}}^{}<\mathtt{w}$, and
$\mathtt{t}_{\mathtt{x}}^{}$, $0\le\mathtt{t}_{\mathtt{x}}^{}<32$,
denotes a lane (a thread) within the warp.  Throughout a thread block,
each warp is in charge of two ``ordinary'' (i.e., not block) columns,
in the global or in the shared memory, but of which two varies between
and within the subphases.
%
%
\subsubsection{Subphase~\textbf{1} (two \texttt{ZHERK} or \texttt{DSYRK} like operations)}\label{sss:3.5.1}
%
%
The task of this subphase is to form $\widehat{\mathsf{A}}_k^{(\ell)}$
and then $\widehat{\mathsf{B}}_k^{(\ell)}$ in the shared memory,
occupying $\Real(\mathsf{F})$ (and $\Imag(\mathsf{F})$), and
$\Real(\mathsf{G})$ (and $\Imag(\mathsf{G})$), respectively, by a
single pass through the pivot block columns of $F_k^{}$ and $G_k^{}$.
The resulting matrices are Hermitian in theory, but unlike in BLAS,
both the strictly lower and the strictly upper triangle of each matrix
are explicitly computed, even though only the lower triangle is read
in the subphase~\textbf{2}, thus avoiding a possible issue with one
triangle not being the exact transpose-conjugate of the other
numerically.

A warp indexed by $\mathtt{t}_{\mathtt{y}}^{}$ is assigned two
column indices, $p_{\mathtt{y};\mathtt{k}}^{(\ell)}$ and
$q_{\mathtt{y};\mathtt{k}}^{(\ell)}$, in the range of the first and
the second pivot block column, respectively, as
\begin{displaymath}
  p_{\mathtt{y};\mathtt{k}}^{(\ell)}\assgn\mathtt{p}_{\mathtt{k}}^{(\ell)}\cdot\mathtt{w}+\mathtt{t}_{\mathtt{y}}^{},\quad
  q_{\mathtt{y};\mathtt{k}}^{(\ell)}\assgn\mathtt{q}_{\mathtt{k}}^{(\ell)}\cdot\mathtt{w}+\mathtt{t}_{\mathtt{y}}^{}.
\end{displaymath}
Each thread holds four register-stored variables,
\begin{displaymath}
  \mathtt{r}[\mathtt{t}_{\mathtt{x}}^{},\mathtt{t}_{\mathtt{y}}^{}],\quad
  \mathtt{r}[\mathtt{t}_{\mathtt{x}}^{},\mathtt{t}_{\mathtt{y}}'],\quad
  \mathtt{i}[\mathtt{t}_{\mathtt{x}}^{},\mathtt{t}_{\mathtt{y}}^{}],\quad
  \mathtt{i}[\mathtt{t}_{\mathtt{x}}^{},\mathtt{t}_{\mathtt{y}}'],
\end{displaymath}
initially set to zero, that hold the real (first two) and the
imaginary (last two) parts of two (partial) dot-products of the
columns of $F_k^{}$ and, in the second instance, of $G_k^{}$, where
$\mathtt{t}_{\mathtt{y}}'\assgn\mathtt{t}_{\mathtt{y}}^{}+\mathtt{w}$.

In a loop over $i$, starting from $i\assgn\mathtt{t}_{\mathtt{x}}^{}$
and terminating when $i\ge m_F^{}$, with $i\assgn i^{}+64$, in each
step two consecutive chunks of 32 rows (i.e., 64 rows) of the columns
$p_{\mathtt{y};\mathtt{k}}^{(\ell)}$ and
$q_{\mathtt{y};\mathtt{k}}^{(\ell)}$ are read from $\Real(F_k^{})$ and
$\Imag(F_k^{})$ into $\Real(\mathsf{F}_{64\times 32}^{})$ and
$\Imag(\mathsf{F}_{64\times 32}^{})$.  Each lane reads an element from
the global memory and writes it into the shared memory, both in the
coalesced manner, four times per chunk.  The elements of the column
$p_{\mathtt{y};\mathtt{k}}^{(\ell)}$ are stored into the
$\mathtt{t}_{\mathtt{y}}^{}$th column, and those of the column
$q_{\mathtt{y};\mathtt{k}}^{(\ell)}$ are stored into the
$\mathtt{t}_{\mathtt{y}}'$th column of the shared memory buffer.  The
elements of the first chunk are stored into the
$\mathtt{t}_{\mathtt{x}}^{}$th row, and of the second chunk into the
$(\mathtt{t}_{\mathtt{x}}^{}+32)$th row of the buffer.  The thread
block is then synchronized, to complete filling the buffer by all
warps.

An unrolled inner loop over $j$, $0 \le j < 64$, followed by a
synchronization call, updates the local partial dot-products.

For each $j$, let
$\mathtt{t}_{\mathtt{x}}'\assgn(\mathtt{t}_{\mathtt{x}}^{}+j)\!\!\mod 64$,
and
\begin{displaymath}
  z_{\mathtt{y}}^{}\assgn(\mathtt{r}[\mathtt{t}_{\mathtt{x}}^{},\mathtt{t}_{\mathtt{y}}^{}],\mathtt{i}[\mathtt{t}_{\mathtt{x}}^{},\mathtt{t}_{\mathtt{y}}^{}]),\quad
  z_{\mathtt{y}}'\assgn(\mathtt{r}[\mathtt{t}_{\mathtt{x}}^{},\mathtt{t}_{\mathtt{y}}'],\mathtt{i}[\mathtt{t}_{\mathtt{x}}^{},\mathtt{t}_{\mathtt{y}}']),
\end{displaymath}
\begin{displaymath}
  \begin{aligned}
    \mathsf{z}_{\mathtt{x}}^{\ast}&\assgn(\Real(\mathsf{F}_{64\times 32}^{})[\mathtt{t}_{\mathtt{x}}',\mathtt{t}_{\mathtt{x}}^{}],-\Imag(\mathsf{F}_{64\times 32}^{})[\mathtt{t}_{\mathtt{x}}',\mathtt{t}_{\mathtt{x}}^{}]),\\
    \mathsf{z}_{\mathtt{y}}^{}&\assgn(\Real(\mathsf{F}_{64\times 32}^{})[\mathtt{t}_{\mathtt{x}}',\mathtt{t}_{\mathtt{y}}^{}],\Imag(\mathsf{F}_{64\times 32}^{})[\mathtt{t}_{\mathtt{x}}',\mathtt{t}_{\mathtt{y}}^{}]),\\
    \mathsf{z}_{\mathtt{y}}'&\assgn(\Real(\mathsf{F}_{64\times 32}^{})[\mathtt{t}_{\mathtt{x}}',\mathtt{t}_{\mathtt{y}}'],\Imag(\mathsf{F}_{64\times 32}^{})[\mathtt{t}_{\mathtt{x}}',\mathtt{t}_{\mathtt{y}}']).
  \end{aligned}
\end{displaymath}
Two fused multiply-add operations perform the updates
\begin{displaymath}
  z_{\mathtt{y}}^{}\assgn\mathop{\mathtt{zfma}}(\mathsf{z}_{\mathtt{x}}^{\ast},\mathsf{z}_{\mathtt{y}}^{},z_{\mathtt{y}}^{}),\quad
  z_{\mathtt{y}}'\assgn\mathop{\mathtt{zfma}}(\mathsf{z}_{\mathtt{x}}^{\ast},\mathsf{z}_{\mathtt{y}}',z_{\mathtt{y}}').
\end{displaymath}
The first updates constitute a computation of the dot-product of
$\mathtt{t}_{\mathtt{x}}^{}$th and $\mathtt{t}_{\mathtt{y}}^{}$th
column of $\mathsf{F}_{64\times 32}^{}$ and updating the partial sum
$z_{\mathtt{y}}^{}$ with it, while the second ones form the
dot-product of the $\mathtt{t}_{\mathtt{x}}^{}$th and
$\mathtt{t}_{\mathtt{y}}'$th column and update $z_{\mathtt{y}}'$ with
it.  Note that all the rows of the buffer are read exactly once,
albeit in the modular (circular) fashion throughout the loop, with the
different starting offsets in each column to minimize the shared
memory bank conflicts.

When the outer loop over $i$ terminates, $z_{\mathtt{y}}^{}$ and
$z_{\mathtt{y}}'$ are stored into $\mathsf{F}_{64\times 32}^{}$ at
the corresponding indices, and a synchronization barrier is reached,
thus finalizing the formation of $\widehat{\mathsf{A}}_k^{(\ell)}$.
The same procedure is repeated with $G_k^{}$ instead of $F_k^{}$ to
obtain $\widehat{\mathsf{B}}_k^{(\ell)}$, substituting $G$ and
$\mathsf{G}$ for $F$ and $\mathsf{F}$, respectively, in the procedure
described above.  Note that $\mathsf{F}_{96\times 32}^{}$ could
(however, unclear if it should) be used instead of
$\mathsf{F}_{64\times 32}^{}$, i.e., three chunks instead of two would
be read into the buffer and the dot-products of the columns of length
96 instead of 64 would be computed.  That would not be possible,
though, for $\mathsf{G}$, since $\widehat{\mathsf{A}}_k^{(\ell)}$,
once formed, must not be overwritten until the next subphase.

In Figure~\ref{fig:zAhA} the arguments \texttt{A0D}, \texttt{A0J},
\texttt{A1D}, \texttt{A1J}, \texttt{AD}, and \texttt{AJ} stand for the
real and the imaginary planes of the
$p_{\mathtt{y};\mathtt{k}}^{(\ell)}$th and the
$q_{\mathtt{y};\mathtt{k}}^{(\ell)}$th columns of $F_k^{}$, and for
$\Real(\mathsf{F}_{64\times 32}^{})$ and
$\Imag(\mathsf{F}_{64\times 32}^{})$, respectively, in the first call
of the device function.  The same holds for $G_k^{}$ and
$\mathsf{G}_{64\times 32}^{}$ in the second call.  The indices
\texttt{x}, \texttt{y0}, and \texttt{y1} correspond to
$\mathtt{t}_{\mathtt{x}}^{}$, $\mathtt{t}_{\mathtt{y}}^{}$, and
$\mathtt{t}_{\mathtt{y}}'$, respectively, while \texttt{m} is the
number of rows of $F$ or $G$.
\begin{figure}
{\small\begin{verbatim}// F??(A, i, j) = A[?? * j + i] (??=32|64)
// cuD: real, cuJ: imaginary part (double)
__device__ __forceinline__ void zAhA
(const cuD *const __restrict__ A0D,
 const cuJ *const __restrict__ A0J,
 const cuD *const __restrict__ A1D,
 const cuJ *const __restrict__ A1J,
 volatile cuD *const __restrict__ AD,
 volatile cuJ *const __restrict__ AJ,
 const unsigned m, const unsigned x,
 const unsigned y0, const unsigned y1)
{
  cuD y0xD = 0.0, y1xD = 0.0;
  cuJ y0xJ = 0.0, y1xJ = 0.0;
  const unsigned x32 = x + 32u;

  for (unsigned i = x; i < m; i += 32u) {
    // read the 1st 32 x 32 chunk from RAM
    F64(AD, x, y0) = A0D[i];
    F64(AJ, x, y0) = A0J[i];
    F64(AD, x, y1) = A1D[i];
    F64(AJ, x, y1) = A1J[i];

    i += 32u;
    // read the 2nd 32 x 32 chunk from RAM
    F64(AD, x32, y0) = A0D[i];
    F64(AJ, x32, y0) = A0J[i];
    F64(AD, x32, y1) = A1D[i];
    F64(AJ, x32, y1) = A1J[i];
    __syncthreads();

    #pragma unroll
    for (unsigned j = 0u; j < 64u; ++j) {
      const unsigned x_64 =
        (x + j) & 0x3Fu; // (x + j) % 64u
      const cuD _x_hD =  F64(AD, x_64, x);
      const cuJ _x_hJ = -F64(AJ, x_64, x);
      const cuD _y0_D =  F64(AD, x_64, y0);
      const cuJ _y0_J =  F64(AJ, x_64, y0);
      const cuD _y1_D =  F64(AD, x_64, y1);
      const cuJ _y1_J =  F64(AJ, x_64, y1);
      // [complex] y0x = _x_h * _y0_ + y0x
      Zfma(y0xD, y0xJ, _x_hD, _x_hJ,
        _y0_D, _y0_J, y0xD, y0xJ);
      // [complex] y1x = _x_h * _y1_ + y1x
      Zfma(y1xD, y1xJ, _x_hD, _x_hJ,
        _y1_D, _y1_J, y1xD, y1xJ);
    }
    __syncthreads();
  }

  // A^H * A stored into the shared memory
  F32(AD, x, y0) = y0xD;
  F32(AJ, x, y0) = y0xJ;
  F32(AD, x, y1) = y1xD;
  F32(AJ, x, y1) = y1xJ;
  __syncthreads();
}\end{verbatim}}
\caption{A CUDA implementation of the subphase~\textbf{1} ($\mathbb{C}$).}
\label{fig:zAhA}
\end{figure}
%
%
\subsubsection{Subphase~\textbf{2} (two \texttt{ZPOTRF} or \texttt{DPOTRF} like operations)}\label{sss:3.5.2}
%
%
The Cholesky factorization of
$\mathsf{A}\assgn\widehat{\mathsf{A}}_k^{(\ell)}$ or
$\mathsf{B}\assgn\widehat{\mathsf{B}}_k^{(\ell)}$ consists of two
similar, unrolled loops over $j$.  The matrix (in Fortran array order)
is accessed and transformed columnwise to avoid the shared memory bank
conflicts, but then a transpose-conjugate operation must follow on the
computed lower triangular factor to obtain the corresponding upper
triangular one.  Along with the transposition-conjugation, the strictly
lower triangle is zeroed-out, since the following subphase makes no
assumptions about the triangularity of the initial matrices.

The first loop iterates over $0 \le j < \mathtt{w}$.  First, the $j$th
diagonal element of $\Real(\mathsf{A})$, $a_{\null\!jj}^{}$, is read
(the imaginary part is assumed to be zero) if
$\mathtt{t}_{\mathtt{y}}^{}=j$ and $\mathtt{t}_{\mathtt{x}}^{}\ge j$
(i.e., in the threads of the $j$th warp which correspond to the lower
triangle, called the ``active'' threads), and the thread block is then
synchronized.

The active threads then scale the $j$th column below the diagonal,
each thread the real and the imaginary part of its element in the
$\mathtt{t}_{\mathtt{x}}^{}$th row, by $1/\!\sqrt{a_{\null\!jj}^{}}$,
while the diagonal is set to $(\sqrt{a_{\null\!jj}^{}},0)$, and the
thread block is then synchronized.

Next, the columns to the right of the $j$th have to be updated, with
all warps (but not all their threads) participating in the update.
Let $j'\assgn(j+1)+\mathtt{t}_{\mathtt{y}}^{}$. Then, if
$\mathtt{t}_{\mathtt{x}}^{}\ge j'$,
\begin{displaymath}
  \mathsf{A}[\mathtt{t}_{\mathtt{x}}^{},j']\assgn\mathop{\mathtt{zfma}}(-\mathsf{A}[\mathtt{t}_{\mathtt{x}}^{},j],\overline{\mathsf{A}[j',j]},\mathsf{A}[\mathtt{t}_{\mathtt{x}}^{},j']),
\end{displaymath}
and the thread block is synchronized.  However, this only updates the
columns from $j+1$ to $j+\mathtt{w}$.  The same update has to be
performed with $j''\assgn j'+\mathtt{w}$ instead of $j'$, i.e., if
$\mathtt{t}_{\mathtt{x}}^{}\ge j'$ (which also ensures that $j''<32$),
\begin{displaymath}
  \mathsf{A}[\mathtt{t}_{\mathtt{x}}^{},j'']\assgn\mathop{\mathtt{zfma}}(-\mathsf{A}[\mathtt{t}_{\mathtt{x}}^{},j],\overline{\mathsf{A}[j'',j]},\mathsf{A}[\mathtt{t}_{\mathtt{x}}^{},j'']),
\end{displaymath}
and another thread synchronization occurs.

The second loop over $\mathtt{w}\le j<32$ is identical to the first
one, except that $\mathtt{t}_{\mathtt{y}}'$ is used instead of
$\mathtt{t}_{\mathtt{y}}^{}$ and the second updates (of the $j''$th
columns) are not needed since $j''\ge 32$.

The ensuing transpose-conjugate with zeroing-out of the strictly lower
triangle is performed by reading
$\mathsf{A}[\mathtt{t}_{\mathtt{x}}^{},\mathtt{t}_{\mathtt{y}}^{}]$
and
$\mathsf{A}[\mathtt{t}_{\mathtt{x}}^{},\mathtt{t}_{\mathtt{y}}']$ into
the register of the
$[\mathtt{t}_{\mathtt{x}}^{},\mathtt{t}_{\mathtt{y}}^{}]$th thread if
$\mathtt{t}_{\mathtt{x}}^{}\ge\mathtt{t}_{\mathtt{y}}^{}$ and
$\mathtt{t}_{\mathtt{x}}^{}\ge\mathtt{t}_{\mathtt{y}}'$,
respectively (i.e., the indices belong to the lower triangle of
$\mathsf{A}$).  Otherwise, those registers are set to 0.  After
negating the imaginary parts in the former case, the values are
written to
$\mathsf{A}[\mathtt{t}_{\mathtt{y}}^{},\mathtt{t}_{\mathtt{x}}^{}]$
and
$\mathsf{A}[\mathtt{t}_{\mathtt{y}}',\mathtt{t}_{\mathtt{x}}^{}]$,
respectively, unfortunately requiring the shared memory bank
conflicts, and the thread block is synchronized, yielding
$\mathsf{F}\assgn\widehat{\mathsf{F}}_k^{(\ell)}$.  The same procedure
is then repeated with $\mathsf{B}$ instead of $\mathsf{A}$, yielding
$\mathsf{G}\assgn\widehat{\mathsf{G}}_k^{(\ell)}$.
%
%
\subsubsection{Subphase~\textbf{3} (the pointwise one-sided algorithm)}\label{sss:3.5.3}
%
%
The pointwise implicit Hari--Zimmermann algorithm, described in
section~\ref{s:2}, subsections~\ref{ss:2.1}, \ref{ss:2.2}, and the
relevant parts of subsections~\ref{ss:2.4} and \ref{ss:2.5}, is
implemented as follows.

The $\mathtt{t}_{\mathtt{y}}^{}$th warp transforms the pairs of
columns of $\mathsf{F}$, $\mathsf{G}$, and $\mathsf{Z}$ in each inner
step $l'\ge 0$.  Let $l\assgn l'\!\!\mod 31$, since the \textsc{me}
strategy is used exclusively at the inner level in the tests.  Each of
the three pivot pairs comprise the columns indexed by
$\mathsf{p}_{\mathtt{y};l}^{}$ and $\mathsf{q}_{\mathtt{y};l}^{}$,
where the indices are read from the $l$th row of the inner strategy
table at the position $\mathtt{t}_{\mathtt{y}}^{}$.  Within a warp,
the $\mathtt{t}_{\mathtt{x}}^{}$th thread is responsible for the
elements in the $\mathtt{t}_{\mathtt{x}}^{}$th row of those columns.

First, $\mathsf{Z}$ is initialized similarly to the procedure
described in subsection~\ref{ss:3.3}, but on the shared memory level.
In the variants 2, 3, 6, and 7, the diagonal of $\Real(\mathsf{Z})$ is
set to unity, and the rest to zero, by the threads in charge of those
elements.  In the variants 0 and 4, the sum of squares of the
magnitudes of the elements of the columns $\mathsf{g}_{\null\!j}^{}$,
i.e., $\|\mathsf{g}_{\null\!j}^{}\|_F^2$, where
$j\in\{\mathsf{p}_{\mathtt{y};l}^{},\mathsf{q}_{\mathtt{y};l}^{}\}$
and $l=0$, is computed by a sum-reduction as in Appendix~\ref{s:C}.
The thread block is then synchronized.  For each of the two indices
$j$, $\Imag(\mathsf{Z})[\mathtt{t}_{\mathtt{x}}^{},j]$ is set to zero,
as well as $\Real(\mathsf{Z})[\mathtt{t}_{\mathtt{x}}^{},j]$, except
when $\mathtt{t}_{\mathtt{x}}^{}=j$, where
$\Real(\mathsf{Z})[j,j]\assgn 1/\!\sqrt{\|\mathsf{g}_{\null\!j}^{}\|_F^2}$
if $\|\mathsf{g}_{\null\!j}^{}\|_F^2\ne 1$, and one otherwise.
The columns $\mathsf{f}_{\null\!j}^{}$ and $\mathsf{g}_{\null\!j}^{}$
are scaled by $\Real(\mathsf{Z})[j,j]$ if it is not unity, and the
thread block is synchronized.   The similar procedure is applied in
the variants 1 and 5, except that the partial sums of squares are
computed as in Appendix~\ref{s:A} (see subsection~\ref{ss:3.3}), and
summed by a routine from Appendix~\ref{s:C}.

Having thus obtained $\mathsf{F}_0^{}$, $\mathsf{G}_0^{}$, and
$\mathsf{Z}_0^{}$, the iterative part of the algorithm starts, with at
most 30 (\textsc{fb}) or 1 (\textsc{bo}) inner sweeps.  At the start
of each sweep two per-sweep counters, of the ``big'' ($\mathfrak{b}$)
and of all ($\mathfrak{s}$) transformations applied, are reset to
zero.  The counters are kept in each thread, but their values are
synchronized across all threads in a thread block.

In the step $l$ and the warp $\mathtt{t}_{\mathtt{y}}^{}$, let
$i\assgn\mathsf{p}_{\mathtt{y};l}^{}$ and
$j\assgn\mathsf{q}_{\mathtt{y};l}^{}$.  The elements of the three
pivot column pairs are loaded into the registers by each thread
reading its row from the shared memory, after which the thread block
is synchronized.  For each original element, there are two variables
for its real and imaginary parts, and two more variables to hold the
value of the new element after transformation, since the old value
is used twice in computing the new one and thus cannot be
overwritten.  For example,
$\Imag(\mathsf{F})[\mathtt{t}_{\mathtt{x}}^{},i]$ has
$\Imag(\mathsf{F}')[\mathtt{t}_{\mathtt{x}}^{},i]$ as its
counterpart.

The $2\times 2$ pivot submatrices $\widehat{A}_l'$ and
$\widehat{B}_l'$ are then formed.  The diagonal elements are obtained
by computing the squares of the column norms as above, and the
off-diagonal ones are given by the dot-products, either ordinary
(i.e., by sum-reducing the real and the imaginary parts of the
products of an element of the $i$th column conjugated and the
corresponding element of the $j$th column) or enhanced (as
in Appendix~\ref{s:A}) ones.

However, $\widehat{A}_l^{}$ and $\widehat{B}_l^{}$ thus obtained have
to be multiplied by $\widehat{D}_l^{}$ from the left and right in the
variants 2, 3, 6, and 7 to get $\widehat{A}_l'$ and $\widehat{B}_l'$.
If $\widehat{B}_{11;l}^{}\ne 1$, then
$\widehat{A}_{11;l}'\assgn\widehat{A}_{11;l}^{}/\widehat{B}_{11;l}^{}$,
$\widehat{D}_{11;l}^{}\assgn 1/\!\sqrt{\widehat{B}_{11;l}^{}}$, and
$\widehat{A}_{12;l}^{}$, $\widehat{B}_{12;l}^{}$ are scaled by
$\widehat{D}_{11;l}^{}$; otherwise, $\widehat{D}_{11;l}^{}=1$, as it
is in the variants 0, 1, 4, and 5.  If $\widehat{B}_{22;l}^{}\ne 1$,
then
$\widehat{A}_{22;l}'\assgn\widehat{A}_{22;l}^{}/\widehat{B}_{22;l}^{}$,
$\widehat{D}_{22;l}^{}\assgn 1/\!\sqrt{\widehat{B}_{22;l}^{}}$, and
$\widehat{A}_{12;l}^{}$, $\widehat{B}_{12;l}^{}$ are scaled by
$\widehat{D}_{22;l}^{}$; otherwise, $\widehat{D}_{22;l}^{}=1$.

All threads in a warp now have the elements of the pivot submatrices
held in their register-stored variables, and the elements' values are
identical across the warp.  Therefore, the subsequent computation of
$\widehat{Z}_l^{}$ on a per-thread basis also has to produce the same
transformation across the warp.

First it has to be established whether a transformation is warranted.
If the relative orthogonality criterion is satisfied,
$\hat{\mathfrak{s}}$ is set to zero, else to one.  All threads in a
thread block agree if there is some computational work (apart from
merely the optional column sorting) to be done in the current step by
uniformly incrementing $\mathfrak{s}$,
\begin{displaymath}
  \mathfrak{s}\assgn\mathfrak{s}+\mathop{\text{\texttt{\_\_syncthreads\_count}}}(\hat{\mathfrak{s}})/32,
\end{displaymath}
by the number of the thread block's warps with the non-trivial
transformations to be applied.

If $\hat{\mathfrak{s}}=0$ and
$\widehat{A}_{11;l}'<\widehat{A}_{22;l}'$, then
$\mathop{\mathrm{V}}(\mathsf{Y}')[\mathtt{t}_{\mathtt{x}}^{},i]\assgn\mathop{\mathrm{V}}(\mathsf{Y})[\mathtt{t}_{\mathtt{x}}^{},j]$
and
$\mathop{\mathrm{V}}(\mathsf{Y}')[\mathtt{t}_{\mathtt{x}}^{},j]\assgn\mathop{\mathrm{V}}(\mathsf{Y})[\mathtt{t}_{\mathtt{x}}^{},i]$,
where $\mathrm{V}\in\{\Real,\Imag\}$ and
$\mathsf{Y}\in\{\mathsf{F},\mathsf{G},\mathsf{Z}\}$ in
\texttt{bstep1s}.  Then, the values of the new variables are stored in
the shared memory.  When $\hat{\mathfrak{s}}=0$ in \texttt{bstep1n}, or
in \texttt{bstep1s} and $\widehat{A}_{11;l}'\ge\widehat{A}_{22;l}'$,
the new variables take the value of the corresponding old ones, i.e.,
no column swapping occurs.

Otherwise, for $\hat{\mathfrak{s}}=1$, $\widehat{Z}_l'$ is computed
according to a procedure described either in
subsection~\ref{sss:2.1.1} for the complex, or in
subsection~\ref{sss:2.1.2} for the real case.  Then, it is established
whether the criterion C1 (for the variants 0, 1, 2, and 3) or the
criterion C2 (for the variants 4, 5, 6, and 7) indicates that the
transformation is ``small''.  If so, $\hat{\mathfrak{b}}\assgn 0$;
else, $\hat{\mathfrak{b}}\assgn 1$.

If $\widehat{D}_{11;l}^{}\ne 1$, the first row of $\widehat{Z}_l'$ is
scaled by $\widehat{D}_{11;l}^{}$.  If $\widehat{D}_{22;l}^{}\ne 1$,
the second row of $\widehat{Z}_l'$ is scaled by
$\widehat{D}_{22;l}^{}$.  Now the completed transformation
$\widehat{Z}_l^{}$ has to be applied to the pivot columns:
\begin{displaymath}
  \begin{aligned}
    \mathsf{Y}'[\mathtt{t}_{\mathtt{x}}^{},i]&\assgn\mathop{\mathtt{zfma}}(\mathsf{Y}[\mathtt{t}_{\mathtt{x}}^{},j],\widehat{Z}_{21;l}^{},\mathsf{Y}[\mathtt{t}_{\mathtt{x}}^{},i]\cdot\Real(\widehat{Z}_{11;l}^{})),\\
    \mathsf{Y}'[\mathtt{t}_{\mathtt{x}}^{},j]&\assgn\mathop{\mathtt{zfma}}(\mathsf{Y}[\mathtt{t}_{\mathtt{x}}^{},i],\widehat{Z}_{12;l}^{},\mathsf{Y}[\mathtt{t}_{\mathtt{x}}^{},j]\cdot\Real(\widehat{Z}_{22;l}^{})),
  \end{aligned}
\end{displaymath}
where $\mathsf{Y}\in\{\mathsf{F},\mathsf{G},\mathsf{Z}\}$.  If one or
both scaled cosines lying on the diagonal of $\widehat{Z}_l^{}$ are
equal to one, the transformation can be (and is) simplified by
removing the corresponding multiplications without numerically
affecting the result.

In \texttt{bstep1s}, to determine if the column swap is required,
the squares of the norms of the transformed columns of $\mathsf{F}$
are computed as the sum-reduced sums of squares of the magnitudes of
the new ($\mathsf{F}'$) elements, depending on the variant.  Those two
values are however not stored for the next step, because that would
require an additional shared memory workspace that might not be
available on all supported architectures.

In the real case it is easy to compute instead the
transformed diagonal elements of the first pivot submatrix
directly~\cite{Novakovic-Singer-Singer-2015}:
\begin{displaymath}
  \begin{aligned}
    a_{11}''&\assgn\widehat{Z}_{11;l}^{\null\,2}\widehat{A}_{11;l}'
    + 2\widehat{Z}_{11;l}^{}\widehat{Z}_{21;l}^{}\widehat{A}_{12;l}'
    + \widehat{Z}_{21;l}^{\null\,2}\widehat{A}_{22;l}',\\
    a_{22}''&\assgn\widehat{Z}_{12;l}^{\null\,2}\widehat{A}_{11;l}'
    + 2\widehat{Z}_{22;l}^{}\widehat{Z}_{12;l}^{}\widehat{A}_{12;l}'
    + \widehat{Z}_{22;l}^{\null\,2}\widehat{A}_{22;l}',
  \end{aligned}
\end{displaymath}
and to swap the $i$th and the $j$th column when $a_{11}''<a_{22}''$.

If the norm of the $i$th column is smaller than the norm of the $j$th
column, then the values of $\mathsf{Y}'[\mathtt{t}_{\mathtt{x}}^{},i]$
and $\mathsf{Y}'[\mathtt{t}_{\mathtt{x}}^{},j]$ are swapped via an
intermediary variable.  Else, or in \texttt{bstep1n}, no swaps occur.
The values of the new variables are then stored in the shared memory,
and $\mathfrak{b}$ is uniformly incremented across the thread block,
\begin{displaymath}
  \mathfrak{b}\assgn\mathfrak{b}+\mathop{\text{\texttt{\_\_syncthreads\_count}}}(\hat{\mathfrak{b}})/32.
\end{displaymath}
The $l$th step is now complete.

At the end of a sweep, if $\hat{s}=0$, the loop is terminated.  Else,
the counters $\mathfrak{S}$ and $\mathfrak{B}$, set at the start of
this subphase to zero, are incremented by $\mathfrak{s}$ and
$\mathfrak{b}$, respectively.

The same rescaling as in subsection~\ref{ss:3.4} with
\texttt{f=false}, but performed on the shared memory, yields
$\widetilde{\mathsf{Z}}_k^{(\ell)}$.  Using the last values of
$\mathsf{F}'[\mathtt{t}_{\mathtt{x}}^{},i]$ and
$\mathsf{G}'[\mathtt{t}_{\mathtt{x}}^{},i]$, the squares of the norms
of the $i$th column of $\mathsf{F}'$ and $\mathsf{G}'$, respectively,
are computed.  Then, $\mathsf{Z}'[\mathtt{t}_{\mathtt{x}}^{},i]$ is
read (or its last value is used), scaled by
$1/\!\sqrt{\|\mathsf{f}_i'\|_F^2+\|\mathsf{g}_i'\|_F^2}$, stored, and
the thread block is synchronized.  The same procedure is repeated with
$j$ instead of $i$, giving
$\widetilde{\mathsf{Z}}_k^{(\ell)}\assgn\mathsf{Z}'$.

A thread with
$\mathtt{t}_{\mathtt{x}}^{}=\mathtt{t}_{\mathtt{y}}^{}=0$
stores $\mathfrak{S}$ and $\mathfrak{B}$ into $C$ as
\begin{displaymath}
  C[2\cdot\mathtt{b}_{\mathtt{x}}^{}]\assgn\mathfrak{S},\quad
  C[2\cdot\mathtt{b}_{\mathtt{x}}^{}+1]\assgn\mathfrak{B},
\end{displaymath}
and finally the thread block is synchronized.
%
%
\subsubsection{Subphase~\textbf{4} (three postmultiplications)}\label{sss:3.5.4}
%
%
In this subphase the pivot block columns of $F_k^{}$, $G_k^{}$, and
$Z_k^{}$ are multiplied by $\widetilde{\mathsf{Z}}_k^{(\ell)}$ and
overwritten by the respective results.

Each multiplication of a pair of pivot block columns (residing in the
global memory) by $\widetilde{\mathsf{Z}}_k^{(\ell)}$ (residing in the
shared memory in $\mathsf{Z}$) and the following update are performed
by a single pass over (i.e., a single read from and a single write to)
the block columns, using the Cannon-like algorithm~\cite{Cannon-69}
for parallel multiplication of two square matrices.

Reading the chunks of a block column pair from the global memory is
identical to the one from the subphase~\textbf{1} in
subsection~\ref{sss:3.5.1}, except that in each iteration of the outer
loop (over $i$) only one chunk is read to $\mathsf{Y}$, instead of two
(which would also be a possibility).  The number of loop iterations
(in parenthesis) depends on the number of rows of $F_k^{}$
($m_F^{}/32$), $G_k^{}$ ($m_G^{}/32$), and $Z_k^{}$ ($n/32$).  Here,
$\mathsf{Y}$ is $\mathsf{F}$ when updating $F_k^{}$ and $Z_k^{}$, and
$\mathsf{G}$ when updating $G_k^{}$.  The thread block is then
synchronized.

The per-thread variables to hold the product of the current chunk with
$\mathsf{Z}$ are set to zero.  Each thread is in charge of forming the
elements with indices
$[\mathtt{t}_{\mathtt{x}}^{},\mathtt{t}_{\mathtt{y}}^{}]$ and
$[\mathtt{t}_{\mathtt{x}}^{},\mathtt{t}_{\mathtt{y}}']$ of the
product $\Pi$.

The initial skews are defined as
$\imath\assgn(\mathtt{t}_{\mathtt{y}}^{}+\mathtt{t}_{\mathtt{x}}^{})\!\!\mod 32$
and
$\imath'\assgn(\mathtt{t}_{\mathtt{y}}'+\mathtt{t}_{\mathtt{x}}^{})\!\!\mod 32$.
Then, in each iteration of the unrolled inner loop over $0\le j<32$
the local elements of $\Pi$ are updated,
\begin{displaymath}
  \begin{aligned}
    \Pi[\mathtt{t}_{\mathtt{x}}^{},\mathtt{t}_{\mathtt{y}}^{}]&\assgn\mathop{\mathtt{zfma}}(\mathsf{Y}[\mathtt{t}_{\mathtt{x}}^{},\imath],\mathsf{Z}[\imath,\mathtt{t}_{\mathtt{y}}^{}],\Pi[\mathtt{t}_{\mathtt{x}}^{},\mathtt{t}_{\mathtt{y}}^{}]),\\
    \Pi[\mathtt{t}_{\mathtt{x}}^{},\mathtt{t}_{\mathtt{y}}']&\assgn\mathop{\mathtt{zfma}}(\mathsf{Y}[\mathtt{t}_{\mathtt{x}}^{},\imath'],\mathsf{Z}[\imath',\mathtt{t}_{\mathtt{y}}'],\Pi[\mathtt{t}_{\mathtt{x}}^{},\mathtt{t}_{\mathtt{y}}']),
  \end{aligned}
\end{displaymath}
and $\imath$ and $\imath'$ are cyclically shifted as
$\imath\assgn(\imath+1)\!\!\mod 32$ and
$\imath'\assgn(\imath'+1)\!\!\mod 32$.  When the inner loop
terminates, the thread block is synchronized.

The local values of $\Pi$ now have to be written back to the global
memory, where
$\Pi[\mathtt{t}_{\mathtt{x}}^{},\mathtt{t}_{\mathtt{y}}^{}]$
overwrites $Y_k^{}[i,p_{\mathtt{y};\mathtt{k}}^{(\ell)}]$, while
$\Pi[\mathtt{t}_{\mathtt{x}}^{},\mathtt{t}_{\mathtt{y}}']$ overwrites
$Y_k^{}[i,q_{\mathtt{y};\mathtt{k}}^{(\ell)}]$, for $Y$ being one of
$\{F,G,Z\}$.  The thread block is then synchronized and the next outer
iteration, if any are left, follows.

This procedure is called thrice to update $F_k^{}$, $G_k^{}$ and
$Z_k^{}$, after which the kernel execution (i.e., the $k$th outer
step) terminates and the control returns to the CPU\@.

Figure~\ref{fig:zPostMult} shows the postmultiplication device
function, where the arguments \texttt{A0D}, \texttt{A0J},
\texttt{A1D}, and \texttt{A1J} have the same meaning as in
Figure~\ref{fig:zAhA}, but the columns of $Z_k^{}$ are also expected.
A shared memory buffer, in which the $32\times 32$ chunks of a block
column pair are loaded and packed, is pointed to by \texttt{AD} and
\texttt{AJ}, while \texttt{BD} and \texttt{BJ} point to the
accumulated transformation matrix from the subphase~\textbf{3}, by
which the postmultiplication has to take place.  The product of
matrices \texttt{A} and \texttt{B} overwrites the respective chunk of
the original block columns before another chunk is loaded.
\begin{figure}
{\small\begin{verbatim}__device__ __forceinline__ void zPostMult
(cuD *const __restrict__ A0D,
 cuJ *const __restrict__ A0J,
 cuD *const __restrict__ A1D,
 cuJ *const __restrict__ A1J,
 volatile cuD *const __restrict__ AD,
 volatile cuJ *const __restrict__ AJ,
 volatile const cuD *const __restrict__ BD,
 volatile const cuJ *const __restrict__ BJ,
 const unsigned x, const unsigned y0,
 const unsigned y1, const unsigned m)
{
  // Cannon-like C = A * B
  for (unsigned i = x; i < m; i += 32u) {
    F32(AD, x, y0) = A0D[i];
    F32(AJ, x, y0) = A0J[i];
    F32(AD, x, y1) = A1D[i];
    F32(AJ, x, y1) = A1J[i];
    __syncthreads();

    cuD Cxy0D = 0.0, Cxy1D = 0.0;
    cuJ Cxy0J = 0.0, Cxy1J = 0.0;
    unsigned // skew (mod 32)
      p0 = ((y0 + x) & 0x1Fu),
      p1 = ((y1 + x) & 0x1Fu);

    // multiply and cyclic shift (mod 32)
    #pragma unroll
    for (unsigned k = 0u; k < 32u; ++k) {
      Zfma(Cxy0D, Cxy0J,
        F32(AD, x, p0), F32(AJ, x, p0),
        F32(BD, p0, y0), F32(BJ, p0, y0),
        Cxy0D, Cxy0J);
      Zfma(Cxy1D, Cxy1J,
        F32(AD, x, p1), F32(AJ, x, p1),
        F32(BD, p1, y1), F32(BJ, p1, y1),
        Cxy1D, Cxy1J);
      p0 = (p0 + 1u) & 0x1Fu;
      p1 = (p1 + 1u) & 0x1Fu;
    }
    __syncthreads();

    A0D[i] = Cxy0D; A0J[i] = Cxy0J;
    A1D[i] = Cxy1D; A1J[i] = Cxy1J;
    __syncthreads();
  }
}\end{verbatim}}
\caption{A CUDA C implementation of the subphase~\textbf{4} ($\mathbb{C}$).}
\label{fig:zPostMult}
\end{figure}
%
%
\subsubsection{Dataflow and the shared memory perspective}\label{sss:3.5.5}
\begin{figure*}[h!btp]
  \centering
  \includegraphics[height=.9515\textheight,keepaspectratio]{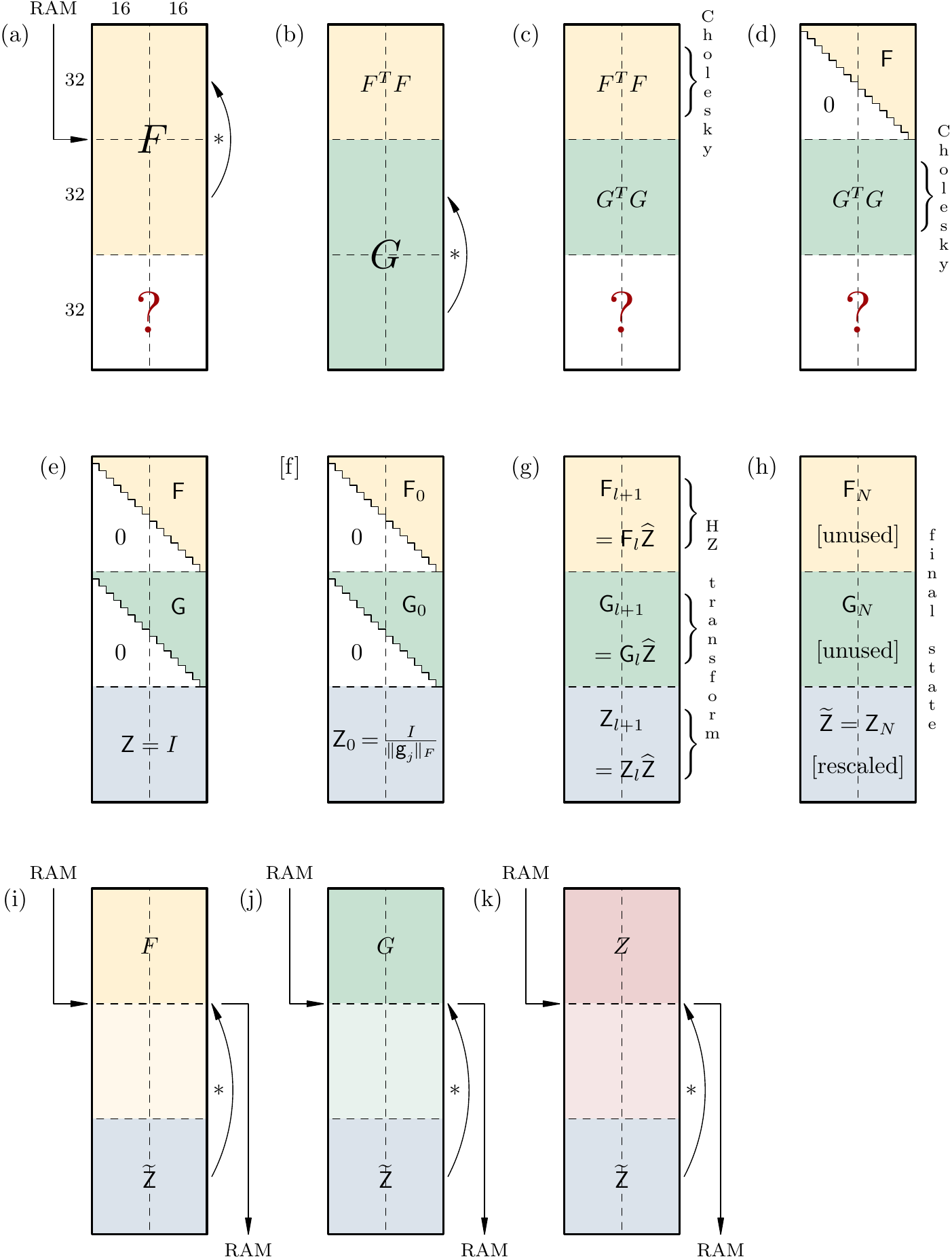}
  \caption{The sequence of operations performed on the shared memory
    of a GPU's multiprocessor by a thread block of the main
    computational kernel.}
  \label{fig:DHZ-shmem1}
\end{figure*}

%
%
In Figure~\ref{fig:DHZ-shmem1} the subphases in the simpler, real case
are summarized from a perspective of the data in the shared memory and
the transformations that it undergoes.
The subfigures (a--b) show the chunks of data coming from the global
memory to form $\mathsf{A}$ and $\mathsf{B}$.  Two Cholesky
factorizations are shown in the subfigures (c--d), after which the
subfigures (e--h) correspond to the subphase~\textbf{3}.  The
subfigure~[f] shows the optional (depending on the variant) prescaling
of $\mathsf{F}$ and $\mathsf{G}$.  The final three subfigures, (i--k),
show the three postmultiplications taking place, with a chunk
of data read from and written to the global memory.  Two different
hues of the first and the middle parts of the shared memory depiction
indicate that the chunks with 64 instead of 32 rows might be used
there.
%
%
\subsection{The CPU part of the algorithm}\label{ss:3.6}
%
%
Algorithm~\ref{alg:cpu1} summarizes the CPU-side actions with a single
GPU\@.
\begin{algorithm*}[h!bt]
  \SetKwFunction{MemCpy}{cudaMemcpyAsync}
  \SetKwFunction{MemSet}{cudaMemsetAsync}
  \SetKwFunction{StreamSync}{cudaStreamSynchronize}
  \SetKwFunction{InitFGZ}{initFGZ}
  \SetKwFunction{BstepI}{bstep1(s/n)}
  \SetKwFunction{Rescale}{rescale}
  \SetKw{Break}{break}
  \InitFGZ{}\tcp*[r]{execute once on a single GPU in the stream \hbox{\sc s} to get $F_0^{},\;G_0^{},\;Z_0^{}$}
  $\widetilde{\mathfrak{S}}\assgn 0;\quad\widetilde{\mathfrak{B}}\assgn 0$\tcp*[r]{initialize the global convergence statistics}
  \For(\tcp*[f]{outer sweep $c$}){$0\le c<S$}{
    \MemSet{$C,0,\mathtt{n}\cdot\hbox{\tt sizeof*}C,\hbox{\sc s}$}\tcp*[r]{zero-out $C$}
    \For(\tcp*[f]{a loop over the outer steps}){$0\le k<\mathtt{n}'$}{
      \BstepI{$k$}\tcp*[r]{each main kernel call in \hbox{\sc s} transforms $\mathtt{n}/2$ block pivots}
    }
    \MemCpy{$\widetilde{C},C,\mathtt{n}\cdot\hbox{\tt sizeof*}C,\hbox{\tt cudaMemcpyDeviceToHost},\hbox{\sc s}$}\tcp*[r]{retrieve $C$}
    \StreamSync{\textsc{s}}\tcp*[r]{synchronize \hbox{\sc s}}
    $\tilde{\mathfrak{s}}\assgn 0;\quad\tilde{\mathfrak{b}}\assgn 0$\tcp*[r]{initialize the sweep convergence statistics}
    \For(\tcp*[f]{gather the sweep convergence statistics}){$0\le i<\mathtt{n}/2$}{
      $\tilde{\mathfrak{s}}\assgn\tilde{\mathfrak{s}}+\widetilde{C}[2\cdot i];\quad\tilde{\mathfrak{b}}\assgn\tilde{\mathfrak{b}}+\widetilde{C}[2\cdot i+1]$\;
    }
    $\widetilde{\mathfrak{S}}\assgn\widetilde{\mathfrak{S}}+\tilde{\mathfrak{s}};\quad\widetilde{\mathfrak{B}}\assgn\widetilde{\mathfrak{B}}+\tilde{\mathfrak{b}}$\tcp*[r]{update the global convergence statistics}
    \lIf(\tcp*[f]{no \underline{big} transformations performed in the sweep}){$\tilde{\mathfrak{b}}=0$}{\Break}
    \Rescale{$\hbox{\tt false}$}\tcp*[r]{a fast rescaling of $Z$ in \hbox{\sc s}}
  }
  \Rescale{$\hbox{\tt f}$}\tcp*[r]{queued in \hbox{\sc s} with \texttt{f=true} on one GPU and \texttt{f=false} on multiple GPUs}
  \caption{The CPU part of the single-GPU implicit Hari--Zimmermann algorithm with data in place.}
  \label{alg:cpu1}
\end{algorithm*}
The same routine is called within the multi-GPU algorithm, except for
allowing $S$ to be a parameter (not the constant $30$ as it is assumed
here), a variation of the final rescaling of $Z$, and some differences 
regarding the \texttt{initFGZ} call, described in
subsections~\ref{ss:3.3} and \ref{ss:4.2}.  The copy-ins and copy-outs
of the majority of data, as well as the initialization of the constant
memory, are left out from Algorithm~\ref{alg:cpu1} but are included in
the single-GPU timing in subsection~\ref{s:5}.

Apart from the copy-ins, copy-outs, and zeroings of data, there is no
scope for using more than one CUDA stream.  All GPU operations can be
performed in any predefined stream \textsc{s} (e.g., in the default
one if no other has been explicitly chosen).  Also, as no GPU
computation, except the fast \texttt{rescale}, can be overlapped with
any CPU task, the execution time of the algorithm depends almost
solely on the GPU performance and the time required to set up a kernel
call.

Each kernel invocation and each sequence of memory copy/set operations
is followed in the testing code by a (generally redundant, except
where noted in Algorithm~\ref{alg:cpu1})
\texttt{cudaStreamSynchronize} call on the chosen stream \textsc{s},
to keep the CPU-side timing consistent (but maybe higher than it is
necessary).

Except the reductions of $\tilde{\mathfrak{s}}$ and
$\tilde{\mathfrak{b}}$ for very large matrices, no other part of the
algorithm might benefit from being executed multi-threadedly on the
CPU\@.  From the CPU perspective, the algorithm is therefore purely
sequential.

Note that the algorithm stops when $\tilde{\mathfrak{b}}=0$, i.e.,
when no big transformations occurred in an outer sweep.  The ``global''
counts $\widetilde{\mathfrak{S}}$ and $\widetilde{\mathfrak{B}}$ of
all and of big transformations applied during the execution of
Algorithm~\ref{alg:cpu1} are only informative here, but they are
consulted in the multi-GPU algorithm's stopping criterion.
%
%
\section{The multi-GPU implementation}\label{s:4}
%
%
When the input data is larger than the available RAM on a single GPU,
it is necessary to either split the workload among the several GPUs,
or resort to some out-of-core technique for swapping the parts of data
in and out of the GPU as the computation progresses.  Here, only the
former approach is followed, since it is simpler, more efficient and
widely applicable now when the multi-GPU installations are becoming
ubiquitous.  In the case when not enough GPUs are available for the
input data to be distributed across them, see an outline of a possible
out-of-core single-GPU algorithm in Appendix~\ref{s:B}.

There is no single, best and straightforward way of generalizing a
single-GPU algorithm to multiple GPUs.  For the (ordinary and
hyperbolic) SVD, the approach in~\cite{Novakovic-2015} was to
distribute the matrix over the GPUs, shorten the assigned part of the
matrix (the Grammian formation being done by cuBLAS, and the ensuing
Cholesky factorization by MAGMA~\cite{Tomov-Dongarra-Baboulin-2010}),
run the single-GPU algorithm on the shortened part, update the
original (non-shortened) columns, and redistribute the parts.  Despite
its decent performance, such a three-level algorithm suffered from the
increased memory usage and some numerical difficulties, both with the
stopping criteria and with the relative accuracy obtained.

A different approach is taken here, to achieve the optimal GPU
\emph{and\/} CPU memory usage (without any work arrays) and a better
accuracy, but with a possible performance penalty induced by
transforming the tall and skinny parts of the matrices directly,
without any shortening.  As no floating-point computation is performed
on the CPU, while the GPU computation still does not rely on any
numerical libraries, it is guaranteed that the results stay bitwise
identical in the repeated runs over the same input data with the same
parameters and the same number of GPUs.
%
%
\subsection{Algorithm setup}\label{ss:4.1}
%
%
In this subsection the multi-GPU computational environment, the input
and the output data distribution across the CPUs and the GPUs, the
communication-aware Jacobi strategies, and the algorithm's
initialization are explained.
%
%
\subsubsection{MPI environment}\label{sss:4.1.1}
%
%
Unlike in~\cite{Novakovic-2015}, where the multiple GPUs were assumed
to belong to the same node, and thus a separate CPU thread of a single
process could be dedicated for driving the computation on each GPU,
here a more flexible solution has been chosen, by assigning to a GPU a
single-threaded MPI~\cite{mpi31} process.  The GPUs can thus be
selected from any number of nodes, with a different number of them on
each node.  Also, the GPUs are not required to be architecturally
identical or even similar across the processes in an MPI job, as long
as they all have enough RAM available.

The count of GPUs, and thus the governing MPI processes
($\mathbf{s}$), for the multi-GPU algorithm is not constrained in
principle, save for being at least two (otherwise, the single-GPU
algorithm is sufficient), and small enough so that at least two (but
for the reasons of performance, a multiple of 32) columns of each
matrix are available to each process, when the matrices are divided
among them columnwise, as described below.  The upper bound on the
number of GPUs is a tunable parameter in the code, while the MPI
implementation might have its own limit on the number of processes in
a job.

\looseness=-1
The MPI processes need not be arranged in any special topology.  Only
the predefined \texttt{MPI\_COMM\_WORLD} communicator is used.  A GPU
and its governing process are jointly referred to by the process' rank
($\mathbf{r}$) in that communicator.

It is assumed that the MPI distribution is CUDA-aware in a sense that
sending data from the GPU RAM of one process and receiving it in the
CPU RAM of another (including the same) process is possible with the
standard MPI point-to-point communication routines (i.e., no manual
CPU buffering of the GPU data is necessary).

Also, the number of elements of each local submatrix has to be at most
\texttt{INT\_MAX}, which at present is the upper limit on the count of
elements that can be transferred in a single MPI
operation~\cite{Hammond-Schaefer-Latham-2014}.  That limit is easily
circumvented by transferring the data in several smaller chunks, but
such chunking has not been implemented since it was not needed for the
amount of RAM ($16\,\mathrm{GiB}$) of the GPUs used for testing.  That
issue will have to be addressed for the future GPUs.
%
%
\subsubsection{Data distribution}\label{sss:4.1.2}
%
%
The matrices $F$, $G$, and $Z$, and the vectors $\Sigma_F^{}$,
$\Sigma_G^{}$, and $\Sigma$, are assumed to always stay distributed
among the MPI processes, i.e., at no moment they are required to be
present in entirety in any subset of the processes.  The amount of the
CPU and the GPU memory required is identical (i.e., not depending on
$\mathbf{r}$), constant throughout the computation, and derivable in
advance from $m_F^{}$, $m_G^{}$, $m_Z^{}\assgn n$, and $\mathbf{s}$
for all processes.

If $n\!\!\mod\mathbf{s}\ne 0$ or $(n/\mathbf{s})\!\!\mod 32\ne 0$, the
matrices $F$, $G$, and $Z$ are bordered as described in
subsection~\ref{ss:3.1}, but requiring that the enlarged $n$ satisfy
$n\!\!\mod(32\cdot\mathbf{s})=0$.  Similarly, the bordering is
required if $m_F^{}\!\!\mod 64\ne 0$ or $m_G^{}\!\!\mod 64\ne 0$.

The columns of the bordered matrices can be distributed evenly among
the processes, such that each process is assigned
$\mathbf{n}\assgn n/\mathbf{s}$ columns.  Let
$\mathbf{w}\assgn\mathbf{n}/2$ consecutive columns of an entire matrix
be called a stripe, to avoid reusing the term ``block column''.  Then,
a process holds two, not necessarily consecutive, stripes of each
matrix, logically separate but with their real parts physically joined
in the same memory allocation, as well as their imaginary parts.  The
dimensions of two joined stripes, one following the other in the
Fortran array order, of $\Real(F)/\Imag(F)$
($m_F^{}\times\mathbf{n}$), $\Real(G)/\Imag(G)$
($m_G^{}\times\mathbf{n}$), and $\Real(Z)/\Imag(Z)$
($n\times\mathbf{n}$), fit the requirements for the input data of the
single-GPU algorithm.

The CPU RAM of the $\mathbf{r}$th process thus holds two
joined stripes in $\Real(F^{(\mathbf{r})})$,
$\Imag(F^{(\mathbf{r})})$, $\Real(G^{(\mathbf{r})})$,
$\Imag(G^{(\mathbf{r})})$, $\Real(Z^{(\mathbf{r})})$, and
$\Imag(Z^{(\mathbf{r})})$ allocations.  The same memory arrangement is
present in the GPU RAM, which holds the same stripes undergoing the
transformations and joined in the allocations
$\Real(F^{[\mathbf{r}]})$, $\Imag(F^{[\mathbf{r}]})$,
$\Real(G^{[\mathbf{r}]})$, $\Imag(G^{[\mathbf{r}]})$,
$\Real(Z^{[\mathbf{r}]})$, and $\Imag(Z^{[\mathbf{r}]})$.  The first
stripe within an allocation is denoted by the index $\mathbf{0}$, and
the second one by the index $\mathbf{1}$; e.g.,
$\Imag(G_{\mathbf{1}}^{[\mathbf{r}]})$ is the second stripe in
$\Imag(G^{[\mathbf{r}]})$.

A similar distribution is in place for $\Sigma_F^{}$, $\Sigma_G^{}$,
and $\Sigma$, where each process holds $\Sigma_F^{(\mathbf{r})}$,
$\Sigma_G^{(\mathbf{r})}$, and $\Sigma^{(\mathbf{r})}$ in the CPU RAM,
and $\Sigma_F^{[\mathbf{r}]}$, $\Sigma_G^{[\mathbf{r}]}$, and
$\Sigma^{[\mathbf{r}]}$ in the GPU RAM, where each allocation is of
length $\mathbf{n}$ and is unused in the algorithm until after the
last step.  Each process also has its convergence vectors
$C^{(\mathbf{r})}$ and $C^{[\mathbf{r}]}$, of length
$\mathbf{n}/\mathtt{w}$, in the CPU and in the GPU RAM, respectively.
%
%
\subsubsection{Communication-aware Jacobi strategies}\label{sss:4.1.3}
%
%
The parallel Jacobi strategies, as defined in subsection~\ref{ss:2.2},
do not contain any explicit information on how to progress from one
step to another by exchanging the pivot (block) columns among the
tasks in a distributed memory environment.  However, such a
communication pattern can be easily retrieved by looking at each two
successive steps, $k$ and $k'\assgn(k+1)\!\!\mod s$, and for each task
$\ell$ in the $k$th step finding the tasks $\ell'$ and $\ell''$ in the
$k'$th step that are to hold either $i_k^{(\ell)}$th or
$j_k^{(\ell)}$th (block) column.

Therefore, given either \textsc{me} or \textsc{mm} strategy table for
the order $\mathfrak{n}\assgn n/\mathbf{w}$ (with the stripes seen as
the block columns), each process independently computes and encodes
the strategy's communication pattern before the start of the
algorithm.  Such a computation requires $O(\mathfrak{n}^3)$
comparisons, but since $\mathfrak{n}$ is a small number an
unacceptable overhead is not incurred.  The computation can be (but it
has not been) parallelized on a CPU, e.g., by turning the outer loop
over $k$ into a parallel one.  The multi-GPU algorithm then references
the following encoded mapping when progressing from one step to the
next.

After each outermost (multi-GPU) step $\mathbf{k}$ of the algorithm,
the first of each two joined stripes on the $\mathbf{r}$th GPU has to
be transferred to either the first or the second stripe on the
$\mathbf{t}^{\{\mathbf{0}\}}$th CPU, for some
$\mathbf{t}^{\{\mathbf{0}\}}$.  Similarly, the second stripe on the
$\mathbf{r}$th GPU has to be transferred to either the first or the
second stripe on the $\mathbf{t}^{\{\mathbf{1}\}}$th CPU, for some
$\mathbf{t}^{\{\mathbf{1}\}}$, establishing a mapping
\begin{displaymath}
  (\mathbf{k},\mathbf{r})\mapsto(\mathbf{p}_{\mathbf{k},\mathbf{r}}^{},\mathbf{q}_{\mathbf{k},\mathbf{r}}^{},\mathfrak{t}_{\mathbf{k},\mathbf{r}}^{\{\mathbf{0}\}},\mathfrak{t}_{\mathbf{k},\mathbf{r}}^{\{\mathbf{1}\}}),
\end{displaymath}
where
$\mathfrak{t}_{\mathbf{k},\mathbf{r}}^{\{\mathbf{0}\}}\assgn\pm(\mathbf{t}_{\mathbf{k},\mathbf{r}}^{\{\mathbf{0}\}}+1)$,
$\mathfrak{t}_{\mathbf{k},\mathbf{r}}^{\{\mathbf{1}\}}\assgn\pm(\mathbf{t}_{\mathbf{k},\mathbf{r}}^{\{\mathbf{1}\}}+1)$,
while $\mathbf{p}_{\mathbf{k},\mathbf{r}}^{}$ and
$\mathbf{q}_{\mathbf{k},\mathbf{r}}^{}$ are indices of the first and
the second stripe in the entire matrix,
$0\le\mathbf{p}_{\mathbf{k},\mathbf{r}}^{}<\mathbf{q}_{\mathbf{k},\mathbf{r}}^{}<\mathfrak{n}$,
and
$\mathbf{t}_{\mathbf{k},\mathbf{r}}^{\{\mathbf{0}\}}\ne\mathbf{t}_{\mathbf{k},\mathbf{r}}^{\{\mathbf{1}\}}$
are the MPI ranks.  If the $\mathbf{p}_{\mathbf{k},\mathbf{r}}^{}$th
stripe globally (i.e., the first locally) has to be transferred to the
first stripe in the
$\mathbf{t}_{\mathbf{k},\mathbf{r}}^{\{\mathbf{0}\}}$th process, that
is encoded as
$-(\mathbf{t}_{\mathbf{k},\mathbf{r}}^{\{\mathbf{0}\}}+1)$, else the
second stripe of the target is encoded as
$(\mathbf{t}_{\mathbf{k},\mathbf{r}}^{\{\mathbf{0}\}}+1)$.  Similarly,
if the $\mathbf{q}_{\mathbf{k},\mathbf{r}}^{}$th stripe globally
(i.e., the second locally) has to be transferred to the first stripe
in the $\mathbf{t}_{\mathbf{k},\mathbf{r}}^{\{\mathbf{1}\}}$th
process, that is encoded as
$-(\mathbf{t}_{\mathbf{k},\mathbf{r}}^{\{\mathbf{1}\}}+1)$, else the
second stripe of the target is encoded as
$(\mathbf{t}_{\mathbf{k},\mathbf{r}}^{\{\mathbf{1}\}}+1)$.  Adding $1$
to the rank ensures that the rank $0$ can be encoded as either $1$ or
$-1$.

The number of steps in a sweep, denoted by $\mathfrak{n}'$, is
$\mathfrak{n}-1$ for \textsc{me}, and $\mathfrak{n}$ for \textsc{mm}.
The strategy mapping, once computed, can be reused for multiple runs
of the algorithm, as long as the strategy kind, $n$ (after bordering),
and $\mathbf{s}$ do not change between the runs.
%
%
\subsubsection{Algorithm initialization}\label{sss:4.1.4}
%
%
First, the CPU memory is allocated in each process, and the data is
loaded (e.g., from a file), assuming $\mathbf{k}=0$, i.e., the
$\mathbf{r}$th process contains the
$\mathbf{p}_{\mathbf{0},\mathbf{r}}^{}$th and
$\mathbf{q}_{\mathbf{0},\mathbf{r}}$th stripes of $F$, $G$, and $Z$.

Then, the device memory is allocated, if it is not already available,
and an MPI barrier is reached.  Timing of the algorithm includes
everything that occurs from this barrier on, except the optional
deallocation of the device memory.

The constant memory on each GPU is set up, and the stripes are copied
to the device (global) memory, all of which could be done
asynchronously.  The involved stream(s) are then synchronized,
depending on the way the copies have been performed.

It remains to be decided how many sweeps $S$ in
Algorithm~\ref{alg:cpu1} to allow.  As with the pointwise level, there
are two obvious choices: either some reasonably large number, e.g.,
$30$ (as in \textsc{fb}), or $1$ (as in \textsc{bo}).  Now a variant
of the multi-GPU algorithm is specified by the selected variant of the
single-GPU algorithm, with the outermost strategy and the
choice of $S$ added; e.g., ZHZ0-(\textsc{me}-\textsc{bo},
\textsc{me}-\textsc{fb}-\textsc{me}) for \textsc{me} and $S=1$,
respectively, using ZHZ0-(\textsc{me}-\textsc{fb}-\textsc{me}) at the
single-GPU level.

As shown in subsection~\ref{sss:5.3.2}, the imbalance of the
computational time each GPU requires with \textsc{fb} (i.e., one GPU
may need more sweeps in Algorithm~\ref{alg:cpu1} than another to reach
convergence within an outermost step) is significantly detrimental to
the overall performance---contrary to the single-GPU case (see
subsection~\ref{ss:3.5}).  Unlike there, where such imbalance between
the thread blocks' sweep counts is offset by a large number of thread
blocks to be scheduled on a small number of multiprocessors, here in
the multi-GPU case there is a one-to-one correspondence between the
number of tasks to perform and the number of execution units (GPUs) to
perform them, so the time required for an outermost step depends on
the slowest run of Algorithm~\ref{alg:cpu1} within it.  Thus,
\textsc{bo} is recommended here instead.
%
%
\subsection{The main part of the algorithm}\label{ss:4.2}
%
%
In the pre-iterative part of the algorithm, \texttt{initFGZ} kernel is
called (see subsection~\ref{ss:3.3}), once in each process, in the
chosen stream $\hbox{\sc s}_{\mathbf{r}}$.  It is \emph{not\/} called
again in the context of Algorithm~\ref{alg:cpu1}.  Here, the row
offset $l$ in \texttt{initFGZ} is calculated according to the logical
(not physical) index of a column, i.e.,
$l\assgn j+\mathbf{p}_{\mathbf{0},\mathbf{r}}^{}\cdot\mathbf{w}$ 
if $j<\mathbf{w}$, and
$l\assgn j-\mathbf{w}+\mathbf{q}_{\mathbf{0},\mathbf{r}}^{}\cdot\mathbf{w}$
otherwise, with $\mathbf{p}_{\mathbf{0},\mathbf{r}}^{}$ and
$\mathbf{q}_{\mathbf{0},\mathbf{r}}^{}$ sent to the kernel as
parameters.  The stripes of $F_0^{}$, $G_0^{}$, and $Z_0^{}$ are then
ready on each GPU (copying them to the CPU is not needed) for the
iterative part of the algorithm.
%
%
\subsubsection{Point-to-point communication and reductions}\label{sss:4.2.1}
%
%
Except for a single collective \texttt{MPI\_Allreduce} operation
required per an outermost sweep, all other communication in the
algorithm is of the non-blocking, point-to-point kind, occurring in
every outermost step.  The communication parts of the algorithm, from
a given process' perspective, are formalized in
Algorithms~\ref{alg:recv}, \ref{alg:send}, \ref{alg:wait}, and
\ref{alg:ared}, and put together in Algorithm~\ref{alg:mgpu}.
\begin{algorithm*}[h!btp]
  \SetKwFunction{Irecv}{MPI\_Irecv}
  $\hbox{\tt tag}\assgn 1;\quad\mathtt{i}\assgn 0$\tcp*[r]{\texttt{tag} tells which stripe from a sender has to be received}
  \ForEach(\tcp*[f]{$\mathbf{o}$ indexes the first or the second destination's stripe}){$\mathbf{o}\in\{\mathbf{0},\mathbf{1}\}$}{
    \ForEach(\tcp*[f]{$Y$ denotes the destination's host matrix}){$Y\in\{F_{\mathbf{o}}^{(\mathbf{r})},G_{\mathbf{o}}^{(\mathbf{r})},Z_{\mathbf{o}}^{(\mathbf{r})}\}$}{
      \ForEach(\tcp*[f]{$V$ refers to the real or the imaginary part}){$V\in\{\Real,\Imag\}$}{
        \Irecv{$V(Y),m_Y^{}\!\cdot\!\mathbf{w},\hbox{\tt MPI\_DOUBLE},\hbox{\tt MPI\_ANY\_SOURCE},\hbox{\tt tag},\hbox{\tt MPI\_COMM\_WORLD},\mathtt{r}[\mathtt{i}]$}\;
        $\hbox{\tt tag}\assgn\hbox{\tt tag}+1;\quad\mathtt{i}\assgn\mathtt{i}+1$\tcp*[r]{increment \texttt{tag} and $\mathtt{i}$, which indexes requests}
      }
    }
  }
  \caption{The non-blocking receives in the $\mathbf{k}$th step of the $\mathbf{r}$th process.}
  \label{alg:recv}
\end{algorithm*}
\begin{algorithm*}[h!btp]
  \SetKwFunction{Isend}{MPI\_Isend}
  \tcp{{\rm variable {\tt i} is assumed to hold the last value assigned to it in Algorithm~\ref{alg:recv} in the {\bf k}th step}}
  \ForEach(\tcp*[f]{$\mathbf{o}$ indexes the first or the second source's stripe}){$\mathbf{o}\in\{\mathbf{0},\mathbf{1}\}$}{
    $\mathtt{j}\assgn 1$\tcp*[r]{$\mathtt{j}$ is a \underline{base} tag}
    \ForEach(\tcp*[f]{$Y$ denotes the source's device matrix}){$Y\in\{F_{\mathbf{o}}^{[\mathbf{r}]},G_{\mathbf{o}}^{[\mathbf{r}]},Z_{\mathbf{o}}^{[\mathbf{r}]}\}$}{
      \ForEach(\tcp*[f]{$V$ refers to the real or the imaginary part}){$V\in\{\Real,\Imag\}$}{
        \tcp{tag: base + offset 0 or 6 (first or second stripe at destination)}
        $\hbox{\tt tag}\assgn\mathtt{j}+(\sign(\mathfrak{t}_{\mathbf{k},\mathbf{r}}^{\{\mathbf{o}\}})+1)\cdot 3$\;
        \tcp{send to destination $\mathbf{t}_{\mathbf{k},\mathbf{r}}^{\{\mathbf{o}\}}$}
        \Isend{$V(Y),m_Y^{}\cdot\mathbf{w},\hbox{\tt MPI\_DOUBLE},\mathbf{t}_{\mathbf{k},\mathbf{r}}^{\{\mathbf{o}\}},\hbox{\tt tag},\hbox{\tt MPI\_COMM\_WORLD},\mathtt{r}[\mathtt{i}]$}\;
        $\mathtt{j}\assgn\mathtt{j}+1;\quad\mathtt{i}\assgn\mathtt{i}+1$\tcp*[r]{increment $\mathtt{j}$ and $\mathtt{i}$}
      }
    }
  }
  \caption{The non-blocking sends in the $\mathbf{k}$th step of the $\mathbf{r}$th process.}
  \label{alg:send}
\end{algorithm*}
\begin{algorithm*}[h!btp]
  \SetKwFunction{Waitall}{MPI\_Waitall}
  \SetKwFunction{Barrier}{MPI\_Barrier}
  \tcp{{\rm variable {\tt i} is assumed to hold the last value assigned to it in Algorithm~\ref{alg:send} in the {\bf k}th step}}
  \Waitall{$\mathtt{i},\mathtt{r},\hbox{\tt statuses}$}\tcp*[r]{wait for all pending MPI requests to complete}
  \ForEach(\tcp*[f]{$(W,Y)$ are (host,device) matrices}){$(W,Y)\in\{(F^{(\mathbf{r})},F^{[\mathbf{r}]}),(G^{(\mathbf{r})},G^{[\mathbf{r}]}),(Z^{(\mathbf{r})},Z^{[\mathbf{r}]})\}$}{
    \ForEach(\tcp*[f]{$V$ refers to the real or the imaginary part}){$V\in\{\Real,\Imag\}$}{
      copy $V(W)$ to $V(Y)$ using \texttt{cudaMemcpy2DAsync}\tcp*[r]{in the appropriate stream(s)}
    }
  }
  synchronize the stream(s) used for copying and call \Barrier{$\hbox{\tt MPI\_COMM\_WORLD}$}\;
  \caption{Completion of the communication and the host-to-device transfers in the $\mathbf{k}$th step of the $\mathbf{r}$th process.}
  \label{alg:wait}
\end{algorithm*}
\begin{algorithm*}[h!btp]
  \SetKwFunction{AllRed}{MPI\_Allreduce}
  \SetKw{Break}{break}
  \AllRed{$\{\widehat{\mathfrak{S}}_{\mathbf{c}}^{},\widehat{\mathfrak{B}}_{\mathbf{c}}^{}\},\{\sum_{\mathbf{s}}^{}\widehat{\mathfrak{S}}_{\mathbf{c}}^{},\sum_{\mathbf{s}}^{}\widehat{\mathfrak{B}}_{\mathbf{c}}^{}\},2,\hbox{\tt MPI\_UNSIGNED\_LONG},\hbox{\tt MPI\_SUM},\hbox{\tt MPI\_COMM\_WORLD}$}\;
  \tcp{Is the sum of all per-process, per-sweep big transformation counters 0?}
  \lIf{$\sum_{\mathbf{s}}^{}\widehat{\mathfrak{B}}_{\mathbf{c}}^{}=0$}{\Break}
  \caption{Convergence criterion checking in the $\mathbf{c}$th sweep of the $\mathbf{r}$th process.}
  \label{alg:ared}
\end{algorithm*}

The first guiding principle for such a design of the communication is
to keep it as general as possible.  Any process topology (including no
topology in particular), suggested by the communication pattern of the
chosen Jacobi strategy can be accommodated with equal ease.

The second principle is to facilitate hiding the communication
overhead behind the GPU computation.  Before a call of
Algorithm~\ref{alg:cpu1} occurs within an outermost step of a given
process, the non-blocking receives to the CPU stripes are started in
anticipation of an early finish of the GPU work of the step in the
processes that are to send their transformed stripes to the process in
question.  That way, while the given GPU still computes, its CPU can
in theory start or even complete receiving one or both transformed
stripes required in the following step.  There remains an issue with
several slowly progressing processes that might keep the rest of them
idle, but at least the point-to-point data transfers can happen soon
after the data is ready, not waiting for a massive data exchange with
all processes communicating at the same time.

The third principle is to minimize the memory requirements of both the
CPUs and the GPUs by sending the transformed data from the GPU RAM of
one process to the CPU RAM of another two.  That way, no separate,
``shadow'' GPU buffers are required to receive the data.  The CPU
stripes have to be present anyhow, to load the inputs and to collect
the outputs, so they are reused as the communication buffers, with a
penalty of the additional CPU-to-GPU data transfers after the main
data exchange.

Matching a stripe to be sent from one process to a stripe that has to
be received in another process is accomplished by MPI tags annotating
the messages.  In the complex case there are twelve stripes in total
(six in the real case, without the imaginary stripes) to be received
by a process in a single outermost step (see Algorithm~\ref{alg:recv}
for their tag numbers).

When a message comes to a process, from any sender, it is only
accepted if it bears a valid tag (between $1$ and $12$, inclusive) and
the message data is stored in the corresponding stripe, as in
Algorithm~\ref{alg:recv}.  Likewise, when a stripe has to be sent, the
strategy mapping is consulted to get the destination process' rank,
and decide if the stripe should become the first or the second one at
the destination.  According to that information, the message's tag is
calculated as in Algorithm~\ref{alg:send}.
%
%
\subsubsection{The CPU part of the algorithm}\label{sss:4.2.2}
%
%
The pre-iterative, iterative, and post-iterative parts of the
algorithm are shown in Algorithm~\ref{alg:mgpu}.
\begin{algorithm*}[h!btp]
  \SetKw{Break}{break}
  \InitFGZ{$\mathbf{p}_{\mathbf{0},\mathbf{r}}^{},\mathbf{q}_{\mathbf{0},\mathbf{r}}^{}$}\tcp*[r]{compute $F_0^{[\mathbf{r}]},\;G_0^{[\mathbf{r}]},\;Z_0^{[\mathbf{r}]}$ in the stream $\hbox{\sc s}_{\mathbf{r}}$}
  \For(\tcp*[f]{outermost sweep $\mathbf{c}$}){$0\le\mathbf{c}<30$}{
    $\{\widehat{\mathfrak{S}}_{\mathbf{c}}^{},\widehat{\mathfrak{B}}_{\mathbf{c}}^{}\}\assgn\{0,0\}$\tcp*[r]{reset the per-process, per-sweep transformation counters}
    \For(\tcp*[f]{outermost step $\mathbf{k}$}){$0\le\mathbf{k}<\mathfrak{n}'$}{
      start receiving into $F^{(\mathbf{r})},G^{(\mathbf{r})},Z^{(\mathbf{r})}$ as in Algorithm~\ref{alg:recv}\;
      call the single-GPU Algorithm~\ref{alg:cpu1} with $\hbox{\sc s}=\hbox{\sc s}_{\mathbf{r}}$ on $F^{[\mathbf{r}]},G^{[\mathbf{r}]},Z^{[\mathbf{r}]}$ with the chosen $S$\;
      \tcp{increment the transformation counters by those from Algorithm~\ref{alg:cpu1}}
      $\widehat{\mathfrak{S}}_{\mathbf{c}}^{}\assgn\widehat{\mathfrak{S}}_{\mathbf{c}}^{}+\widetilde{\mathfrak{S}};\quad\widehat{\mathfrak{B}}_{\mathbf{c}}^{}\assgn\widehat{\mathfrak{B}}_{\mathbf{c}}^{}+\widetilde{\mathfrak{B}}$\;
      start sending the transformed $F^{[\mathbf{r}]},G^{[\mathbf{r}]},Z^{[\mathbf{r}]}$ as in Algorithm~\ref{alg:send}\;
      complete the communication and copy the received $F^{(\mathbf{r})},G^{(\mathbf{r})},Z^{(\mathbf{r})}$ to $F^{[\mathbf{r}]},G^{[\mathbf{r}]},Z^{[\mathbf{r}]}$, as in Algorithm~\ref{alg:wait}\;
    }
    reduce the transformation counters across the communicator\;
    \Break if the convergence has been reached, as in Algorithm~\ref{alg:ared}\;
  }
  \Rescale{$\hbox{\tt true}$};$\quad$\StreamSync{$\hbox{\sc s}_{\mathbf{r}}$}\tcp*[r]{full rescaling of $Z^{[\mathtt{r}]}$ in $\hbox{\sc s}_{\mathbf{r}}$}
  optionally, copy $F^{[\mathbf{r}]},G^{[\mathbf{r}]},Z^{[\mathbf{r}]},\Sigma_F^{[\mathbf{r}]},\Sigma_G^{[\mathbf{r}]},\Sigma^{[\mathbf{r}]}$ to $F^{(\mathbf{r})},G^{(\mathbf{r})},Z^{(\mathbf{r})},\Sigma_F^{(\mathbf{r})},\Sigma_G^{(\mathbf{r})},\Sigma^{(\mathbf{r})}$ and synchronize the stream(s)\;
  \Barrier{$\hbox{\tt MPI\_COMM\_WORLD}$}\tcp*[r]{completion of the algorithm and its timing}
  \caption{The CPU part of the multi-GPU implicit Hari--Zimmermann algorithm (for the $\mathbf{r}$th process).}
  \label{alg:mgpu}
\end{algorithm*}
The final full rescaling with the extraction of the generalized
singular values happens only once (i.e., not in the context of
Algorithm~\ref{alg:cpu1}).  As the convergence criterion relies on
sum-reducing the per-sweep counters of the big transformations applied
in all processes, an implicit synchronization point at the end of a
sweep is introduced.
%
%
\section{Numerical testing}\label{s:5}
%
%
The purpose of the numerical testing of the single-GPU and the
multi-GPU algorithms is twofold.  First, it has been meant to compare
the variants of the algorithms in terms of performance and accuracy
and discover which (if any) variant stands out as the best one in
either aspect.  Second, it should inform the potential users about the
algorithms' behavior on two sets of realistic, small and
medium-to-large sized problems.

By performance it is meant the wall execution time.  Counting FLOPS
(floating-point operations per second) rate makes less sense here than
in the algorithms (such as the matrix multiplication) that solely
depend on a subset of the arithmetic operations of a similar execution
complexity, such as additions, subtractions, multiplications, and
FMAs.  Instead, the algorithms presented here necessarily involve a
substantial amount of divisions and (reciprocal) square roots.
Moreover, the majority of performance gains compared to a simple,
pointwise algorithm come from a careful usage of the fast shared
memory and the GPU registers, as it is also shown
in~\cite{Novakovic-2015}, and not from tweaking the arithmetic
intensity.  The wall time should therefore be more informative than
FLOPS about the expected behavior of the algorithm on present-day
hardware, and about the differences in the algorithm's variants, since
the future performance is very hard to predict without a complex
model that takes into account all levels of the memory hierarchy, not
only the arithmetic operations and the amount of parallelism
available.

Accuracy of the algorithm can be assessed in several ways.  In both
the real and the complex case the relative normwise errors of the
decompositions of $F$ and $G$,
\begin{displaymath}
  \|F - U \Sigma_F X\|_F / \|F\|_F,\quad \|G - V \Sigma_G X\|_F / \|G\|_F,
\end{displaymath}
were computed the same way as
in~\cite{Singer-DiNapoli-Novakovic-Caklovic-2020}.  Namely, $X$ had to
be explicitly obtained as $Z^{-1}$ by solving the linear system
$ZX=I$.  First, the LU factorization of $Z$ with complete pivoting was
performed by the LAPACK routine \texttt{DGETC2} (or \texttt{ZGETC2}),
followed by the system solving using the routine \texttt{DGESC2} (or
\texttt{ZGESC2}).

The ensuing matrix multiplications and the Frobenius norm computations
were using Intel 80-bit hardware-supported extended precision
(\texttt{REAL(KIND=10)} in GNU Fortran), to reduce the effects of the
rounding errors on the final result while avoiding the expensive,
emulated quadruple (128-bit) precision.

The numerical orthogonality of the left generalized singular vectors
$U$ and $V$ was computed in the extended precision as
$\|U^{\ast}U-I\|_F^{}$ and $\|V^{\ast}V-I\|_F^{}$, respectively.

When the (almost) exact generalized singular values $\Sigma$ are known
in advance, as is the case with the small real dataset (see
subsection~\ref{sss:5.1.2}), the maximal relative error in the
computed $\widehat{\Sigma}$ can be obtained as
\begin{displaymath}
  \max_{1\le i\le n}|(\sigma_i-\hat{\sigma}_i)/\sigma_i|.
\end{displaymath}
%
%
\subsection{Testing environment and data}\label{ss:5.1}
%
%
The testing environment was the same for all tests, as described in
subsection~\ref{sss:5.1.1}.  Apart from the GPU compute architecture
7.0, some tests have been repeated on a Maxwell GPU (GeForce GTX TITAN
X, architecture 5.2) and a Kepler GPU (GeForce GT 730, architecture
3.5), to verify the portability of the code and the numerical
reproducibility of the results.  Also, a few sample runs of the
multi-GPU algorithm on a small matrix have been tried on a combination
of those two GPUs, with the code built for both architectures, to
ensure that the algorithm functions correctly in such a heterogeneous
environment.

The testing data is synthetic (not from any application domain) and is
described in subsection~\ref{sss:5.1.2}.  Please see the supplementary
material for its availability.
%
%
\subsubsection{Testing environment}\label{sss:5.1.1}
%
%
The testing environment comprises two Intel Xeon Silver 4114 CPUs,
$384\,\mathrm{GiB}$ of RAM, and four NVIDIA Tesla V100-SXM2-16GB
(Volta) GPUs per node, with a 64-bit Linux (CentOS 7.5.1804), the GCC
4.8.5 C++ compiler, CUDA 10.0, and a build of Open MPI 3.0.0
distribution with the CUDA support.
%
%
\subsubsection{Testing data}\label{sss:5.1.2}
%
%
Two datasets have been generated: a ``small'' and a ``large'' one,
with their names referring to the orders of the square matrices
forming the pairs contained in them.  The small dataset contains both
the real and the complex matrix pairs, with each matrix stored in (and
then read from) its unformatted binary file, while the large dataset
contains only the complex matrix pairs.

The small dataset has 19 matrix pairs for each datatype, with the
orders of the matrices ranging from $512$ to $9728$ in steps of $512$.
The large dataset has 3 matrix pairs, with the orders of the matrices
being $18\cdot 1024=18432$, $24\cdot 1024=24576$, and
$36\cdot 1024=36864$, so that the GPU RAM requirements do not exceed
the memory provided by one, two, and four GPUs, respectively.  No
matrices in either dataset require bordering.

The real matrix pairs in the small dataset were generated in quadruple
datatype (\texttt{REAL(KIND=16)} in Intel Fortran) and rounded to
double precision datatype.  The same test generation method was
employed as in~\cite{Novakovic-Singer-Singer-2015}.  The required BLAS
and LAPACK routines had been adapted as required.  The core of the
generation are two quadruple-adapted LAPACK testing routines:
\texttt{xLAGSY}, that generates a pseudorandom symmetric matrix, here
of the full bandwidth, from a given diagonal prescribing the
eigenvalues of the matrix; and \texttt{xLAROR}, that here multiplies a
given matrix from the left by a pseudorandom orthogonal matrix.  The
diagonals of $\Sigma_F^{}$, $\Sigma_G^{}$, and $\Lambda_X^{}$ were
generated by \texttt{DLARND}, a standard LAPACK's pseudorandom number
generator, here with the uniform probability distribution on $(0,1)$,
such that only those values returned by it that had been greater than
$10^{-10}$ were accepted.  Then, $U\Sigma_F^{}$ was generated from
$\Sigma_F^{}$, and $V\Sigma_G^{}$ from $\Sigma_G^{}$, both using
\texttt{xLAROR}, while $X$ was obtained from $\Lambda_X^{}$ using
\texttt{xLAGSY}.  After that, $F\assgn U\Sigma_F^{}\cdot X$ and
$G\assgn V\Sigma_G^{}\cdot X$ (in quadruple).  The generator finishes
with a call to the preprocessing part, \texttt{DGGSVP3}, of the
LAPACK's GSVD method (see~\cite{Anderson-et-al-99} and the routine's
comments), to make the data usable for comparison with \texttt{DTGSJA}
Kogbetliantz-type GSVD routine from LAPACK\@.  On the small dataset
the relative errors in the generalized singular values computed on the
CPU by LAPACK and on a single GPU by the proposed algorithm were
compared.

The complex matrix pairs in both datasets were generated by a much
simpler procedure, described
in~\cite{Singer-DiNapoli-Novakovic-Caklovic-2020}.  Namely, each
matrix in a pair was generated by a call to \texttt{ZLATMS} LAPACK
testing routine as Hermitian and positive definite, with its
pseudorandom eigenvalues uniformly distributed in $(0,1)$.
%
%
\subsection{Results with the single-GPU algorithm}\label{ss:5.2}
%
%
When presenting the performance of several variants, a decision has
been made to show the execution time plots and tables for the fastest
variant on the largest matrix pair in the small dataset, separately in
the real and in the complex case.  To compare those results with the
ones from other variants, a useful measure is the relative slowdown on
a given matrix order.  Fixing the reference variant $r$, as suggested
by the results in the plots, for another variant $v$ and a matrix
order $n$ let $T_n^{(r)}$ and $T_n^{(v)}$ be the execution times of
$r$ and $v$ on a matrix pairs with the matrices of order $n$,
respectively.  The relative slowdown $S_n^{(v:r)}$ of $v$ compared to
$r$ on $n$, given in percentages of the execution time of $r$, is
\begin{displaymath}
  S_n^{(v:r)}\assgn\frac{T_n^{(v)}-T_n^{(r)}}{T_n^{(r)}}\cdot 100.
\end{displaymath}
The average relative slowdown across the entire dataset (with 19
matrix pairs) is given as
$S_{\text{avg}}^{(v:r)}\assgn\sum_n^{}S_n^{(v:r)}/19$.
%
%
\subsubsection{Performance in the real case}\label{sss:5.2.1}
\begin{figure*}[h!btp]
  \centering
  \includegraphics[keepaspectratio,width=.4975\textwidth]{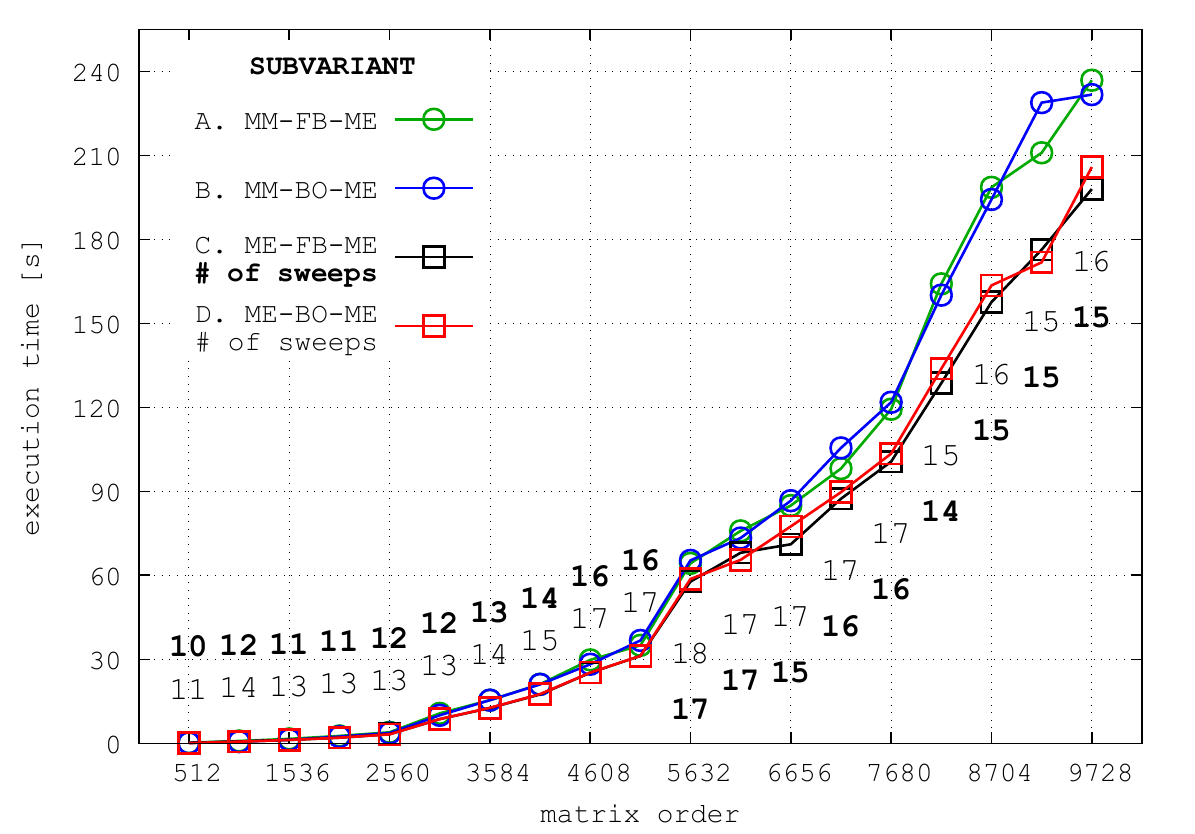}
  \includegraphics[keepaspectratio,width=.4975\textwidth]{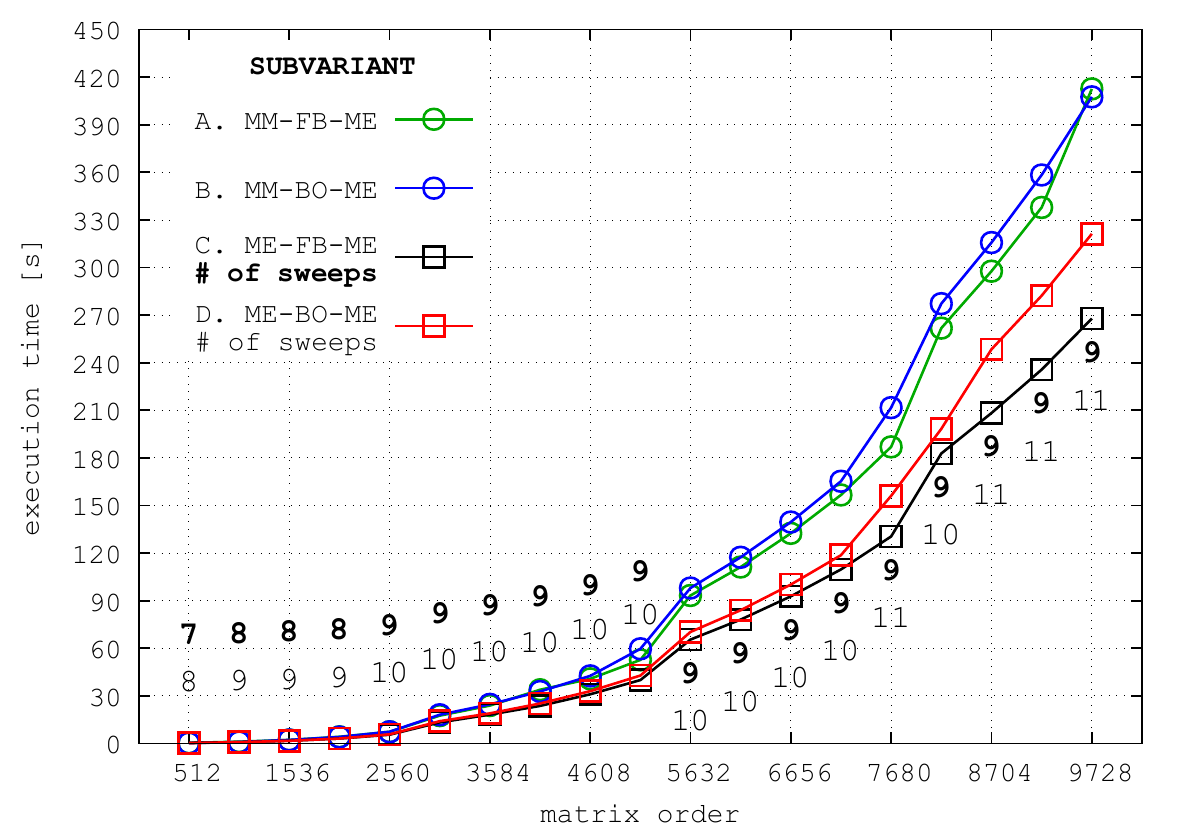}
  \caption{The wall execution time of four subvariants of DHZ0 (left)
    and ZHZ4 (right) and the number of outer sweeps for two fastest
    subvariants on the small real (left) and complex (right)
    datasets.}
  \label{fig:HZ-ts}
\end{figure*}

%
%
In the left subfigure of Figure~\ref{fig:HZ-ts} the wall execution
time of four subvariants of DHZ0 (the fastest of eight variants from
Table~\ref{tbl:var} on the largest matrix pair) and the number of
outer sweeps for two fastest subvariants on the small real dataset are
shown.
Table~\ref{tbl:QRvsCholesky} in Appendix~\ref{s:E} contains the wall
time for DHZ0-(\textsc{me}-\textsc{fb}-\textsc{me}), which is
generally the fastest of the four subvariants.

In Table~\ref{tbl:Dslow} the intervals of relative slowdown of other
real single-GPU (\textsc{me}-\textsc{fb}-\textsc{me}) variants
compared to DHZ0-(\textsc{me}-\textsc{fb}-\textsc{me}) (the reference
variant) are given.
\begin{table}[h!bt]
  \caption{The intervals of relative slowdown of other real single-GPU
    variants compared to DHZ0, all of
    (\textsc{me}-\textsc{fb}-\textsc{me}) subvariant.  A negative
    slowdown is a speedup.}
  \label{tbl:Dslow}
  \centering
  \addtolength{\tabcolsep}{-2pt}
  \begin{tabular}{@{}cccc@{}}
    \toprule
    \multirow{2}{*}{ID} & maximal relative & minimal relative & average relative\\
    & slowdown~[\%] & slowdown~[\%] & slowdown~[\%]\\
    \midrule
    1 & 11.107853 & \hphantom{-}0.984035 & 3.093490\\
    2 & 17.137007 & \hphantom{-}1.722547 & 4.870879\\
    3 & 32.907217 & \hphantom{-}3.144621 & 9.204879\\
    4 & \hphantom{0}0.328320 & -0.678588 & 0.139650\\
    5 & 11.353652 & \hphantom{-}1.395604 & 3.373963\\
    6 & 17.962363 & \hphantom{-}2.065662 & 5.199807\\
    7 & 32.290521 & \hphantom{-}3.511386 & 9.413653\\
    \bottomrule
  \end{tabular}
\end{table}
%
%
\subsubsection{Performance in the complex case}\label{sss:5.2.2}
%
%
In the right subfigure of Figure~\ref{fig:HZ-ts} the wall execution
time of four subvariants of ZHZ4 (the fastest of eight variants from
Table~\ref{tbl:var} on the largest matrix pair) and the number of
outer sweeps for two fastest subvariants on the small complex dataset
are shown.

In Table~\ref{tbl:Zslow} the intervals of relative slowdown of other
complex single-GPU (\textsc{me}-\textsc{fb}-\textsc{me}) variants
compared to ZHZ4-(\textsc{me}-\textsc{fb}-\textsc{me}) (the reference
variant) are given.  Since ZHZ0-(\textsc{me}-\textsc{fb}-\textsc{me})
is the fastest variant \emph{on average\/}, differing from ZHZ4 only
subtly (in the convergence criterion), ZHZ0 is used instead of ZHZ4 in
the multi-GPU algorithm and in subsection~\ref{sss:5.2.6}.
Table~\ref{tbl:GPUvsCPU} in Appendix~\ref{s:F} shows the wall time for
ZHZ0-(\textsc{me}-\textsc{fb}-\textsc{me}).
\begin{table}[h!bt]
  \caption{The intervals of relative slowdown of other complex
    single-GPU variants compared to ZHZ4, all of
    (\textsc{me}-\textsc{fb}-\textsc{me}) subvariant.  A negative
    slowdown is a speedup; e.g., ZHZ0 is faster than ZHZ4 on average,
    even if it is slower sometimes.}
  \label{tbl:Zslow}
  \centering\addtolength{\tabcolsep}{-2pt}
  \begin{tabular}{@{}cccc@{}}
    \toprule
    \multirow{2}{*}{ID} & maximal relative & minimal relative & average relative\\
    & slowdown~[\%] & slowdown~[\%] & slowdown~[\%]\\
    \midrule
    0 & \hphantom{0}0.488301 &            -6.403801 &            -0.493666\\
    1 & \hphantom{0}2.309416 &            -0.241201 & \hphantom{-}0.841493\\
    2 & \hphantom{0}7.135216 &            -1.206508 & \hphantom{-}1.756809\\
    3 &            12.632933 & \hphantom{-}0.866139 & \hphantom{-}3.167284\\
    5 & \hphantom{0}1.366689 &            -0.442658 & \hphantom{-}0.312471\\
    6 & \hphantom{0}8.423992 & \hphantom{-}0.720393 & \hphantom{-}2.287251\\
    7 &            11.257148 & \hphantom{-}1.092003 & \hphantom{-}3.240513\\
    \bottomrule
  \end{tabular}
\end{table}
%
%
\subsubsection{Detailed timings of the main kernel's subphases}\label{sss:5.2.3}
%
%
Tables~\ref{tbl:Dsubph} and~\ref{tbl:Zsubph} show the percentages of
time each subphase of the main kernel takes on a Volta GPU for the
selected subvariants of the real and the complex single-GPU algorithm,
respectively.  As the problem size increases, the first and the fourth
subphases (i.e., the matrix multiplications) start to dominate over
the others.  With a suitable modification for the recent GPU
architectures like Ampere, the Hari--Zimmermann algorithm could
therefore benefit from their dedicated hardware and instructions for
speeding up the \texttt{GEMM}-like operations.
\begin{table}[h!bt]
  \caption{Percentage of time (rounded to the nearest per mil) spent
    in the subphases~\textbf{1} to \textbf{4} in all invocations of
    the main kernel for the subvariants $\mathtt{C}_0^{\mathbb{R}}$,
    i.e., DHZ0-(\textsc{me}-\textsc{fb}-\textsc{me}) and
    $\mathtt{D}_0^{\mathbb{R}}$, i.e.,
    DHZ0-(\textsc{me}-\textsc{bo}-\textsc{me}), on the small real
    dataset.}
  \label{tbl:Dsubph}
  \centering\addtolength{\tabcolsep}{-1.25pt}
  \begin{tabular}{@{}c|cccc|cccc@{}}
    \toprule
    \multirow{2}{*}{$n$} & \multicolumn{4}{c}{subphases of $\mathtt{C}_0^{\mathbb{R}}$ [\%]} & \multicolumn{4}{c}{subphases of $\mathtt{D}_0^{\mathbb{R}}$ [\%]}\\
     & \textbf{1} & \textbf{2} & \textbf{3} & \textbf{4} & \textbf{1} & \textbf{2} & \textbf{3} & \textbf{4}\\
    \midrule
    \hphantom{0}512 & 12.3 &            10.0 &            53.8 & 23.9 & 20.9 &            17.0 &            21.3 & 40.7\\
               1024 & 17.5 & \hphantom{0}7.1 &            41.1 & 34.4 & 25.7 &            10.4 &            13.1 & 50.8\\
               1536 & 20.6 & \hphantom{0}5.3 &            34.1 & 40.0 & 28.2 & \hphantom{0}7.3 & \hphantom{0}9.3 & 55.2\\
               2048 & 22.8 & \hphantom{0}4.4 &            28.2 & 44.6 & 29.4 & \hphantom{0}5.7 & \hphantom{0}7.3 & 57.6\\
               2560 & 24.9 & \hphantom{0}3.7 &            23.3 & 48.0 & 30.6 & \hphantom{0}4.6 & \hphantom{0}5.9 & 59.0\\
               3072 & 26.1 & \hphantom{0}3.3 &            20.5 & 50.1 & 31.1 & \hphantom{0}3.9 & \hphantom{0}5.0 & 59.9\\
               3584 & 26.5 & \hphantom{0}2.9 &            18.7 & 51.8 & 31.2 & \hphantom{0}3.4 & \hphantom{0}4.4 & 61.0\\
               4096 & 27.3 & \hphantom{0}2.6 &            17.1 & 53.0 & 31.6 & \hphantom{0}3.0 & \hphantom{0}3.9 & 61.4\\
               4608 & 27.8 & \hphantom{0}2.4 &            15.4 & 54.4 & 31.7 & \hphantom{0}2.7 & \hphantom{0}3.5 & 62.1\\
               5120 & 28.2 & \hphantom{0}2.2 &            14.5 & 55.2 & 32.0 & \hphantom{0}2.5 & \hphantom{0}3.2 & 62.4\\
               5632 & 28.8 & \hphantom{0}2.0 &            13.2 & 56.0 & 32.2 & \hphantom{0}2.2 & \hphantom{0}2.9 & 62.7\\
               6144 & 29.2 & \hphantom{0}1.9 &            12.0 & 56.9 & 32.2 & \hphantom{0}2.1 & \hphantom{0}2.7 & 63.0\\
               6656 & 29.4 & \hphantom{0}1.8 &            11.7 & 57.2 & 32.4 & \hphantom{0}1.9 & \hphantom{0}2.5 & 63.2\\
               7168 & 29.7 & \hphantom{0}1.7 &            10.8 & 57.9 & 32.5 & \hphantom{0}1.8 & \hphantom{0}2.3 & 63.4\\
               7680 & 30.0 & \hphantom{0}1.7 &            10.1 & 58.3 & 32.5 & \hphantom{0}1.8 & \hphantom{0}2.1 & 63.6\\
               8192 & 29.5 & \hphantom{0}1.6 & \hphantom{0}9.9 & 59.0 & 32.0 & \hphantom{0}1.6 & \hphantom{0}2.1 & 64.3\\
               8704 & 30.3 & \hphantom{0}1.5 & \hphantom{0}8.9 & 59.2 & 32.6 & \hphantom{0}1.7 & \hphantom{0}1.9 & 63.8\\
               9216 & 30.4 & \hphantom{0}1.4 & \hphantom{0}8.7 & 59.5 & 32.7 & \hphantom{0}1.5 & \hphantom{0}1.9 & 64.0\\
               9728 & 30.6 & \hphantom{0}1.4 & \hphantom{0}8.2 & 59.9 & 32.7 & \hphantom{0}1.5 & \hphantom{0}1.7 & 64.1\\
    \bottomrule
  \end{tabular}
\end{table}
\begin{table}[h!bt]
  \caption{Percentage of time (rounded to the nearest per mil) spent
    in the subphases~\textbf{1} to \textbf{4} in all invocations of
    the main kernel for the subvariants $\mathtt{C}_0^{\mathbb{C}}$,
    i.e., ZHZ0-(\textsc{me}-\textsc{fb}-\textsc{me}) and
    $\mathtt{D}_0^{\mathbb{C}}$, i.e.,
    ZHZ0-(\textsc{me}-\textsc{bo}-\textsc{me}), on the small complex
    dataset.}
  \label{tbl:Zsubph}
  \centering\addtolength{\tabcolsep}{-1pt}
  \begin{tabular}{@{}c|cccc|cccc@{}}
    \toprule
    \multirow{2}{*}{$n$} & \multicolumn{4}{c}{subphases of $\mathtt{C}_0^{\mathbb{C}}$ [\%]} & \multicolumn{4}{c}{subphases of $\mathtt{D}_0^{\mathbb{C}}$ [\%]}\\
     & \textbf{1} & \textbf{2} & \textbf{3} & \textbf{4} & \textbf{1} & \textbf{2} & \textbf{3} & \textbf{4}\\
    \midrule
    \hphantom{0}512 & 13.9 & 6.5 &            51.6 & 28.1 & 23.2 &            10.9 &            18.7 & 47.2\\
               1024 & 19.6 & 4.7 &            35.8 & 39.9 & 27.2 & \hphantom{0}6.5 &            10.8 & 55.5\\
               1536 & 21.9 & 3.5 &            29.3 & 45.4 & 28.5 & \hphantom{0}4.6 & \hphantom{0}7.7 & 59.3\\
               2048 & 24.7 & 2.8 &            23.3 & 49.2 & 30.2 & \hphantom{0}3.5 & \hphantom{0}6.0 & 60.3\\
               2560 & 27.3 & 2.4 &            19.2 & 51.1 & 32.1 & \hphantom{0}2.8 & \hphantom{0}4.6 & 60.5\\
               3072 & 28.4 & 2.1 &            16.3 & 53.2 & 32.5 & \hphantom{0}2.4 & \hphantom{0}3.9 & 61.2\\
               3584 & 29.1 & 1.9 &            14.3 & 54.7 & 32.8 & \hphantom{0}2.1 & \hphantom{0}3.4 & 61.8\\
               4096 & 29.5 & 1.7 &            12.9 & 55.9 & 32.9 & \hphantom{0}1.8 & \hphantom{0}3.0 & 62.3\\
               4608 & 29.9 & 1.5 &            12.0 & 56.5 & 33.0 & \hphantom{0}1.7 & \hphantom{0}2.7 & 62.6\\
               5120 & 30.1 & 1.4 &            11.2 & 57.3 & 33.1 & \hphantom{0}1.5 & \hphantom{0}2.4 & 63.0\\
               5632 & 30.6 & 1.3 &            10.5 & 57.5 & 33.5 & \hphantom{0}1.4 & \hphantom{0}2.2 & 62.9\\
               6144 & 30.8 & 1.2 & \hphantom{0}9.7 & 58.3 & 33.4 & \hphantom{0}1.3 & \hphantom{0}2.1 & 63.2\\
               6656 & 31.1 & 1.1 & \hphantom{0}9.2 & 58.6 & 33.6 & \hphantom{0}1.2 & \hphantom{0}1.9 & 63.3\\
               7168 & 31.3 & 1.1 & \hphantom{0}8.5 & 59.2 & 33.6 & \hphantom{0}1.1 & \hphantom{0}1.8 & 63.5\\
               7680 & 31.3 & 1.0 & \hphantom{0}8.0 & 59.6 & 33.5 & \hphantom{0}1.0 & \hphantom{0}1.6 & 63.8\\
               8192 & 32.2 & 1.0 & \hphantom{0}7.4 & 59.4 & 34.2 & \hphantom{0}1.0 & \hphantom{0}1.6 & 63.3\\
               8704 & 31.7 & 0.9 & \hphantom{0}7.3 & 60.1 & 33.7 & \hphantom{0}1.0 & \hphantom{0}1.5 & 63.9\\
               9216 & 31.7 & 0.9 & \hphantom{0}7.1 & 60.3 & 33.7 & \hphantom{0}0.9 & \hphantom{0}1.4 & 64.0\\
               9728 & 32.0 & 0.8 & \hphantom{0}6.7 & 60.5 & 33.8 & \hphantom{0}0.9 & \hphantom{0}1.3 & 64.0\\
    \bottomrule
  \end{tabular}
\end{table}
%
%
\subsubsection{Intensities of the floating-point arithmetic operations}\label{sss:5.2.4}
%
%
Tables~\ref{tbl:Dops} and~\ref{tbl:Zops} contain the relative
intensity of floating-point operations with rounding across all
invocations of all kernels (not only the main one, though other
kernels' contributions are negligible) on a Maxwell GPU\@.  The
columns names correspond to the double precision CUDA arithmetic
intrinsics, while $\diamond$ represents the calls to $\mathtt{hypot}$
and $\mathtt{rsqrt\_rn}$ functions (their constituent operations were
not counted separately).  The results strongly correlate with those
from Tables~\ref{tbl:Dsubph} and~\ref{tbl:Zsubph}; namely, the
intensity of $\mathtt{fma}$ increases as the matrix multiplications
take the larger portions of the overall time.  For the smaller inputs
the combined amount of divisions/reciprocals, square roots, and the
$\diamond$ function calls is in single percents, but according
to~\cite[Table~IV]{Arafa-et-al-2019} even such an amount has a
considerable influence on the execution time; e.g., a division has on
average close to 20 times the latency of a simple instruction (like
$\mathtt{fma}$) on the Volta GPUs.
\begin{table}[h!bt]
  \caption{Percentage of the total number of the floating-point
    operations with rounding performed in the invocations of all
    kernels (rounded to the nearest per myriad) for the subvariant
    DHZ0-(\textsc{me}-\textsc{fb}-\textsc{me}) on the small real
    dataset.}
  \label{tbl:Dops}
  \centering\addtolength{\tabcolsep}{-2.5pt}
  \begin{tabular}{@{}ccccccccc@{}}
    \toprule
    $n$ & $\diamond$ & $\mathtt{add}$ & $\mathtt{sub}$ & $\mathtt{mul}$ & $\mathtt{fma}$ & $\mathtt{div}$ & $\mathtt{rcp}$ & $\mathtt{sqrt}$\\
    \midrule
    \hphantom{0}512 & 0.83 & 22.68 & 1.32 & 16.33 & 54.20 & 1.71 & 0.57 & 2.35\\
               1024 & 0.55 & 15.86 & 0.91 & 11.18 & 68.29 & 1.17 & 0.39 & 1.64\\
               1536 & 0.43 & 12.34 & 0.71 & \hphantom{0}8.70 & 75.33 & 0.91 & 0.30 & 1.27\\
               2048 & 0.34 & \hphantom{0}9.99 & 0.57 & \hphantom{0}7.01 & 80.06 & 0.74 & 0.25 & 1.03\\
               2560 & 0.27 & \hphantom{0}8.00 & 0.45 & \hphantom{0}5.53 & 84.16 & 0.58 & 0.19 & 0.82\\
               3072 & 0.23 & \hphantom{0}6.86 & 0.38 & \hphantom{0}4.72 & 86.44 & 0.49 & 0.16 & 0.71\\
               3584 & 0.20 & \hphantom{0}6.16 & 0.35 & \hphantom{0}4.25 & 87.80 & 0.45 & 0.15 & 0.63\\
               4096 & 0.17 & \hphantom{0}5.43 & 0.31 & \hphantom{0}3.71 & 89.31 & 0.39 & 0.13 & 0.56\\
               4608 & 0.16 & \hphantom{0}4.89 & 0.28 & \hphantom{0}3.34 & 90.36 & 0.35 & 0.12 & 0.50\\
               5120 & 0.14 & \hphantom{0}4.51 & 0.26 & \hphantom{0}3.09 & 91.10 & 0.33 & 0.11 & 0.46\\
               5632 & 0.12 & \hphantom{0}3.94 & 0.22 & \hphantom{0}2.67 & 92.27 & 0.28 & 0.09 & 0.40\\
               6144 & 0.11 & \hphantom{0}3.62 & 0.20 & \hphantom{0}2.45 & 92.89 & 0.26 & 0.09 & 0.37\\
               6656 & 0.11 & \hphantom{0}3.52 & 0.20 & \hphantom{0}2.41 & 93.06 & 0.26 & 0.09 & 0.36\\
               7168 & 0.10 & \hphantom{0}3.13 & 0.18 & \hphantom{0}2.13 & 93.85 & 0.22 & 0.07 & 0.32\\
               7680 & 0.09 & \hphantom{0}2.91 & 0.16 & \hphantom{0}1.96 & 94.30 & 0.21 & 0.07 & 0.30\\
               8192 & 0.09 & \hphantom{0}2.85 & 0.16 & \hphantom{0}1.94 & 94.39 & 0.21 & 0.07 & 0.29\\
               8704 & 0.08 & \hphantom{0}2.59 & 0.15 & \hphantom{0}1.75 & 94.92 & 0.18 & 0.06 & 0.27\\
               9216 & 0.08 & \hphantom{0}2.48 & 0.14 & \hphantom{0}1.68 & 95.13 & 0.18 & 0.06 & 0.25\\
               9728 & 0.07 & \hphantom{0}2.34 & 0.13 & \hphantom{0}1.58 & 95.41 & 0.17 & 0.06 & 0.24\\
    \bottomrule
  \end{tabular}
\end{table}
\begin{table}[h!bt]
  \caption{Percentage of the total number of the floating-point
    operations with rounding performed in the invocations of all
    kernels (rounded to the nearest per myriad) for the subvariant
    ZHZ0-(\textsc{me}-\textsc{fb}-\textsc{me}) on the small complex
    dataset.}
  \label{tbl:Zops}
  \centering\addtolength{\tabcolsep}{-2.5pt}
  \begin{tabular}{@{}ccccccccc@{}}
    \toprule
    $n$ & $\diamond$ & $\mathtt{add}$ & $\mathtt{sub}$ & $\mathtt{mul}$ & $\mathtt{fma}$ & $\mathtt{div}$ & $\mathtt{rcp}$ & $\mathtt{sqrt}$\\
    \midrule
    \hphantom{0}512 & 0.92 & 12.81 & 0.42 & 11.63 & 71.86 & 0.42 & 0.83 & 1.12\\
               1024 & 0.53 & \hphantom{0}7.46 & 0.24 & \hphantom{0}6.63 & 83.78 & 0.24 & 0.47 & 0.65\\
               1536 & 0.39 & \hphantom{0}5.42 & 0.17 & \hphantom{0}4.83 & 88.21 & 0.17 & 0.34 & 0.47\\
               2048 & 0.30 & \hphantom{0}4.16 & 0.13 & \hphantom{0}3.68 & 90.98 & 0.13 & 0.26 & 0.36\\
               2560 & 0.23 & \hphantom{0}3.27 & 0.10 & \hphantom{0}2.84 & 92.97 & 0.10 & 0.20 & 0.28\\
               3072 & 0.19 & \hphantom{0}2.70 & 0.08 & \hphantom{0}2.33 & 94.21 & 0.08 & 0.16 & 0.23\\
               3584 & 0.17 & \hphantom{0}2.34 & 0.07 & \hphantom{0}2.02 & 94.98 & 0.07 & 0.14 & 0.20\\
               4096 & 0.15 & \hphantom{0}2.03 & 0.06 & \hphantom{0}1.74 & 95.66 & 0.06 & 0.12 & 0.17\\
               4608 & 0.13 & \hphantom{0}1.88 & 0.06 & \hphantom{0}1.61 & 95.99 & 0.06 & 0.11 & 0.16\\
               5120 & 0.12 & \hphantom{0}1.71 & 0.05 & \hphantom{0}1.46 & 96.35 & 0.05 & 0.10 & 0.15\\
               5632 & 0.11 & \hphantom{0}1.55 & 0.05 & \hphantom{0}1.32 & 96.70 & 0.05 & 0.09 & 0.13\\
               6144 & 0.10 & \hphantom{0}1.40 & 0.04 & \hphantom{0}1.19 & 97.02 & 0.04 & 0.08 & 0.12\\
               6656 & 0.09 & \hphantom{0}1.32 & 0.04 & \hphantom{0}1.12 & 97.20 & 0.04 & 0.08 & 0.11\\
               7168 & 0.09 & \hphantom{0}1.20 & 0.04 & \hphantom{0}1.03 & 97.44 & 0.04 & 0.07 & 0.10\\
               7680 & 0.08 & \hphantom{0}1.13 & 0.03 & \hphantom{0}0.96 & 97.60 & 0.03 & 0.07 & 0.10\\
               8192 & 0.07 & \hphantom{0}1.03 & 0.03 & \hphantom{0}0.88 & 97.81 & 0.03 & 0.06 & 0.09\\
               8704 & 0.07 & \hphantom{0}1.01 & 0.03 & \hphantom{0}0.86 & 97.86 & 0.03 & 0.06 & 0.09\\
               9216 & 0.07 & \hphantom{0}0.96 & 0.03 & \hphantom{0}0.81 & 97.96 & 0.03 & 0.06 & 0.08\\
               9728 & 0.07 & \hphantom{0}0.91 & 0.03 & \hphantom{0}0.78 & 98.06 & 0.03 & 0.05 & 0.08\\
    \bottomrule
  \end{tabular}
\end{table}
%
%
\subsubsection{Accuracy in the real case}\label{sss:5.2.5}
\begin{figure*}[h!btp]
  \centering
  \includegraphics[keepaspectratio,width=.4975\textwidth]{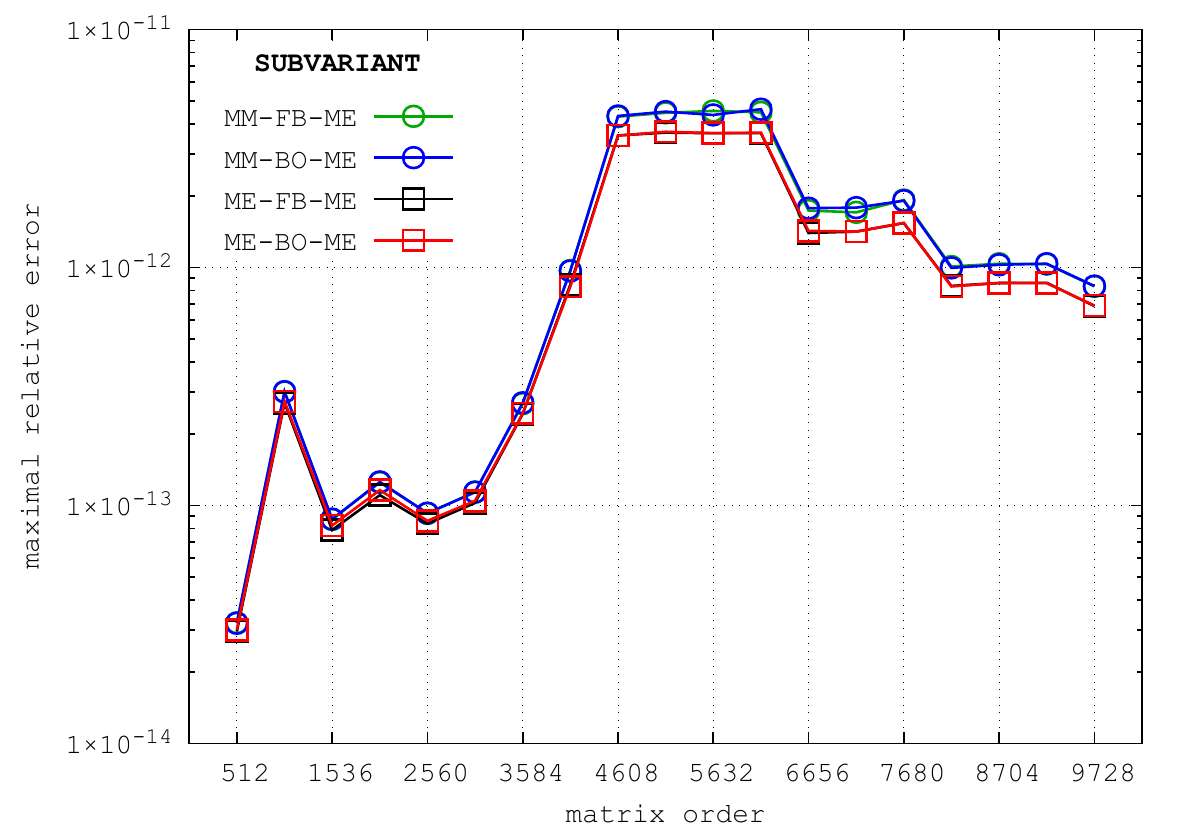}
  \includegraphics[keepaspectratio,width=.4975\textwidth]{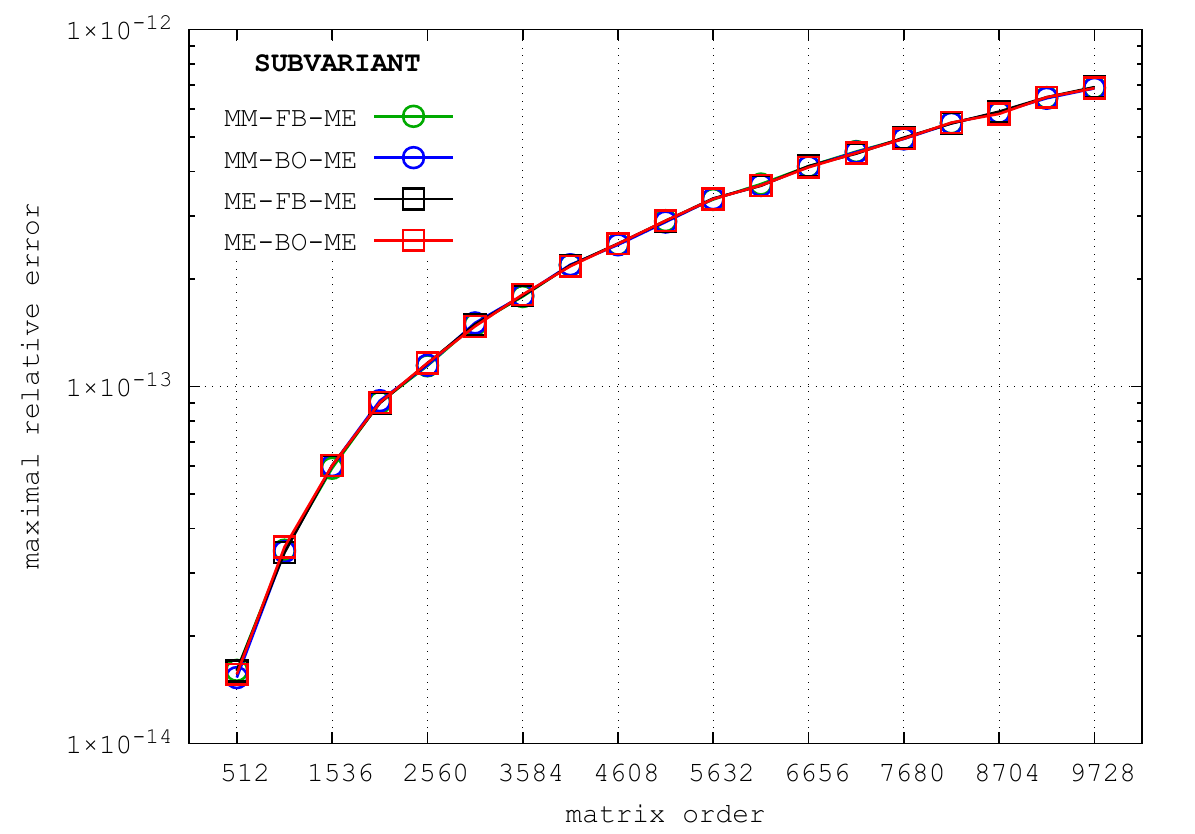}
  \caption{The relative normwise errors,
    $\|F-U\Sigma_F^{}X\|_F/\|F\|_F$, of four subvariants of DHZ0
    (left) and ZHZ0 (right) on the small real (left) and complex
    (right) datasets.}
  \label{fig:HZ0-nre1}
\end{figure*}

%
%
In the left subfigure of Figure~\ref{fig:HZ0-nre1} the normwise
relative errors of four subvariants of DHZ0 are shown in a logarithmic
scale, while in Table~\ref{tbl:DHZ0-nre1} the spectral condition
numbers $\mathop{\kappa_2^{}}(F)$ in the small real dataset are given,
offering a justification for the non-monotonicity of the error
graphs.
\begin{table}[h!bt]
  \caption{The spectral condition numbers $\mathop{\kappa_2^{}}(F)$ in
    the small real dataset.}
  \label{tbl:DHZ0-nre1}
  \centering
  \begin{tabular}{@{}cccc@{}}
    \toprule
    $n$ & $\mathop{\kappa_2^{}}(F)$ & $n$ & $\mathop{\kappa_2^{}}(F)$\\
    \midrule
    \hphantom{0}512 & $7.19081\cdot 10^{4}$ & 5632 & $2.33003\cdot 10^{6}$\\
    1024 & $1.29837\cdot 10^{6}$ & 6144 & $3.69195\cdot 10^{6}$\\
    1536 & $2.35302\cdot 10^{5}$ & 6656 & $1.90070\cdot 10^{6}$\\
    2048 & $1.48022\cdot 10^{5}$ & 7168 & $7.72806\cdot 10^{5}$\\
    2560 & $1.69855\cdot 10^{5}$ & 7680 & $5.33036\cdot 10^{5}$\\
    3072 & $7.28415\cdot 10^{4}$ & 8192 & $3.59719\cdot 10^{5}$\\
    3584 & $2.78307\cdot 10^{5}$ & 8704 & $3.54038\cdot 10^{5}$\\
    4096 & $1.43141\cdot 10^{6}$ & 9216 & $1.68270\cdot 10^{6}$\\
    4608 & $1.35132\cdot 10^{6}$ & 9728 & $3.90607\cdot 10^{5}$\\
    5120 & $2.33209\cdot 10^{6}$ & -- & --\\
    \bottomrule
  \end{tabular}
\end{table}
In Table~\ref{tbl:DHZ0-nre2} the spectral condition numbers
$\mathop{\kappa_2^{}}(G)$ in the small real dataset are given, without
a corresponding error figure, which is almost indistinguishable from the
left subfigure of Figure~\ref{fig:HZ0-nre1}.
\begin{table}[h!bt]
  \caption{The spectral condition numbers $\mathop{\kappa_2^{}}(G)$ in
    the small real dataset.}
  \label{tbl:DHZ0-nre2}
  \centering
  \begin{tabular}{@{}cccc@{}}
    \toprule
    $n$ & $\mathop{\kappa_2^{}}(G)$ & $n$ & $\mathop{\kappa_2^{}}(G)$\\
    \midrule
    \hphantom{0}512 & $7.79269\cdot 10^{3}$ & 5632 & $1.94774\cdot 10^{6}$\\
    1024 & $5.17836\cdot 10^{4}$ & 6144 & $1.14571\cdot 10^{7}$\\
    1536 & $2.41272\cdot 10^{4}$ & 6656 & $3.21194\cdot 10^{7}$\\
    2048 & $1.91692\cdot 10^{4}$ & 7168 & $6.84297\cdot 10^{6}$\\
    2560 & $3.30063\cdot 10^{4}$ & 7680 & $2.24983\cdot 10^{7}$\\
    3072 & $1.90896\cdot 10^{4}$ & 8192 & $2.89085\cdot 10^{6}$\\
    3584 & $6.56610\cdot 10^{4}$ & 8704 & $1.20004\cdot 10^{6}$\\
    4096 & $2.08194\cdot 10^{5}$ & 9216 & $3.35081\cdot 10^{6}$\\
    4608 & $8.83907\cdot 10^{5}$ & 9728 & $3.58109\cdot 10^{6}$\\
    5120 & $1.53531\cdot 10^{6}$ & -- & --\\
    \bottomrule
  \end{tabular}
\end{table}

\begin{figure*}[h!btp]
  \centering
  \includegraphics[keepaspectratio,width=.4975\textwidth]{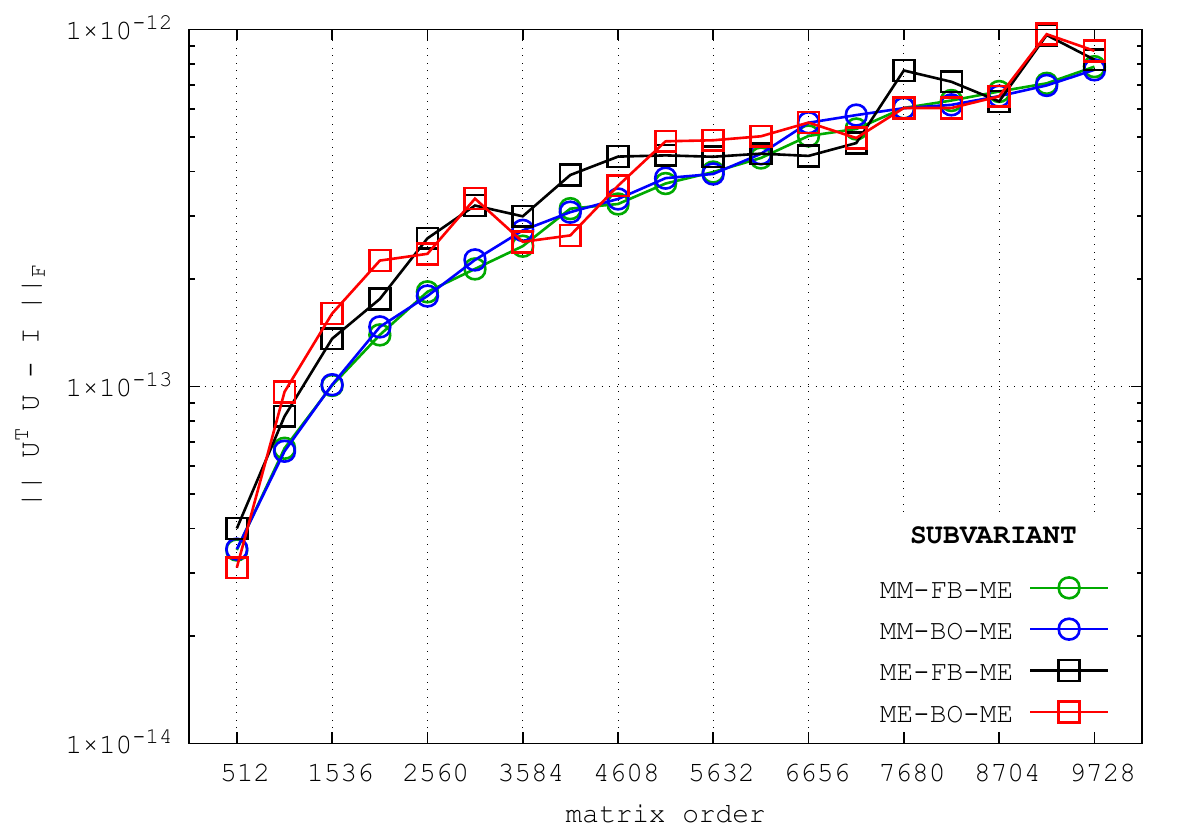}
  \includegraphics[keepaspectratio,width=.4975\textwidth]{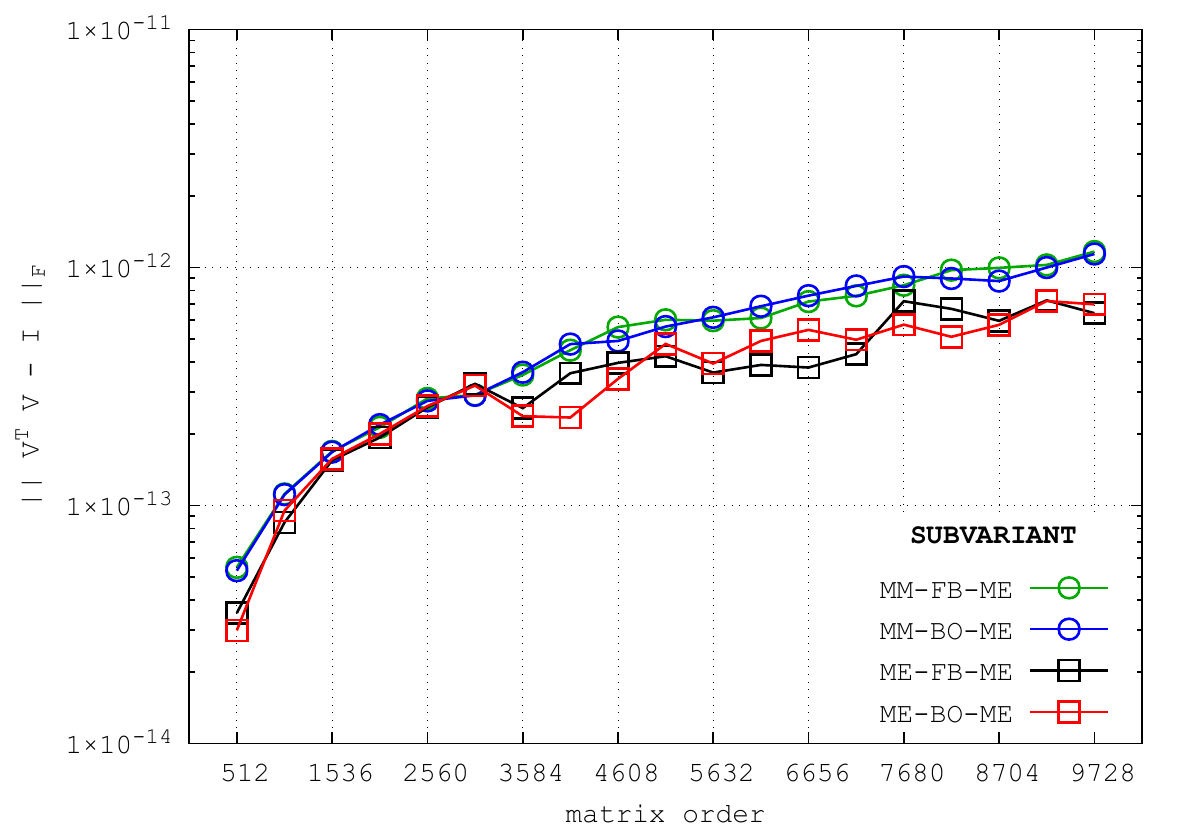}
  \caption{The numerical orthogonality $\|U^TU-I\|_F^{}$ (left) and
    $\|V^TV-I\|_F^{}$ (right) of the left generalized singular vectors
    $U$ and $V$, respectively, achieved by four subvariants of DHZ0 on
    the small real dataset.}
  \label{fig:DHZ0-ortho}
\end{figure*}

Figure~\ref{fig:DHZ0-ortho} shows the numerical orthogonality of the
left generalized singular vectors $U$ (left subfigure) and $V$ (right
subfigure) across the small real dataset achieved by DHZ0.

Tables~\ref{tbl:DHZ0-re} and \ref{tbl:DTGSJA-re} show the maximal
relative errors in the generalized singular values computed by
DHZ0-(\textsc{me}-\textsc{fb}-\textsc{me}) and \texttt{DTGSJA} from
the Intel Math Kernel Library (version 2020.1.217 on an Intel Xeon Phi
7210 CPU), respectively, on the small real dataset.
\begin{table}[h!bt]
  \caption{\looseness=-1 The relative errors in $\widehat{\Sigma}$,
    the generalized singular val\-ues computed by
    DHZ0-(\textsc{me}-\textsc{fb}-\textsc{me}) on the small real
    dataset.}
  \label{tbl:DHZ0-re}
  \centering
  \addtolength{\tabcolsep}{-1pt}
  \begin{tabular}{@{}cccc@{}}
    \toprule
    $n$ & $\max_i|(\sigma_i-\hat{\sigma}_i)/\sigma_i|$ & $n$ & $\max_i|(\sigma_i-\hat{\sigma}_i)/\sigma_i|$\\
    \midrule
    \hphantom{0}512 & $4.20077\cdot 10^{-13}$ & 5632 & $6.66337\cdot 10^{-11}$\\
               1024 & $2.00541\cdot 10^{-11}$ & 6144 & $6.85154\cdot 10^{-11}$\\
               1536 & $7.29365\cdot 10^{-12}$ & 6656 & $3.39017\cdot 10^{-10}$\\
               2048 & $8.80861\cdot 10^{-13}$ & 7168 & $5.65600\cdot 10^{-11}$\\
               2560 & $1.08576\cdot 10^{-12}$ & 7680 & $8.18645\cdot 10^{-11}$\\
               3072 & $1.12990\cdot 10^{-12}$ & 8192 & $1.51962\cdot 10^{-11}$\\
               3584 & $6.02566\cdot 10^{-12}$ & 8704 & $8.99801\cdot 10^{-12}$\\
               4096 & $1.64618\cdot 10^{-11}$ & 9216 & $1.81947\cdot 10^{-11}$\\
               4608 & $2.01553\cdot 10^{-11}$ & 9728 & $4.48431\cdot 10^{-11}$\\
               5120 & $4.85422\cdot 10^{-11}$ & -- & --\\
    \bottomrule
  \end{tabular}
\end{table}
\begin{table}[h!bt]
  \caption{\looseness=-1 The relative errors in $\widetilde{\Sigma}$,
    the generalized singular val\-ues computed by LAPACK's
    \texttt{DTGSJA} on the small real dataset.}
  \label{tbl:DTGSJA-re}
  \centering
  \addtolength{\tabcolsep}{-1pt}
  \begin{tabular}{@{}cccc@{}}
    \toprule
    $n$ & $\max_i|(\sigma_i-\tilde{\sigma}_i)/\sigma_i|$ & $n$ & $\max_i|(\sigma_i-\tilde{\sigma}_i)/\sigma_i|$\\
    \midrule
    \hphantom{0}512 & $4.24173\cdot 10^{-13}$ & 5632 & $6.65480\cdot 10^{-11}$\\
               1024 & $2.00948\cdot 10^{-11}$ & 6144 & $6.85305\cdot 10^{-11}$\\
               1536 & $7.28953\cdot 10^{-12}$ & 6656 & $3.38985\cdot 10^{-10}$\\
               2048 & $9.06332\cdot 10^{-13}$ & 7168 & $5.65602\cdot 10^{-11}$\\
               2560 & $1.04699\cdot 10^{-12}$ & 7680 & $8.18993\cdot 10^{-11}$\\
               3072 & $1.10638\cdot 10^{-12}$ & 8192 & $1.51537\cdot 10^{-11}$\\
               3584 & $6.00819\cdot 10^{-12}$ & 8704 & $9.00031\cdot 10^{-12}$\\
               4096 & $1.64629\cdot 10^{-11}$ & 9216 & $1.81785\cdot 10^{-11}$\\
               4608 & $2.03583\cdot 10^{-11}$ & 9728 & $4.48437\cdot 10^{-11}$\\
               5120 & $4.84123\cdot 10^{-11}$ & -- & --\\
    \bottomrule
  \end{tabular}
\end{table}

The average relative errors in the computed generalized singular
values in both the GPU and the CPU case are two to three orders of
magnitude smaller than the maximal ones, as can be seen for the former
in the supplementary material.  By comparing Tables~\ref{tbl:DHZ0-re}
and \ref{tbl:DTGSJA-re} it can be concluded that, in this sense and
instance, the proposed GPU algorithm exhibits accuracy similar to the
LAPACK's \texttt{DTGSJA}.
%
%
\subsubsection{Accuracy in the complex case}\label{sss:5.2.6}
%
%
In the right subfigure of Figure~\ref{fig:HZ0-nre1} the normwise
relative errors of four subvariants of ZHZ0 are shown in a logarithmic
scale, while in Table~\ref{tbl:ZHZ0-nre1} the spectral condition
numbers $\mathop{\kappa_2^{}}(F)$ in the small complex dataset are
presented, offering a justification for a smooth shape of the error
graphs.
\begin{table}[h!bt]
  \caption{The spectral condition numbers $\mathop{\kappa_2^{}}(F)$ in
    the small complex dataset.}
  \label{tbl:ZHZ0-nre1}
  \centering
  \begin{tabular}{@{}cccc@{}}
    \toprule
    $n$ & $\mathop{\kappa_2^{}}(F)$ & $n$ & $\mathop{\kappa_2^{}}(F)$\\
    \midrule
    \hphantom{0}512 & $2.30011\cdot 10^{3}$ & 5632 & $8.81464\cdot 10^{3}$\\
    1024 & $2.30011\cdot 10^{3}$ & 6144 & $8.81464\cdot 10^{3}$\\
    1536 & $2.30083\cdot 10^{3}$ & 6656 & $8.81464\cdot 10^{3}$\\
    2048 & $2.30083\cdot 10^{3}$ & 7168 & $8.81464\cdot 10^{3}$\\
    2560 & $4.92822\cdot 10^{3}$ & 7680 & $8.81464\cdot 10^{3}$\\
    3072 & $4.92822\cdot 10^{3}$ & 8192 & $8.81464\cdot 10^{3}$\\
    3584 & $4.92822\cdot 10^{3}$ & 8704 & $8.81464\cdot 10^{3}$\\
    4096 & $4.92822\cdot 10^{3}$ & 9216 & $8.81464\cdot 10^{3}$\\
    4608 & $7.69922\cdot 10^{3}$ & 9728 & $8.81464\cdot 10^{3}$\\
    5120 & $8.81464\cdot 10^{3}$ & -- & --\\
    \bottomrule
  \end{tabular}
\end{table}
In Table~\ref{tbl:ZHZ0-nre2} the spectral condition numbers
$\mathop{\kappa_2^{}}(G)$ in the small complex dataset are given.  A
corresponding error figure would be almost identical to the right
subfigure of Figure~\ref{fig:HZ0-nre1}.
\begin{table}[h!bt]
  \caption{The spectral condition numbers $\mathop{\kappa_2^{}}(G)$ in
    the small complex dataset.}
  \label{tbl:ZHZ0-nre2}
  \centering
  \begin{tabular}{@{}cccc@{}}
    \toprule
    $n$ & $\mathop{\kappa_2^{}}(G)$ & $n$ & $\mathop{\kappa_2^{}}(G)$\\
    \midrule
    \hphantom{0}512 & $5.62075\cdot 10^{2}$ & 5632 & $7.79104\cdot 10^{3}$\\
    1024 & $1.54312\cdot 10^{3}$ & 6144 & $4.69173\cdot 10^{3}$\\
    1536 & $8.58280\cdot 10^{2}$ & 6656 & $1.79255\cdot 10^{4}$\\
    2048 & $1.69892\cdot 10^{3}$ & 7168 & $7.57668\cdot 10^{3}$\\
    2560 & $2.66973\cdot 10^{5}$ & 7680 & $1.12171\cdot 10^{5}$\\
    3072 & $4.12692\cdot 10^{4}$ & 8192 & $1.38621\cdot 10^{4}$\\
    3584 & $9.24332\cdot 10^{3}$ & 8704 & $9.11971\cdot 10^{3}$\\
    4096 & $1.67599\cdot 10^{4}$ & 9216 & $8.54830\cdot 10^{3}$\\
    4608 & $6.50253\cdot 10^{3}$ & 9728 & $5.97562\cdot 10^{3}$\\
    5120 & $8.40210\cdot 10^{3}$ & -- & --\\
    \bottomrule
  \end{tabular}
\end{table}

\begin{figure*}[h!btp]
  \centering
  \includegraphics[keepaspectratio,width=.4975\textwidth]{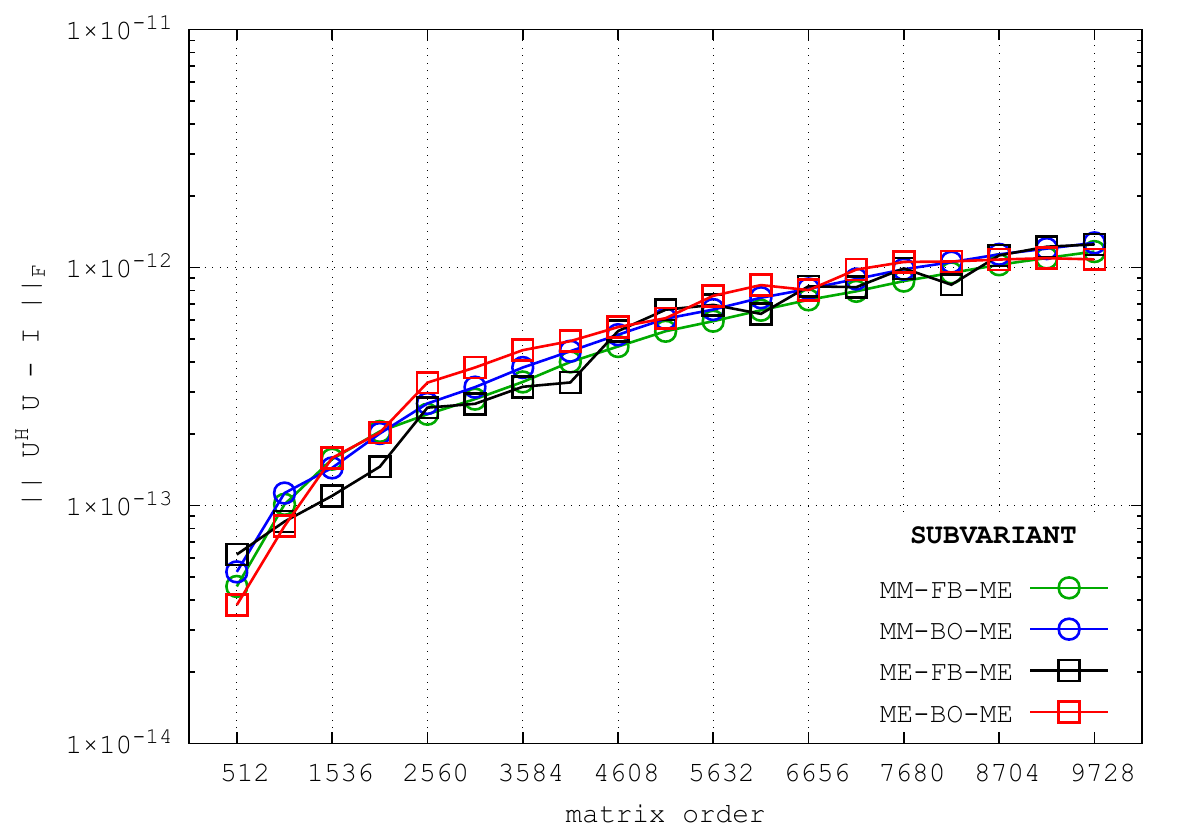}
  \includegraphics[keepaspectratio,width=.4975\textwidth]{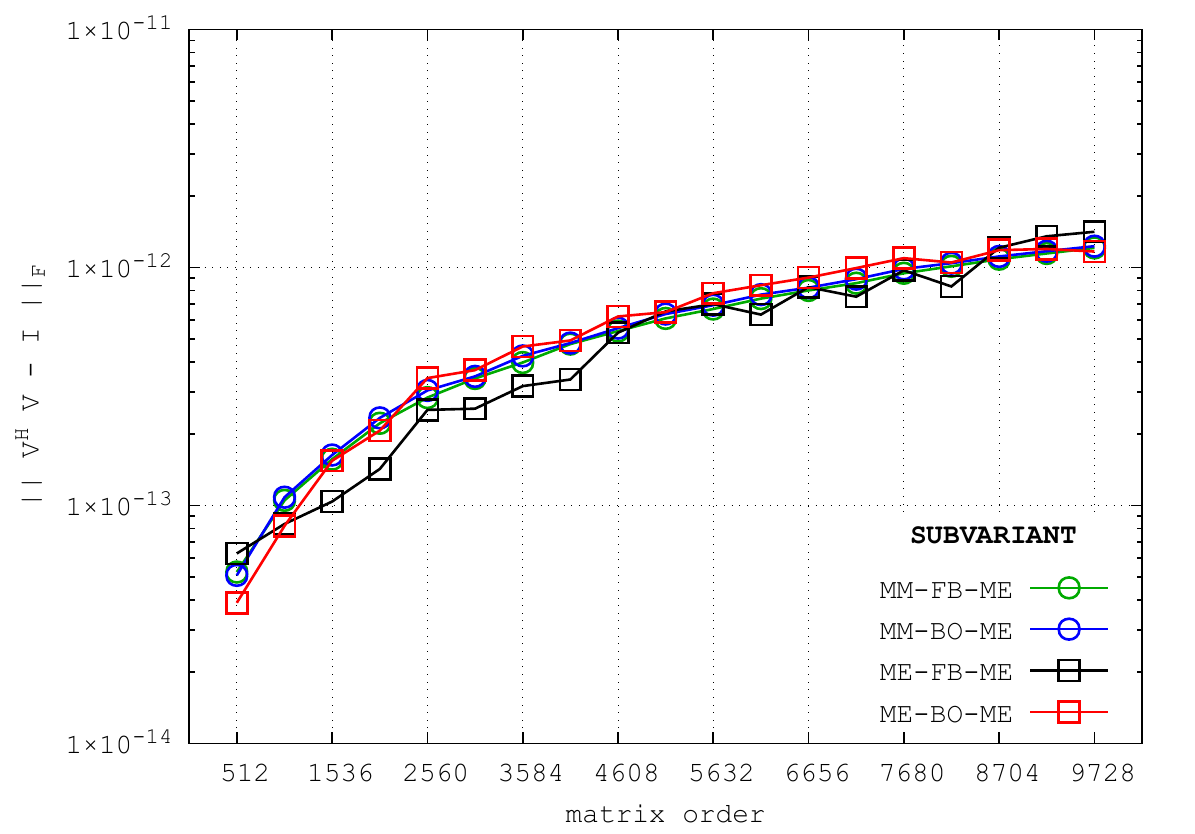}
  \caption{The numerical orthogonality $\|U^{\ast}U-I\|_F^{}$ (left) and
    $\|V^{\ast}V-I\|_F^{}$ (right) of the left generalized singular
    vectors $U$ and $V$, respectively, achieved by four subvariants of
    ZHZ0 on the small complex dataset.}
  \label{fig:ZHZ0-ortho}
\end{figure*}

Figure~\ref{fig:ZHZ0-ortho} shows the numerical orthogonality of the
left generalized singular vectors $U$ (left subfigure) and $V$ (right
subfigure) on the small complex dataset achieved by ZHZ0.
%
%
\subsubsection{Conclusions}\label{sss:5.2.7}
%
%
From Figure~\ref{fig:HZ-ts} it is obvious that \textsc{mm} gives a
significantly slower execution than \textsc{me} at the outer level,
and that \textsc{fb} is slightly faster than \textsc{bo}.  From
Tables~\ref{tbl:Dslow} and \ref{tbl:Zslow} it is clear that the
execution times across the variants do not widely differ, and that the
enhanced dot-products from Appendix~\ref{s:A} add from a few percent
to something more than $15\%$ to the wall time.

In both the real and the complex case the variant 0 with
(\textsc{me}-\textsc{fb}-\textsc{me}) is a reasonable choice
performance-wise.

Regarding the normwise relative errors on the matrices of moderate
spectral conditions, as it is in the real case, the \textsc{mm}
subvariants are almost indistinguishable, as are the \textsc{me}
subvariants, with the latter being slightly more accurate than the
former, as shown in the left subfigure of Figure~\ref{fig:HZ0-nre1}.
With the matrices of small spectral condition, as it is in the complex
case, all subvariants are almost indistinguishable, as shown in the
right subfigure of Figure~\ref{fig:HZ0-nre1}.  The spectral condition
numbers were computed by Matlab R2019a.

In Table~\ref{tbl:Mre} the maximal relative normwise errors with
respect to $F$ and $G$ in the real and the complex case on the small
dataset for all variants of the single-GPU algorithm with
(\textsc{me}-\textsc{fb}-\textsc{me}) are given.
\begin{table*}[h!bt]
  \caption{The maximal relative normwise errors with respect to $F$
    and $G$ in the real and the complex case on the small dataset for
    all variants of the single-GPU algorithm with
    (\textsc{me}-\textsc{fb}-\textsc{me}).  The \textbf{minimal}
    value in each column is shown in bold.}
  \label{tbl:Mre}
  \centering
  \begin{tabular}{@{}ccccc@{}}
    \toprule
    & \multicolumn{2}{c}{real case} & \multicolumn{2}{c}{complex case}\\
    \multirow{2}{*}{ID} & max.~relative & max.~relative & max.~relative & max.~relative\\
    & error w.r.t.~$F$ & error w.r.t.~$G$ & error w.r.t.~$F$ & error w.r.t.~$G$\\
    \midrule
    0 & $3.68432\cdot 10^{-12}$ & $3.70732\cdot 10^{-12}$ & $6.89432\cdot 10^{-13}$ & $6.89366\cdot 10^{-13}$\\
    1 & $3.68346\cdot 10^{-12}$ & $3.70057\cdot 10^{-12}$ & $6.89297\cdot 10^{-13}$ & $6.89204\cdot 10^{-13}$\\
    2 & $3.68803\cdot 10^{-12}$ & $3.69833\cdot 10^{-12}$ & $6.89300\cdot 10^{-13}$ & $6.89375\cdot 10^{-13}$\\
    3 & $3.68659\cdot 10^{-12}$ & $3.70446\cdot 10^{-12}$ & $\mathbf{6.86220\cdot 10^{-13}}$ & $6.90340\cdot 10^{-13}$\\
    4 & $3.68606\cdot 10^{-12}$ & $3.72483\cdot 10^{-12}$ & $6.89432\cdot 10^{-13}$ & $6.89366\cdot 10^{-13}$\\
    5 & $\mathbf{3.67729\cdot 10^{-12}}$ & $3.71555\cdot 10^{-12}$ & $6.89297\cdot 10^{-13}$ & $\mathbf{6.89204\cdot 10^{-13}}$\\
    6 & $3.68803\cdot 10^{-12}$ & $\mathbf{3.69833\cdot 10^{-12}}$ & $6.89300\cdot 10^{-13}$ & $6.89375\cdot 10^{-13}$\\
    7 & $3.68659\cdot 10^{-12}$ & $3.70446\cdot 10^{-12}$ & $6.86220\cdot 10^{-13}$ & $6.90340\cdot 10^{-13}$\\
    \bottomrule
  \end{tabular}
\end{table*}
Looking at the minimal value in each data column, it is evident that,
except for $G$ in the real case, the enhanced dot-products offer a
small advantage in accuracy, but not so significant that it would not
be offset by a drop in performance when the focus is on the latter.

Despite the low occupancy, Table~\ref{tbl:profiler} shows that all
warps that can occupy a multiprocessor execute with almost full
efficiency.
\begin{table}[h!bt]
  \caption{Aggregate minimal and maximal values of several profiler
    metrics for \texttt{bstep1s} kernel invocations.  Min/max was
    again taken over the results in four contexts: DHZ0 and ZHZ0, both
    in (\textsc{me}-\textsc{fb}-\textsc{me}) and
    (\textsc{me}-\textsc{bo}-\textsc{me}) subvariants; $n=9728$.}
  \label{tbl:profiler}
  \centering
  \addtolength{\tabcolsep}{-2pt}
  \begin{tabular}{@{}ccc@{}}
    \toprule
    \texttt{nvprof} metric & minimal & maximal\\
    \texttt{sm\_52} architecture & value [\%] & value [\%]\\
    \midrule
    achieved\_occupancy & $25.00$ & $25.00$\\
    branch\_efficiency & $99.75$ & $99.94$\\
    sm\_efficiency & $96.54$ & $97.64$\\
    warp\_execution\_efficiency & $99.92$ & $99.96$\\
    warp\_nonpred\_execution\_efficiency & $98.57$ & $99.90$\\
    \bottomrule
  \end{tabular}
\end{table}
It is therefore expected that, should the present bottleneck on
multiprocessors be removed in the future by increasing the register
file (see subsection~\ref{ss:3.5}) as it has recently been done with
the shared memory, the occupancy might increase as well, and with it
the overall performance, while the efficiency would stay at the same
high level.
%
%
\subsection{Results with the multi-GPU algorithm}\label{ss:5.3}
%
%
\looseness=-1
Due to a limited availability of GPUs in the testing environment, only
the complex case in the variant 0 was tested on the large dataset and
compared with a single-GPU baseline.

Here, the full accuracy testing was skipped, due to the huge
computational demands of the matrix inversions and of the error
calculation in extended precision.  For $n=18432$, the relative errors
in the corresponding outputs with one and with two GPUs (in both cases
in all variants that were timed) were compared and all of them, for
both $F$ and $G$, were found to be less than $1.7\cdot 10^{-12}$.
Also, for a given matrix, the relative errors in all cases differed
less than $10^{-13}$, indicating that the multi-GPU algorithm did not
introduce any instability in the computation.
%
%
\subsubsection{A single-GPU baseline}\label{sss:5.3.1}
%
%
In Table~\ref{tbl:ZHZ0s} the wall time in seconds and the number of
the outermost sweeps are shown for the
\textsc{me}-\textsc{fb}-\textsc{me} and the
\textsc{me}-\textsc{bo}-\textsc{me} subvariants, with and without the
column sorting, of the ZHZ0 single-GPU algorithm, as a baseline for
the comparison with the multi-GPU algorithm.  As the similar benefits
of the column sorting were obvious in other trial tests runs, the
non-sorting version was not considered for the full testing.
\begin{table}[h!bt]
  \caption{Wall time in seconds and the sweep count for the
    \textsc{me}-\textsc{fb}-\textsc{me} and the
    \textsc{me}-\textsc{bo}-\textsc{me} subvariants of the ZHZ0
    single-GPU algorithm, with and without the column sorting, on a
    pair of matrices of order $n=18432$.}
  \label{tbl:ZHZ0s}
  \centering
  \begin{tabular}{@{}ccc@{}}
    \toprule
    \textsc{me}-\textsc{fb}-\textsc{me} & \textsc{me}-\textsc{bo}-\textsc{me} & column sort\\
    \midrule
    2063.19~s; 10 & 2248.48~s; 11 & yes\\
    3509.21~s; 17 & 3482.66~s; 17 & no\\
    \bottomrule
  \end{tabular}
\end{table}
%
%
\subsubsection{The multi-GPU performance}\label{sss:5.3.2}
%
%
In Tables~\ref{tbl:ZHZ0-0},~\ref{tbl:ZHZ0-4},~\ref{tbl:ZHZ0-8},
and~\ref{tbl:ZHZ0-12} the wall time in seconds and the outermost sweep
count are shown for the multi-GPU variants%
\begin{displaymath}
  \begin{aligned}
    \mathcal{A}&\assgn \text{ZHZ0-(\textsc{me}-\textsc{fb},\,\textsc{me}-\textsc{fb}-\textsc{me})},\\
    \mathcal{B}&\assgn \text{ZHZ0-(\textsc{me}-\textsc{bo},\,\textsc{me}-\textsc{fb}-\textsc{me})},\\
    \mathcal{C}&\assgn \text{ZHZ0-(\textsc{me}-\textsc{fb},\,\textsc{me}-\textsc{bo}-\textsc{me})},\\
    \mathcal{D}&\assgn \text{ZHZ0-(\textsc{me}-\textsc{bo},\,\textsc{me}-\textsc{bo}-\textsc{me})},
  \end{aligned}
\end{displaymath}
respectively, run on two, four, and eight GPUs.  The tests on two and
four GPUs required a single node, and those on eight GPUs required two
InfiniBand-connected nodes.  When a particular test was not possible
to be run due to an insufficient amount of the GPU RAM, ``n/a'' is
shown in the test's table cell.
\begin{table}[h!bt]
  \caption{Wall time in seconds and the outermost sweep count for the variant $\mathcal{A}$.}
  \label{tbl:ZHZ0-0}
  \centering
  \addtolength{\tabcolsep}{-1pt}
  \begin{tabular}{@{}cccc@{}}
    \toprule
    $n$ & 2~GPUs & 4~GPUs & 8~GPUs\\
    \midrule
    18432 & \hphantom{0}3796.78~s; \hphantom{0}5 & \hphantom{0}2594.38~s; \hphantom{0}6 & \hphantom{0}1638.43~s; \hphantom{0}7\\
    24576 & \hphantom{0}8869.04~s; \hphantom{0}5 & \hphantom{0}6718.15~s; \hphantom{0}6 & \hphantom{0}5134.84~s; \hphantom{0}7\\
    36864 & n/a & 21271.83~s; \hphantom{0}6 & 12560.23~s; \hphantom{0}7\\
    \bottomrule
  \end{tabular}
\end{table}
\begin{table}[h!bt]
  \caption{Wall time in seconds and the outermost sweep count for the variant $\mathcal{B}$.}
  \label{tbl:ZHZ0-4}
  \centering
  \addtolength{\tabcolsep}{-1pt}
  \begin{tabular}{@{}cccc@{}}
    \toprule
    $n$ & 2~GPUs & 4~GPUs & 8~GPUs \\
    \midrule
    18432 & \hphantom{0}1606.17~s; \hphantom{0}9 & \hphantom{0}1085.34~s; \hphantom{0}9 & \hphantom{0}\hphantom{0}774.97~s; \hphantom{0}9\\
    24576 & \hphantom{0}3536.52~s; \hphantom{0}9 & \hphantom{0}2568.67~s; \hphantom{0}9 & \hphantom{0}2004.23~s; \hphantom{0}9\\
    36864 & n/a & \hphantom{0}7643.07~s; \hphantom{0}9 & \hphantom{0}4870.56~s; \hphantom{0}9\\
    \bottomrule
  \end{tabular}
\end{table}
\begin{table}[h!bt]
  \caption{Wall time in seconds and the outermost sweep count for the variant $\mathcal{C}$.}
  \label{tbl:ZHZ0-8}
  \centering
  \addtolength{\tabcolsep}{-1pt}
  \begin{tabular}{@{}cccc@{}}
    \toprule
    $n$ & 2~GPUs & 4~GPUs & 8~GPUs\\
    \midrule
    18432 & \hphantom{0}4295.87~s; \hphantom{0}5 & \hphantom{0}2903.28~s; \hphantom{0}6 & \hphantom{0}1806.63~s; \hphantom{0}7\\
    24576 & 10098.11~s; \hphantom{0}5 & \hphantom{0}7693.07~s; \hphantom{0}6 & \hphantom{0}5825.82~s; \hphantom{0}7\\
    36864 & n/a & 23611.93~s; \hphantom{0}6 & 14225.79~s; \hphantom{0}7\\
    \bottomrule
  \end{tabular}
\end{table}
\begin{table}[h!bt]
  \caption{Wall time in seconds and the outermost sweep count for the variant $\mathcal{D}$.}
  \label{tbl:ZHZ0-12}
  \centering
  \addtolength{\tabcolsep}{-1pt}
  \begin{tabular}{@{}cccc@{}}
    \toprule
    $n$ & 2~GPUs & 4~GPUs & 8~GPUs\\
    \midrule
    18432 & \hphantom{0}1759.48~s; 10 & \hphantom{0}1201.17~s; 10 & \hphantom{0}\hphantom{0}849.58~s; 10\\
    24576 & \hphantom{0}3926.69~s; 10 & \hphantom{0}2863.27~s; 10 & \hphantom{0}2224.73~s; 10\\
    36864 & n/a & \hphantom{0}9248.23~s; 11 & \hphantom{0}5877.18~s; 11\\
    \bottomrule
  \end{tabular}
\end{table}
%
%
\subsubsection{Conclusions}\label{sss:5.3.3}
%
%
It is clear from Tables~\ref{tbl:ZHZ0-0}, \ref{tbl:ZHZ0-4},
\ref{tbl:ZHZ0-8}, and~\ref{tbl:ZHZ0-12}, that the multi-GPU variant
$\mathcal{B}$ is to be recommended when performance matters.

Dividing the shortest single-GPU baseline wall time from
Table~\ref{tbl:ZHZ0s} for $n=18432$ with the wall times from
Table~\ref{tbl:ZHZ0-4} for the same $n$ but with a different number of
GPUs, it can be derived that the speedup with two GPUs is
$1.28\times$, with four GPUs is $1.90\times$, and with eight GPUs is
$2.66\times$.  These speedups could be even lower on a slower network
or if there were fewer than four GPUs present per node.

It would be interesting to see for what number of GPUs, depending on
the input sizes, the speedup peaks and starts falling, but that is
beyond reach of the testing environment.  In the absence of such
information, a safe rule of thumb would be to use a modest number of
GPUs on a fast interconnect for a given problem, such that they are
fully utilized in the terms of multiprocessors and memory.
%
%
\section{Conclusions and future work}\label{s:6}
%
%
The proposed algorithms compute the generalized SVD efficiently,
accurately, and almost entirely on the GPU(s).  The single-GPU
algorithm requires a CPU only for the controlling purposes.  The
multi-GPU algorithm involves a substantial amount of unavoidable
communication, but scales acceptably as long as each GPU if kept fully
utilized.

Several generalizations of the algorithms' design are possible for the
other implicit Jacobi-type methods that are to be ported to the
GPUs.  One such method is a computation of the generalized hyperbolic
SVD (GHSVD) \cite{Bojanczyk-2003} by a modification of the implicit
Hari--Zimmermann algorithm, as described
in~\cite{Singer-DiNapoli-Novakovic-Caklovic-2020}.
%
%
\begin{acks}
  This research was performed using the resources of computer cluster
  \href{https://www.srce.unizg.hr/isabella/}{Isabella} based in SRCE -
  University of Zagreb University Computing Centre.

  We would like to thank the associate editor and the anonymous
  referees for their recommendations for improving this manuscript.
\end{acks}
\begin{dci}
  The Authors declare that there is no conflict of interest.
\end{dci}
\begin{funding}
  This work has been supported in part by Croatian Science Foundation
  under the project IP--2014--09--3670.
\end{funding}
\begin{sm}
  The full testing results and some details of the chosen Jacobi
  strategies are provided as the supplementary material \texttt{sm.pdf}.

  The source code is available in
  \url{https://github.com/venovako/GPUHZGSVD} repository.  For
  comparing with the accuracy results a tag
  \url{https://github.com/venovako/GPUHZGSVD/tree/rev0} should be
  used.  Note that the actual names of the kernels, variables,
  etc.\ in the code are different.  The names in the paper are chosen
  for the simplicity of referencing.
\end{sm}
%

\begin{thebibliography}{38}
\providecommand{\natexlab}[1]{#1}
\providecommand{\url}[1]{\texttt{#1}}
\providecommand{\urlprefix}{URL }
\expandafter\ifx\csname urlstyle\endcsname\relax
  \providecommand{\doi}[1]{DOI:\discretionary{}{}{}#1}\else
  \providecommand{\doi}{DOI:\discretionary{}{}{}\begingroup
  \urlstyle{rm}\Url}\fi

\bibitem[{Alter et~al.(2003)Alter, Brown and
  Botstein}]{Alter-Brown-Botstein-2003}
Alter O, Brown PO and Botstein D (2003) Generalized singular value
  decomposition for comparative analysis of genome-scale expression data sets
  of two different organisms.
\newblock \emph{P. Natl. Acad. Sci. {USA}} 100(6): 3351--3356.
\newblock \doi{10.1073/pnas.0530258100}.

\bibitem[{Anderson et~al.(1999)Anderson, Bai, Bischof, Blackford, Demmel,
  Dongarra, Du~Croz, Greenbaum, Hammarling, McKenney and
  Sorensen}]{Anderson-et-al-99}
Anderson E, Bai Z, Bischof C, Blackford S, Demmel J, Dongarra J, Du~Croz J,
  Greenbaum A, Hammarling S, McKenney A and Sorensen D (1999) \emph{{LAPACK}
  Users' Guide}.
\newblock $3^{\rm rd}$ edition. Philadelphia, {PA}, {USA}: Society for
  Industrial and Applied Mathematics.
\newblock ISBN 0-89871-447-8 (paperback).

\bibitem[{Arafa et~al.(2019)Arafa, Badawy, Chennupati, Santhi and
  Eidenbenz}]{Arafa-et-al-2019}
Arafa Y, Badawy AA, Chennupati G, Santhi N and Eidenbenz S (2019) Low overhead
  instruction latency characterization for {NVIDIA} {GPGPU}s.
\newblock In: \emph{2019 IEEE High Performance Extreme Computing Conference
  (HPEC)}. pp. 1--8.
\newblock \doi{10.1109/HPEC.2019.8916466}.

\bibitem[{Bai(1994)}]{Bai-94}
Bai Z (1994) A parallel algorithm for computing the generalized singular value
  decomposition.
\newblock \emph{J. Parallel Distrib. Comput.} 20(3): 280--288.
\newblock \doi{10.1006/jpdc.1994.1027}.

\bibitem[{Bojanczyk(2003)}]{Bojanczyk-2003}
Bojanczyk AW (2003) An implicit {J}acobi-like method for computing generalized
  hyperbolic {SVD}.
\newblock \emph{Linear Algebra Appl.} 358(1): 293--307.
\newblock \doi{10.1016/S0024-3795(02)00394-4}.

\bibitem[{Boukaram et~al.(2018)Boukaram, Turkiyyah, Ltaief and
  Keyes}]{Boukaram-et-al-2018}
Boukaram WH, Turkiyyah G, Ltaief H and Keyes DE (2018) Batched {QR} and {SVD}
  algorithms on {GPU}s with applications in hierarchical matrix compression.
\newblock \emph{Parallel Comput.} 74: 19--33.
\newblock \doi{10.1016/j.parco.2017.09.001}.

\bibitem[{Cannon(1969)}]{Cannon-69}
Cannon LE (1969) \emph{A Cellular Computer to Implement the {K}alman Filter
  Algorithm}.
\newblock PhD Thesis, Montana State University, Bozeman, MT, USA.

\bibitem[{Cempel(2009)}]{Cempel-2009}
Cempel C (2009) Generalized singular value decomposition in multidimensional
  condition monitoring of machines--{A} proposal of comparative diagnostics.
\newblock \emph{Mech. Syst. Signal Pr.} 23(3): 701--711.
\newblock \doi{10.1016/j.ymssp.2008.07.004}.

\bibitem[{Drma{\v{c}}(1997)}]{Drmac-97}
Drma{\v{c}} Z (1997) Implementation of {J}acobi rotations for accurate singular
  value computation in floating point arithmetic.
\newblock \emph{{SIAM} J. Sci. Comput.} 18(4): 1200--1222.
\newblock \doi{10.1137/S1064827594265095}.

\bibitem[{Graillat et~al.(2015)Graillat, Lauter, Tang, Yamanaka and
  Oishi}]{Graillat-et-al-2015}
Graillat S, Lauter C, Tang PTP, Yamanaka N and Oishi S (2015) Efficient
  calculations of faithfully rounded {$l_2$}-norms of {$n$}-vectors.
\newblock \emph{{ACM} Trans. Math. Software} 41(4): art.~no.~24.
\newblock \doi{10.1145/2699469}.

\bibitem[{Hammond et~al.(2014)Hammond, Sch{\"{a}}fer and
  Latham}]{Hammond-Schaefer-Latham-2014}
Hammond JR, Sch{\"{a}}fer A and Latham R (2014) To \texttt{INT\_MAX}\ldots\ and
  beyond! {E}xploring large-count support in {MPI}.
\newblock In: \emph{2014 Workshop on Exascale MPI at Supercomputing
  Conference}. pp. 1--8.
\newblock \doi{10.1109/ExaMPI.2014.5}.

\bibitem[{Hari(1984)}]{Hari-84}
Hari V (1984) \emph{On Cyclic {J}acobi Methods for the Positive Definite
  Generalized Eigenvalue Problem}.
\newblock PhD Thesis, Fern{U}niversit{\"{a}}t--{G}esamthochschule, Hagen,
  Germany.

\bibitem[{Hari(2018)}]{Hari-2018}
Hari V (2018) Globally convergent {J}acobi methods for positive definite matrix
  pairs.
\newblock \emph{Numer. Algorithms} 79(1): 221--249.
\newblock \doi{10.1007/s11075-017-0435-5}.

\bibitem[{Hari(2019)}]{Hari-2019}
Hari V (2019) On the global convergence of the complex {HZ} method.
\newblock \emph{SIAM J. Matrix Anal. Appl.} 40(4): 1291--1310.
\newblock \doi{10.1137/19M1265594}.

\bibitem[{Hari et~al.(2010)Hari, Singer and
  Singer}]{Hari-SingerSanja-SingerSasa-2010}
Hari V, Singer S and Singer S (2010) Block-oriented {$J$}-{J}acobi methods for
  {H}ermitian matrices.
\newblock \emph{Linear Algebra Appl.} 433(8--10): 1491--1512.
\newblock \doi{10.1016/j.laa.2010.06.032}.

\bibitem[{Hari et~al.(2014)Hari, Singer and
  Singer}]{Hari-SingerSanja-SingerSasa-2014}
Hari V, Singer S and Singer S (2014) Full block {$J$}-{J}acobi method for
  {H}ermitian matrices.
\newblock \emph{Linear Algebra Appl.} 444: 1--27.
\newblock \doi{10.1016/j.laa.2013.11.028}.

\bibitem[{Howland et~al.(2003)Howland, Jeon and Park}]{Howland-Jeon-Park-2003}
Howland P, Jeon M and Park H (2003) Structure preserving dimension reduction
  for clustered text data based on the generalized singular value
  decomposition.
\newblock \emph{{SIAM} J. Matrix Anal. Appl.} 25(1): 165--179.
\newblock \doi{10.1137/S0895479801393666}.

\bibitem[{Howland et~al.(2006)Howland, Wang and Park}]{Howland-Wang-Park-2006}
Howland P, Wang J and Park H (2006) Solving the small sample size problem in
  face recognition using generalized discriminant analysis.
\newblock \emph{Pattern Recogn.} 39(2): 277--287.
\newblock \doi{10.1016/j.patcog.2005.06.013}.

\bibitem[{Luk(1985)}]{Luk-85}
Luk FT (1985) A parallel method for computing the generalized singular value
  decomposition.
\newblock \emph{J. Parallel Distrib. Comput.} 2(3): 250--260.
\newblock \doi{10.1016/0743-7315(85)90027-9}.

\bibitem[{Mantharam and Eberlein(1993)}]{Mantharam-Eberlein-93}
Mantharam M and Eberlein PJ (1993) Block recursive algorithm to generate
  {J}acobi--sets.
\newblock \emph{Parallel Comput.} 19: 481--496.
\newblock \doi{10.1016/0167-8191(93)90001-2}.

\bibitem[{Mary et~al.(2015)Mary, Yamazaki, Kurzak, Luszczek, Tomov and
  Dongarra}]{Mary-et-al-2015}
Mary T, Yamazaki I, Kurzak J, Luszczek P, Tomov S and Dongarra J (2015)
  Performance of random sampling for computing low-rank approximations of a
  dense matrix on {GPU}s.
\newblock In: \emph{Proceedings of the International Conference for High
  Performance Computing, Networking, Storage and Analysis (SC ’15)}. pp.
  1--11 (art.~no.~60).
\newblock \doi{10.1145/2807591.2807613}.

\bibitem[{{Message Passing Interface Forum}(2015)}]{mpi31}
{Message Passing Interface Forum} (2015) \emph{{MPI}: {A} {M}essage-{P}assing
  {I}nterface {S}tandard, {V}ersion 3.1}.
\newblock High Performance Computing Center Stuttgart (HLRS).

\bibitem[{Novakovi{\'{c}}(2015)}]{Novakovic-2015}
Novakovi{\'{c}} V (2015) A hierarchically blocked {J}acobi {SVD} algorithm for
  single and multiple graphics processing units.
\newblock \emph{{SIAM} J. Sci. Comput.} 37(1): C1--C30.
\newblock \doi{10.1137/140952429}.

\bibitem[{Novakovi{\'{c}}(2017)}]{Novakovic-2017}
Novakovi{\'{c}} V (2017) \emph{Parallel {J}acobi-type algorithms for the
  singular and the generalized singular value decomposition}.
\newblock PhD Thesis, University of Zagreb, Croatia.
\newblock \urlprefix\url{https://urn.nsk.hr/urn:nbn:hr:217:515320}.

\bibitem[{Novakovi{\'{c}} and Singer(2011)}]{Novakovic-SingerSanja-2011}
Novakovi{\'{c}} V and Singer S (2011) A {GPU}-based hyperbolic {SVD} algorithm.
\newblock \emph{{BIT}} 51(4): 1009--1030.
\newblock \doi{10.1007/s10543-011-0333-5}.

\bibitem[{Novakovi{\'{c}} et~al.(2015)Novakovi{\'{c}}, Singer and
  Singer}]{Novakovic-Singer-Singer-2015}
Novakovi{\'{c}} V, Singer S and Singer S (2015) Blocking and parallelization of
  the {H}ari--{Z}immermann variant of the {F}alk--{L}angemeyer algorithm for
  the generalized {SVD}.
\newblock \emph{Parallel Comput.} 49: 136--152.
\newblock \doi{10.1016/j.parco.2015.06.004}.

\bibitem[{{NVIDIA Corp.}(2019)}]{NVidia-2019}
{NVIDIA Corp} (2019) \emph{{CUDA} {C} {P}rogramming {G}uide v10.1.243}.
\newblock
  \urlprefix\url{https://docs.nvidia.com/cuda/archive/10.1/cuda-c-programming-guide/}.

\bibitem[{Paige and Saunders(1981)}]{Paige-Saunders-81}
Paige CC and Saunders MA (1981) Towards a generalized singular value
  decomposition.
\newblock \emph{{SIAM} J. Numer. Anal.} 18(3): 398--405.
\newblock \doi{10.1137/0718026}.

\bibitem[{Quintana-Ort{\'{\i}} et~al.(1998)Quintana-Ort{\'{\i}}, Sun and
  Bischof}]{Quintana-Orti-Sun-Bischof-98}
Quintana-Ort{\'{\i}} G, Sun X and Bischof C (1998) A {BLAS}-3 version of the
  {QR} factorization with column pivoting.
\newblock \emph{SIAM J. Sci. Comp.} 19(5): 1486--1494.
\newblock \doi{10.1137/S1064827595296732}.

\bibitem[{Senaratne and Tellambura(2013)}]{Senaratne-Tellambura-2013}
Senaratne D and Tellambura C (2013) {GSVD} beamforming for two-user {MIMO}
  downlink channel.
\newblock \emph{{IEEE} T. Veh. Technol.} 62(6): 2596--2606.
\newblock \doi{10.1109/TVT.2013.2241091}.

\bibitem[{Singer et~al.(2020)Singer, {Di Napoli}, Novakovi{\'{c}} and
  {\v{C}}aklovi{\'{c}}}]{Singer-DiNapoli-Novakovic-Caklovic-2020}
Singer S, {Di Napoli} E, Novakovi{\'{c}} V and {\v{C}}aklovi{\'{c}} G (2020)
  The {LAPW} method with eigendecomposition based on the {H}ari--{Z}immermann
  generalized hyperbolic {SVD}.
\newblock \emph{SIAM J. Sci. Comput.} 42(5): C265--C293.
\newblock \doi{10.1137/19M1277813}.

\bibitem[{Singer et~al.(2012)Singer, Singer, Novakovi{\'{c}},
  U{\v{s}}{\'{c}}umli{\'{c}} and
  Dunjko}]{SingerSanja-SingerSasa-Novakovic-Uscumlic-Dunjko-2012}
Singer S, Singer S, Novakovi{\'{c}} V, U{\v{s}}{\'{c}}umli{\'{c}} A and Dunjko
  V (2012) Novel modifications of parallel {J}acobi algorithms.
\newblock \emph{Numer. Algorithms} 59: 1--27.
\newblock \doi{10.1007/s11075-011-9473-6}.

\bibitem[{Slapni{\v{c}}ar(1998)}]{Slapnicar-98}
Slapni{\v{c}}ar I (1998) Componentwise analysis of direct factorization of real
  symmetric and {H}ermitian matrices.
\newblock \emph{Linear Algebra Appl.} 272: 227--275.
\newblock \doi{10.1016/S0024-3795(97)00334-0}.

\bibitem[{Tomov et~al.(2010)Tomov, Dongarra and
  Baboulin}]{Tomov-Dongarra-Baboulin-2010}
Tomov S, Dongarra J and Baboulin M (2010) Towards dense linear algebra for
  hybrid {GPU} accelerated manycore systems.
\newblock \emph{Parallel Comput.} 36(5--6): 232--240.
\newblock \doi{10.1016/j.parco.2009.12.005}.

\bibitem[{{Van Loan}(1976)}]{VanLoan-76}
{Van Loan} CF (1976) Generalizing the singular value decomposition.
\newblock \emph{{SIAM} J. Numer. Anal.} 13(1): 76--83.
\newblock \doi{10.1137/0713009}.

\bibitem[{Zhao et~al.(2010)Zhao, Rutkowski, Zhang and
  Cichocki}]{Zhao-et-al-2010}
Zhao Q, Rutkowski TM, Zhang L and Cichocki A (2010) Generalized optimal spatial
  filtering using a kernel approach with application to {EEG} classification.
\newblock \emph{Cogn. Neurodyn.} 4(4): 355--358.
\newblock \doi{10.1007/s11571-010-9125-x}.

\bibitem[{Zhou and van~der Veen(2017)}]{Zhou-vanderVeen-2017}
Zhou M and van~der Veen AJ (2017) Blind separation of partially overlapping
  data packets.
\newblock \emph{Digit. Signal Process.} 68: 154--166.
\newblock \doi{10.1016/j.dsp.2017.06.009}.

\bibitem[{Zimmermann(1969)}]{Zimmermann-69}
Zimmermann K (1969) \emph{{Z}ur {K}onvergenz eines {J}acobiverfahren f{\"{u}}r
  gew{\"{o}}nliche und verallgemeinerte {E}igenwertprobleme}.
\newblock Dissertation no.~4305, {E}idgen{\"{o}}ssische {T}echnische
  {H}ochschule, Z{\"{u}}rich, Switzerland.

\end{thebibliography}

%
\appendix
%
%
\section{Enhanced dot-product computation}\label{s:A}
%
%
In CUDA, the rounding mode can be specified explicitly for each
arithmetic operation using the intrinsic functions.  That makes an
ideal setting for employing a trick from~\cite{Graillat-et-al-2015} to
cheaply compute possibly more accurate real and complex dot-products.

For two real vectors $\mathbf{a}$ and $\mathbf{b}$ of the same length,
their enhanced dot-product would require one FMA and one negation per
a pair of vector elements, in addition to one multiplication.  Also,
adding the partial sums together needs two sum-reductions instead of
one, as follows.

Let $a$ and $b$ be the elements of $\mathbf{a}$ and $\mathbf{b}$ at
the same, arbitrary index, and let \texttt{\_\_dmul\_rd} stand for a CUDA
intrinsic performing a multiplication with rounding towards $-\infty$.
Take
\begin{displaymath}
  c\assgn\mathop{\mathtt{\_\_dmul\_rd}}(a,b),\quad
  d\assgn\mathop{\mathtt{\_\_fma\_rn}}(a,b,-c).
\end{displaymath}
Then, by looking separately at the both possible signs of $c$, it can
be shown that for the rounding error of the multiplication extracted
by the FMA holds $d \ge 0$.  By sum-reducing $d$, no cancellation can
occur. That value may be added to the sum-reduction result on $c$, to
form the final dot-product.

For two complex vectors $\mathbf{a}$ and $\mathbf{b}$ of the same
length, defining an enhanced dot-product is not so unambiguous.  A
special case of computing
$\|\mathbf{a}\|_2^2 = \mathbf{a}^{\ast} \mathbf{a}^{}$ can be handled
as follows.  For an element $a$ of $\mathbf{a}$, take
\begin{displaymath}
  \begin{aligned}
    c_{\mathrm{r}}^{}\assgn&\mathop{\mathtt{\_\_dmul\_rd}}(\Real(a),\Real(a)),\\
    d_{\mathrm{r}}^{}\assgn&\mathop{\mathtt{\_\_fma\_rn}}(\Real(a),\Real(a),-c_{\mathrm{r}}^{}),
  \end{aligned}
\end{displaymath}
considering the real part of $a$, and
\begin{displaymath}
  \begin{aligned}
    c_{\mathrm{i}}^{}\assgn&\mathop{\mathtt{\_\_dmul\_rd}}(\Imag(a),\Imag(a)),\\
    d_{\mathrm{i}}^{}\assgn&\mathop{\mathtt{\_\_fma\_rn}}(\Imag(a),\Imag(a),-c_{\mathrm{i}}^{}),
  \end{aligned}
\end{displaymath}
considering the imaginary part of $a$.  Let
$\mathop{\Sigma}c_{\mathrm{r}}^{}$,
$\mathop{\Sigma}c_{\mathrm{i}}^{}$,
$\mathop{\Sigma}d_{\mathrm{r}}^{}$, and
$\mathop{\Sigma}d_{\mathrm{i}}^{}$ be the sum-reductions of all
$c_{\mathrm{r}}^{}$, $c_{\mathrm{i}}^{}$, $d_{\mathrm{r}}^{}$, and
$d_{\mathrm{i}}^{}$, respectively.  Then, return
$((\sigma_1^{} + \sigma_2^{}) + \sigma_3^{}) + \sigma_4^{}$, where the
$\sigma_{\null\!j}^{}$ quantities are the four $\Sigma$-sums above,
ordered increasingly.

However, such an approach requires four reductions.  To simplify it by
trading off accuracy for speed, let
$e = d_{\mathrm{r}}^{} + d_{\mathrm{i}}^{}$ be the sum of the rounding
errors.  Then, if $c_{\mathrm{r}}^{} \le c_{\mathrm{i}}^{}$, take 
$s = (e + c_{\mathrm{r}}^{}) + c_{\mathrm{i}}^{}$, else let
$s = (e + c_{\mathrm{i}}^{}) + c_{\mathrm{r}}^{}$, and return the
sum-reduction of $s$ as an approximation of
$\mathop{\Sigma}(\bar{a}\cdot a)$.

The real dot-product and the simpler of the two procedures for
computing of the square of the Euclidean norm of a complex vector have
been incorporated in the special variants of the algorithm and tested,
both without a huge slowdown but also without a significant effect on
accuracy of the result, as explained in section~\ref{s:5}.
%
%
\section{A single-GPU out-of-core algorithm}\label{s:B}
%
%
Algorithm~\ref{alg:ooc} gives an overview of a single-GPU
``out-of-core'' algorithm, when the whole data does not fit into the
GPU RAM, but at least one block pair of $F$, $G$, and $Z$ does.
\begin{algorithm}[h!bt]
  \SetKwFunction{MemCpy}{cudaMemcpy}
  \SetKwFunction{InitFGZ}{initFGZ}
  \SetKwFunction{BstepI}{OoCstep1$($s$/$n$)$}
  \For(\tcp*[f]{for all outer steps}){$0\le k<\mathtt{n}'$}{
    \tcc{Let $\mathbf{o}$ be the number of block column pairs of $F$,
      $G$, and $Z$ that fit into the GPU RAM.}
    \For(\tcp*[f]{$\mathbf{o}\bmod(\mathtt{n}/2)=0$}){$0\le o<(\mathtt{n}/2)/\mathbf{o}$}{
      \For(\tcp*[f]{async}){$o\cdot\mathbf{o}\le\ell<(o+1)\cdot\mathbf{o}$}{
        copy the $\ell$th block pivot pairs of $F$, $G$, and $Z$ from the
        CPU to the GPU\;
      }
      \BstepI{$o,k$}\tcp*[r]{transform}
      \For(\tcp*[f]{async}){$o\cdot\mathbf{o}\le\ell<(o+1)\cdot\mathbf{o}$}{
        copy the $\ell$th block pivot pairs of $F$, $G$, $Z$, transformed,
        from the GPU to the CPU\;
      }
    }
  }
  \caption{A split of the implicit $\ell$-loop of a single outer step
    of Algorithm~\ref{alg:cpu1} for the case when the whole data does
    not fit into the GPU RAM\@.}
  \label{alg:ooc}
\end{algorithm}

\looseness=-1
This algorithm has not been implemented, but essentially the
\texttt{bstep1(s/n)} kernel would be converted to a kernel
\texttt{OoCstep1(s/n)} that operates on a suitably sized subset of all
block pivots of an outer step.  The subset is loaded into the GPU RAM
beforehand, and is brought back, transformed, to the CPU\@.  That
process is repeated until all block pivots in the outer step have been
transformed, exactly as they would be by a single call to
\texttt{bstep1(s/n)}.  Apart from being inefficient, such an approach
cannot handle the case when not even a single block pair can fit into
a GPU\@.  For the GPU RAM being $16\,\mathrm{GiB}$ and the matrices
complex and square, $n$ would then have to be larger than 11 million,
in which case the multi-GPU algorithm is a must anyway.
%
%
\section{Warp-shuffle $+$-reduction of 32 doubles}\label{s:C}
%
%
In Figure~\ref{fig:sum32} a warp-shuffle sum-reduction of 32
double-precision values is shown, where each thread in a warp holds
one value at the start, and all threads get the sum at the end.  The
loop over \texttt{i} is manually unrolled in the code.
\begin{figure}[h!bt]
{\small\begin{verbatim}
__device__ __forceinline__
double dSum32(const double x) {
  int lo0, hi0, lo1, hi1;
  double x0 = x, x1;
  for (int i = 16; i; i >>= 1) {
    lo0 = __double2loint(x0);
    hi0 = __double2hiint(x0);
    lo1 = __shfl_xor_sync(~0u, lo0, i);
    hi1 = __shfl_xor_sync(~0u, hi0, i);
    x1 = __hiloint2double(hi1, lo1);
    x0 = __dadd_rn(x0, x1);
  }
  lo0 = __double2loint(x0);
  hi0 = __double2hiint(x0);
  lo1 = __shfl_sync(~0u, lo0, 0);
  hi1 = __shfl_sync(~0u, hi0, 0);
  return __hiloint2double(hi1, lo1);
}\end{verbatim}}
\caption{Sum-reduction of 32 doubles across a warp.}
\label{fig:sum32}
\end{figure}
%
%
\section{An implementation of the main kernel's subphases \textbf{1} and \textbf{4} with the batched matrix multiplication from cuBLAS}\label{s:D}
%
%
The subphases~\textbf{1} and \textbf{4} of the main kernel (see
subsection~\ref{ss:3.5}) could have also been decoupled and
implemented separately, by the batched matrix multiplication routines
from the cuBLAS library.  For that, more than twice the GPU RAM is
required, to maintain two copies of the data.  The first copy, called
``actual'', holds the current state of the computation at the start of
each block step.  The second one, called ``shadow'', is used as a
temporary buffer, in which the new state is assembled at the end of
the step.  These copies then change the roles by a simple swap of the
pointers to them.  Schematically, this layout can be represented as
\begin{displaymath}
  \mathrm{actual}\rightarrow\text{\fbox{\fbox{$F^!$}\,\fbox{$G^!$}\,\fbox{$Z^!$}}}\,,\ \ \
  \mathrm{shadow}\rightarrow\text{\fbox{\fbox{$F^?$}\,\fbox{$G^?$}\,\fbox{$Z^?$}}}\,.
\end{displaymath}

Let $\mathbb{F}\in\{\mathbb{R},\mathbb{C}\}$ be the field in which the
computations take place.  For the complex algorithm, a data layout
with the real and the imaginary components kept separately is no
longer viable because cuBLAS does not support it.

Furthermore, let $\mathtt{b}=\mathtt{n}/2$ be the number of block
column pairs processed at each block step, i.e., the number of thread
blocks launched.  An additional workspace of $\mathtt{b}$ square
matrices $\mathsf{Z}_l^{}$ of order $32$ is required, where
$1\le l\le\mathtt{b}$, as
\begin{displaymath}
  \settowidth{\fbw}{\fbox{$\mathsf{Z}_{\mathtt{b}-1}^{}$}}
  \mathrm{work}\rightarrow\text{\fbox{\framebox[\fbw][c]{$\mathsf{Z}_0^{}$}$\,\cdots\,$\framebox[\fbw][c]{$\mathsf{Z}_{\mathtt{b}-1}^{}$}}}\,.
\end{displaymath}
In this GPU RAM workspace the  transformation matrices
$\mathsf{Z}_l^{}$ accumulated in the subphase~\textbf{3} will be
stored from the shared memory of each thread block.

For each $l$ there is an associated pair of block indices
$(\mathtt{p}_l^{},\mathtt{q}_l^{})$ drawn from the chosen block
strategy for the current block step, with the subscripts $l$ omitted
when they are implied by the context.  Then, the block columns of the
actual copy 
can be logically grouped and denoted as
\begin{displaymath}
  F_l^!=
  \begin{bmatrix}
    F_{\mathtt{p}}^! & F_{\mathtt{q}}^!
  \end{bmatrix},\quad
  G_l^!=
  \begin{bmatrix}
    G_{\mathtt{p}}^! & G_{\mathtt{q}}^!
  \end{bmatrix},\quad
  Z_l^!=
  \begin{bmatrix}
    Z_{\mathtt{p}}^! & Z_{\mathtt{q}}^!
  \end{bmatrix},
\end{displaymath}
and a similar notation is used for the shadow copy as well.

Replacing the subphase~\textbf{1} involves forming the lower
block-triangle of the Grammian matrices $\widehat{\mathsf{A}}_l^{}$
and $\widehat{\mathsf{B}}_l^{}$ as
\begin{displaymath}
  \widehat{\mathsf{A}}_l^{}\assgn
  \begin{bmatrix}
    F_{\mathtt{p}}^{!\ast}\!\cdot\!F_{\mathtt{p}}^! & ???\\[2pt]
    F_{\mathtt{q}}^{!\ast}\!\cdot\!F_{\mathtt{p}}^! & F_{\mathtt{q}}^{!\ast}\!\cdot\!F_{\mathtt{q}}^!
  \end{bmatrix},\quad
  \widehat{\mathsf{B}}_l^{}\assgn
  \begin{bmatrix}
    G_{\mathtt{p}}^{!\ast}\!\cdot\!G_{\mathtt{p}}^! & ???\\[2pt]
    G_{\mathtt{q}}^{!\ast}\!\cdot\!G_{\mathtt{p}}^! & G_{\mathtt{q}}^{!\ast}\!\cdot\!G_{\mathtt{q}}^!
  \end{bmatrix},
\end{displaymath}
by two batched \texttt{GEMM} calls: one for
$\widehat{\mathsf{A}}_l^{}$, and another for
$\widehat{\mathsf{B}}_l^{}$.

Each batched \texttt{GEMM} call forms all $(\mathtt{p},\mathtt{p})$,
$(\mathtt{q},\mathtt{p})$, and $(\mathtt{q},\mathtt{q})$ blocks, since
the inputs are all of the same dimensions, as well as the outputs.
Forming the remaining, upper-right blocks is not necessary, since the
subphase~\textbf{2} reads only the lower triangle as the Grammian
matrices are Hermitian.  The diagonal blocks should be formed by the
batched \texttt{SYRK} or \texttt{HERK} operations instead, but at
present they are not implemented in cuBLAS\@.  Each
$\widehat{\mathsf{A}}_l^{}$ is stored in the first $32$ rows of the
shadow copy $F^?$, starting from its $(32\cdot l)$th column, and
similarly for $\widehat{\mathsf{B}}_l^{}$ and $G^?$, from where they
are loaded by the modified subphase~\textbf{2} into the shared memory
and their Cholesky factorizations are then computed as before.  Note
that this approach involves storing $\widehat{\mathsf{A}}_l^{}$ and
$\widehat{\mathsf{B}}_l^{}$ in the GPU RAM, only to be read
afterwards into the shared memory, what is a redundant memory traffic
avoided by the original implementation.

The modified main kernel, after completing the subphase~\textbf{2},
computes the transformation matrices $\mathsf{Z}_l^{}$ in the
subphase~\textbf{3}, and stores them from the shared memory to the
workspace (another redundant memory traffic).  For the
postmultiplications, i.e., the block column updates, consider each
$\mathsf{Z}_l^{}$ to be logically partitioned as
\begin{displaymath}
  \mathsf{Z}_l^{}=
  \begin{bmatrix}
    \mathsf{Z}_{\mathtt{1}\mathtt{1}}^{} & \mathsf{Z}_{\mathtt{1}\mathtt{2}}^{}\\
    \mathsf{Z}_{\mathtt{2}\mathtt{1}}^{} & \mathsf{Z}_{\mathtt{2}\mathtt{2}}^{}
  \end{bmatrix},\quad
  \mathsf{Z}_{\mathtt{i}\mathtt{j}}^{}\in\mathbb{F}^{16\times 16},\quad
  1\le\mathtt{i},\mathtt{j}\le 2.
\end{displaymath}

Now, compute a partial update of $F^!$ by a single batched
\texttt{GEMM} call (and similarly for $G^!$ and $Z^!$), where for each
$l$,
\begin{displaymath}
  F_{\mathtt{p}}^?\assgn F_{\mathtt{p}}^!\cdot\mathsf{Z}_{\mathtt{1}\mathtt{1}}^{},\quad
  F_{\mathtt{q}}^?\assgn F_{\mathtt{p}}^!\cdot\mathsf{Z}_{\mathtt{1}\mathtt{2}}^{},
\end{displaymath}
and complete the update by another batched \texttt{GEMM} call as
\begin{displaymath}
  F_{\mathtt{p}}^?\assgn F_{\mathtt{q}}^!\cdot\mathsf{Z}_{\mathtt{2}\mathtt{1}}^{}+F_{\mathtt{p}}^?,\quad
  F_{\mathtt{q}}^?\assgn F_{\mathtt{q}}^!\cdot\mathsf{Z}_{\mathtt{2}\mathtt{2}}^{}+F_{\mathtt{q}}^?.
\end{displaymath}

A swap of the pointers to the copies,
$\mathrm{actual}\leftrightarrows\mathrm{shadow}$, completes the block
step.  See the \texttt{var} subdirectory of the code repository for a
double precision implementation.

A benefit of this alternative is that it significantly reduces the
register pressure in the main kernel, e.g., close to halving the
required number of registers on a Volta GPU\@.  However, that is
offset by allocating more than twice the GPU RAM, which is still a
scarce resource, and consequently by reducing the maximal input sizes
per GPU\@.  Also, setting up a huge number of pointers for each
batched \texttt{GEMM} call, as well as the redundant stores to and
loads from the GPU RAM of the contents of the shared memory makes the
prototype implementation about 20 times slower on average than the
original one.  Therefore, this version of the algorithm is not
recommendable on the present hardware.
%
%
\section{The QR factorizations as an alternative to the main kernel's subphases~\textbf{1} and \textbf{2} for the ill-conditioned input matrices}\label{s:E}
%
%
As mentioned in subsection~\ref{ss:2.3}, for a highly ill-conditioned
input matrix the formation of the Grammians of its block-column pairs
by the main kernel's subphase~\textbf{1} might result in several
numerically singular or indefinite blocks, on which the Cholesky
factorizations in the subphase~\textbf{2} are bound to fail.  In such
cases the affected (or, for simplicity, all) block-column pairs have
to be shortened by the QR factorizations instead.  Ideally, the
factorizations would employ column pivoting (see,
e.g.,~\cite{Quintana-Orti-Sun-Bischof-98} and the LAPACK's
\texttt{XGEQP3} routines,
$\text{\texttt{X}}\in\{\text{\texttt{D}},\text{\texttt{Z}}\}$), but it
can be inefficient on GPUs~\cite{Mary-et-al-2015}.

\looseness=-1
The batched non-pivoted QR factorization routine from cuBLAS could be
used with the memory layout proposed for Appendix~\ref{s:D}, where the
$\tau$ vectors would be stored in the $\mathrm{work}$ buffer.
Replacing the subphases~\textbf{1} and \textbf{2} that way would
require to pack the $\mathrm{actual}$ block columns, e.g.,
$F_{\mathtt{p}}^!$ and $F_{\mathtt{q}}^!$, consecutively in the
$\mathrm{shadow}$ copy $F^?$ for each $l$ (similarly for $G^!$).  The
packed block-column pairs in the $\mathrm{shadow}$ copy would then be
subject to the QR factorizations, and the upper triangular
$32\times 32$ factors would be extracted from it by the modified
subphase~\textbf{3}.  The postmultiplication subphase would remain as
in Appendix~\ref{s:D}.  However, this approach would also suffer from
the redundant traffic between the shared and the global memory.

Another approach is to keep the original memory layout and have two
``streaming'' QR factorizations (one for the assigned block-column
pair of $F$, and another for the same-indexed block-column pair of
$G$) embedded into the main kernel.  Each factorization reads from the
global memory only once, in chunks, never writes to it (since the
matrix $Q$ does not have to be applied again or restored), and leaves
the upper triangular factor in the shared memory where the
subphase~\textbf{3} expects it to be.  The full details of this
in-kernel QR factorization can be found
in~\cite[subsection~4.2]{Novakovic-2015} for double precision, while a
complex version would be implemented similarly.

Table~\ref{tbl:QRvsCholesky} demonstrates that close to eightfold
slowdown can be expected on a Volta GPU with the in-kernel QR
factorizations on larger matrices, even when using the simplest,
unsafe Frobenius norm computations inside them.
\begin{table}[h!bt]
  \caption{The wall times and their ratios of
    DHZ0-(\textsc{me}-\textsc{fb}-\textsc{me}) with the Cholesky
    factorizations of the Grammian matrices ($\blacksquare$) and with
    the QR factorizations of the block-column pairs ($\square$), both
    in-kernel, on the small real dataset.}
  \label{tbl:QRvsCholesky}
  \centering\addtolength{\tabcolsep}{-2.75pt}
  \begin{tabular}{@{}cccccccc@{}}
    \toprule
    $n$ & $\square$~[s] & $\blacksquare$~[s] & $\square/\blacksquare$ & $n$ & $\square$~[s] & $\blacksquare$~[s] & $\square/\blacksquare$\\
    \midrule
    \hphantom{0}512 & \hphantom{00}0.83 & \hphantom{0}0.20 & 4.1 & 5632 & \hphantom{0}437.48 & \hphantom{0}57.81 & 7.6\\
               1024 & \hphantom{00}3.57 & \hphantom{0}0.69 & 5.1 & 6144 & \hphantom{0}522.05 & \hphantom{0}68.10 & 7.7\\
               1536 & \hphantom{00}7.30 & \hphantom{0}1.24 & 5.9 & 6656 & \hphantom{0}544.53 & \hphantom{0}71.13 & 7.7\\
               2048 & \hphantom{0}12.87 & \hphantom{0}2.02 & 6.4 & 7168 & \hphantom{0}674.35 & \hphantom{0}87.38 & 7.7\\
               2560 & \hphantom{0}21.86 & \hphantom{0}3.31 & 6.6 & 7680 & \hphantom{0}775.09 &            100.64 & 7.7\\
               3072 & \hphantom{0}61.46 & \hphantom{0}8.69 & 7.1 & 8192 &            1015.82 &            128.65 & 7.9\\
               3584 & \hphantom{0}90.87 &            12.66 & 7.2 & 8704 &            1233.91 &            157.62 & 7.8\\
               4096 &            128.73 &            17.52 & 7.3 & 9216 &            1389.86 &            176.37 & 7.9\\
               4608 &            186.27 &            25.22 & 7.4 & 9728 &            1551.47 &            198.03 & 7.8\\
               5120 &            230.32 &            31.26 & 7.4 &  --  &               --   &              --   &  --\\
    \bottomrule
  \end{tabular}
\end{table}
%
%
\section{Performance-wise comparison of a CPU and a single-GPU implementation of the complex Hari--Zimmermann GSVD}\label{s:F}
%
%
Table~\ref{tbl:GPUvsCPU} shows the speedup of the complex single-GPU
algorithm on the small dataset versus the CPU GHSVD algorithm
from~\cite{Singer-DiNapoli-Novakovic-Caklovic-2020}, with the
signature matrix for the latter being $J=I$.  The GHSVD's
implementation was not tuned for the simpler GSVD computation.  If it
were, it could have been at least twice as fast on a more modern CPU,
but nevertheless noticeably slower than the GPU algorithm on the Volta
architecture.
\begin{table}[h!bt]
  \caption{The wall times and their ratios of
    ZHZ0-(\textsc{me}-\textsc{fb}-\textsc{me}) ($\bullet$) and a close
    CPU analogue of ZHZ6-(\textsc{me}-\textsc{bo}-\textsc{me}) with
    the variable block sizes and 64 threads (one per core) of an Intel
    Xeon Phi 7210 ($\circ$), on the small complex dataset.}
  \label{tbl:GPUvsCPU}
  \centering\addtolength{\tabcolsep}{-2.2875pt}
  \begin{tabular}{@{}cccccccc@{}}
    \toprule
    $n$ & $\circ$~[s] & $\bullet$~[s] & $\circ/\bullet$ & $n$ & $\circ$~[s] & $\bullet$~[s] & $\circ/\bullet$\\
    \midrule
    \hphantom{0}512 & \hphantom{00}1.22 & \hphantom{0}0.26 & \hphantom{0}4.7 & 5632 & \hphantom{0}422.37 & \hphantom{0}65.27 & \hphantom{0}6.5\\
               1024 & \hphantom{00}9.52 & \hphantom{0}0.81 &            11.8 & 6144 & \hphantom{0}671.15 & \hphantom{0}77.85 & \hphantom{0}8.6\\
               1536 & \hphantom{0}21.52 & \hphantom{0}1.68 &            12.8 & 6656 & \hphantom{0}659.89 & \hphantom{0}92.55 & \hphantom{0}7.1\\
               2048 & \hphantom{0}45.40 & \hphantom{0}2.98 &            15.3 & 7168 & \hphantom{0}756.02 &            109.74 & \hphantom{0}6.9\\
               2560 & \hphantom{0}70.81 & \hphantom{0}5.46 &            13.0 & 7680 & \hphantom{0}971.17 &            131.07 & \hphantom{0}7.4\\
               3072 &            114.68 &            13.30 & \hphantom{0}8.6 & 8192 &            1931.45 &            182.53 &            10.6\\
               3584 &            164.41 &            17.96 & \hphantom{0}9.2 & 8704 &            1258.28 &            208.32 & \hphantom{0}6.0\\
               4096 &            397.57 &            23.63 &            16.8 & 9216 &            1520.19 &            236.45 & \hphantom{0}6.4\\
               4608 &            270.74 &            30.77 & \hphantom{0}8.8 & 9728 &            1646.25 &            269.04 & \hphantom{0}6.1\\
               5120 &            360.70 &            40.11 & \hphantom{0}9.0 &  --  &               --   &              --   &  -- \\
    \bottomrule
  \end{tabular}
\end{table}
%
%
\section{Computational pitfalls of the GEVD of $(F^{\ast}F,G^{\ast}G)$ versus the GSVD of $(F,G)$}\label{s:G}
%
%
This small example shows what happens with the GEVD of
$(F^{\ast}F,G^{\ast}G)$ versus the GSVD of $(F,G)$ when one of the
matrices (here, $G$) is even mildly ill-conditioned.  All results were
obtained by Matlab R2020a Update 3.

Let $F$ and $G$, with $0<g_{11}\ll 1$, be given as
\begin{displaymath}
  F=\begin{bmatrix}
  1 & 1 & 1 & 1\\
  0 & 1 & 1 & 1\\
  0 & 0 & 1 & 1\\
  0 & 0 & 0 & 1
  \end{bmatrix},\quad
  G=\begin{bmatrix}
  g_{11} & 1 & 1 & 1\\
  0 & 1 & 1 & 1\\
  0 & 0 & 1 & 1\\
  0 & 0 & 0 & 1
  \end{bmatrix}.
\end{displaymath}
Clearly, $F$ and $G$ are non-singular, so $A\assgn F^T F$ and
$B\assgn G^T G$ are symmetric positive definite.  If, e.g.,
$g_{11}=10^{-10}$ (rounded to double precision) is taken, the
generalized eigenvalues of $(A,B)$, also rounded to double precision,
are
\begin{displaymath}
  \mathop{\Lambda}(\mathtt{A},\mathtt{B})=\begin{bmatrix*}[l]
  5.000000000500000\cdot 10^{-1}\\
  1.000000000000000\cdot 10^0\\
  1.000000000000000\cdot 10^0\\
  1.999999999800000\cdot 10^{20}
  \end{bmatrix*},
\end{displaymath}
where $\mathtt{A}$ and $\mathtt{B}$ are the Grammians of the rational
(symbolic) representations of $F$ and $G$, respectively, while
$\mathop{\Lambda}(\mathtt{A},\mathtt{B})$ is a sorted representation
of the symbolic output of
$\mathop{\mathtt{eig}}(\mathtt{A},\mathtt{B})$.

When computing $\mathop{\Lambda}(A,B)$ in double precision, two routes
can be taken.  The first one finds the generalized singular values of
$F$ and $G$, i.e.,
$\mathop{\Sigma}(F,G)\assgn\mathop{\mathtt{gsvd}}(F,G)$,
\begin{displaymath}
  \mathop{\Sigma}(F,G)=\begin{bmatrix*}[l]
  7.071067812219032\cdot 10^{-1}\\
  9.999999999999997\cdot 10^{-1}\\
  9.999999999999997\cdot 10^{-1}\\
  1.414213562302384\cdot 10^{10}
  \end{bmatrix*},
\end{displaymath}
and squares them to get the generalized eigenvalues $\widehat{\Lambda}$,
\begin{displaymath}
  \mathop{\widehat{\Lambda}}(A,B)=\begin{bmatrix*}[l]
  5.000000000500004\cdot 10^{-1}\\
  9.999999999999993\cdot 10^{-1}\\
  9.999999999999993\cdot 10^{-1}\\
  1.999999999800000\cdot 10^{20}
  \end{bmatrix*},
\end{displaymath}
which are close to $\mathop{\Lambda}(A,B)$, i.e., the exact ones.

The second route computes
$\mathop{\widetilde{\Lambda}}(A,B)\assgn\mathop{\mathtt{eig}}(A,B)$
from $A$ and $B$ by a ``proper'' GEVD routine, getting
\begin{displaymath}
  \mathop{\widetilde{\Lambda}}(A,B)=\begin{bmatrix*}[l]
  0.000000000000000\cdot 10^0\\
  0.000000000000000\cdot 10^0\\
  0.000000000000000\cdot 10^0\\
  1.999999999800000\cdot 10^{20}
  \end{bmatrix*},
\end{displaymath}
with the lower three generalized eigenvalues vanishing.

The spectral condition numbers of the matrices are:
\begin{displaymath}
  \begin{aligned}
    \mathop{\kappa_2^{}}(F) & \approx 5.411474127809772\cdot 10^0\\
    \mathop{\kappa_2^{}}(G) & \approx 3.972825427374361\cdot 10^{10}\\
    \mathop{\kappa_2^{}}(A) & \approx 2.928405223595454\cdot 10^1 & =\mathop{\kappa_2^2}(F)\\
    \mathop{\kappa_2^{}}(B) & \approx 1.578344820257155\cdot 10^{21} & =\mathop{\kappa_2^2}(G)
  \end{aligned},
\end{displaymath}
i.e., the condition numbers of $A$ and $B$ are the squares of those of
$F$ and $G$, respectively, which is easily shown by taking the SVD of
$F$ or $G$ and expressing $A$ or $B$ in terms of it.  Therefore, a
mildly ill-conditioned matrix $G$ causes its Grammian $B$ to be
severely ill-conditioned and the GEVD fails, so the matrix
multiplications have to be avoided in such a case by computing the
GSVD of $(F,G)$ instead.

Motivated by the previous example, 16 pairs $(F,G_j)$ of square real
matrices of order $512$ were generated in quadruple precision by the
adapted symmetric indefinite factorization with complete
pivoting~\cite{Slapnicar-98}, and then rounded to double precision,
such that $A=F^T F$, with $\mathop{\kappa_2}(A)=10$, was fixed, while
$B_j^{}=G_j^T G_j^{}$ varied, with the prescribed values for
$\kappa_{B_j}=\mathop{\kappa_2}(B_j)=10^j$, for $1\le j\le 16$.  The
matrices $A$ and $B_j$ were generated from the given positive
eigenvalues by the LAPACK's testing routine \texttt{DLAROR} adapted to
quadruple precision, ensuring that their condition numbers are as
close to the prescribed ones as possible, and that their ensuing
factorizations were identical in effect to the Cholesky factorizations
with diagonal pivoting.  Then, the generalized singular values
$\widehat{\Sigma}_j$ of $(F_,G_j)$, i.e., the singular values of
$F G_j^{-1}$, were computed in high precision with 32 significant
$\mathtt{digits}$ as
\begin{displaymath}
  \widehat{\Sigma}_j=\mathop{\mathtt{double}}(\mathop{\mathtt{svd}}(\mathop{\mathtt{vpa}}(F)/\mathop{\mathtt{vpa}}(G_j))),
\end{displaymath}
by Matlab R2020b, to be compared against in two ways.

The generalized singular values $\Sigma_j$ of $(F,G_j)$ were
recomputed by the DHZ0-(\textsc{me}-\textsc{fb}-\textsc{me}) variant,
and the generalized eigenvalues $\Lambda_j$ of $(A,B_j)$ were obtained
by the \texttt{cusolverDnDsygvd} routine from the cuSOLVER library, to
assess how well they approximate squares of the generalized singular
values.  The maximal relative errors in the generalized singular
values computed by the last two approaches were then found, with
$1\le i\le 512$, as 
\begin{displaymath}
  \begin{aligned}
    \mathop{\mathrm{mre}}(\Sigma_j)&=\max_i\frac{|\mathop{\mathtt{ext}}(\sigma_{j,i})-\mathop{\mathtt{ext}}(\hat{\sigma}_{j,i})|}{\mathop{\mathtt{ext}}(\hat{\sigma}_{j,i})},\\
    \mathop{\mathrm{mre}}(\Lambda_j)&=\max_i\frac{|\sqrt{\mathop{\mathtt{ext}}(\lambda_{j,i})}-\mathop{\mathtt{ext}}(\hat{\sigma}_{j,i})|}{\mathop{\mathtt{ext}}(\hat{\sigma}_{j,i})},
  \end{aligned}
\end{displaymath}
where $\widehat{\Sigma}_j$, $\Sigma_j$, and $\Lambda_j$ were sorted
descendingly.  Here, $\mathtt{ext}$ denotes a conversion to the
Intel's extended 80-bit datatype.

As it can be seen in Table~\ref{tbl:GSVDvsGEVDa}, the generalized
singular values computed by taking square roots of the generalized
eigenvalues rapidly lose accuracy when the condition of $B_j$
increases beyond about $1/\sqrt{\varepsilon}$, where $\varepsilon$ is
the machine precision.  Entirely invalid results are possible,
indicated by $\mathop{\mathrm{mre}}(\Lambda_{16})=\mathtt{NaN}$, since
the smallest six generalized eigenvalues of $(A,B_{16})$ were computed
as negative and the square roots could not have been taken.  The
generalized singular values computed by a ``proper'' GSVD (the HZ
algorithm) remained relatively accurate in all cases, though.
\begin{table}[h!bt]
  \caption{Relative accuracy of the generalized singular values of
    $(F,G_j)$ computed by the GSVD and the GEVD algorithms.}
  \label{tbl:GSVDvsGEVDa}
  \centering\addtolength{\tabcolsep}{-3.085pt}
  \begin{tabular}{@{}cccccc@{}}
    \toprule
    $\kappa_B$ & $\mathop{\mathrm{mre}}(\Sigma)$ & $\mathop{\mathrm{mre}}(\Lambda)$ & $\kappa_B$ & $\mathop{\mathrm{mre}}(\Sigma)$ & $\mathop{\mathrm{mre}}(\Lambda)$\\
    \midrule
    $10^1$ & $8.9\cdot 10^{-16}$ & $8.1\cdot 10^{-16}$            & $10^{9\hphantom{0}}$ & $1.4\cdot 10^{-15}$ & $2.9\cdot 10^{-9}$\\
    $10^2$ & $9.3\cdot 10^{-16}$ & $5.8\cdot 10^{-15}$            & $10^{10}$            & $1.2\cdot 10^{-15}$ & $1.2\cdot 10^{-7}$\\
    $10^3$ & $1.2\cdot 10^{-15}$ & $9.1\cdot 10^{-15}$            & $10^{11}$            & $1.1\cdot 10^{-15}$ & $1.3\cdot 10^{-6}$\\
    $10^4$ & $2.2\cdot 10^{-15}$ & $7.1\cdot 10^{-14}$            & $10^{12}$            & $2.2\cdot 10^{-15}$ & $3.4\cdot 10^{-4}$\\
    $10^5$ & $1.2\cdot 10^{-15}$ & $7.6\cdot 10^{-13}$            & $10^{13}$            & $1.1\cdot 10^{-15}$ & $3.7\cdot 10^{-3}$\\
    $10^6$ & $2.2\cdot 10^{-15}$ & $1.1\cdot 10^{-11}$            & $10^{14}$            & $1.4\cdot 10^{-15}$ & $2.3\cdot 10^{-2}$\\
    $10^7$ & $2.3\cdot 10^{-15}$ & $8.4\cdot 10^{-11}$            & $10^{15}$            & $1.1\cdot 10^{-15}$ & $3.5\cdot 10^{-2}$\\
    $10^8$ & $1.3\cdot 10^{-15}$ & $1.5\cdot 10^{-9\hphantom{0}}$ & $10^{16}$            & $1.6\cdot 10^{-15}$ & $\mathtt{NaN}$\\
    \bottomrule
  \end{tabular}
\end{table}

However, if it is known in advance that a particular GSVD problem on
$(F,G)$ is extremely well conditioned, solving it as the GEVD of
$(A,B)$ might be significantly faster.  From $A Z = B Z \Lambda$,
$F Z = U \Sigma_F$, and $G Z = V \Sigma_G$, it follows that finding
$Z$ (the generalized eigenvalues) and $\Lambda$ (the generalized
eigenvectors) suffices to compute $U \Sigma_F$, $V \Sigma_G$, and
$\Sigma=\Lambda^{1/2}$.  Normalizing the columns of $U \Sigma_F$ both
$U$ and $\Sigma_F$ (with the extracted norms on its diagonal) are
obtained, and similarly for $V$ and $\Sigma_G$, but $U$ and $V$ might
not be as orthogonal as those that come from a ``proper'' GSVD\@.  For
the left generalized singular vectors to be computed, the original
matrices $F$ and $G$ have to be preserved, which, along with a
non-trivial amount of workspace needed for the GEVD routine call, more
than doubles the memory requirements compared to those of the
single-GPU Hari--Zimmermann GSVD algorithm.

The matrices $A$ and $B$ are formed on a GPU, each by a single call of
the cuBLAS' routine \texttt{cublasDsyrk} (or \texttt{cublasZherk} in
the complex case) from $F$ and $G$, respectively.  Then, the
cuSOLVER's (legacy) routine \texttt{cusolverDnDsygvd} (or
\texttt{cusolverDnZhegvd} in the complex case) solves the GEVD problem
on $(A,B)$.  For $F Z$ and $G Z$ multiplications the cuBLAS'
\texttt{cublasDgemm} (or \texttt{cublasZgemm} in the complex case)
routine can be used.  Finally, the normalizing, as described above,
can be performed by a custom kernel similar to \texttt{rescale}.

The computation as described above, but skipping the normalization,
was compared performance-wise to the
DHZ0-(\textsc{me}-\textsc{fb}-\textsc{me}) and the
ZHZ0-(\textsc{me}-\textsc{fb}-\textsc{me}) variants.  Let \texttt{D}
stand for the ratio of the wall times of the ``proper'' GSVD and the
GSVD-using-GEVD approaches in the real, and \texttt{Z} in the complex
case, on the respective small datasets.  Then, \texttt{D} and
\texttt{Z}, as shown in Table~\ref{tbl:GSVDvsGEVDt}, are relaxed upper
bounds on the expected speedup by the GEVD approach (the actual
speedup should be somewhat smaller) on the GPUs.  For comparison,
about $15$--$35\times$ speedup is expected on the modern CPUs in the
complex case~\cite{Singer-DiNapoli-Novakovic-Caklovic-2020} on a set
of small-to-medium sized matrices.  But it should never be understated
that there are hard limits to the applicability of the GEVD approach
to the GSVD when the condition of the problem is high or unknown in
advance.
\begin{table}[h!bt]
  \caption{The wall time ratios of the ``proper'' (HZ) GSVD and the
    partial GSVD-using-GEVD on a single Volta GPU in the real
    (\texttt{D}) and the complex (\texttt{Z}) cases, on the small
    datasets.}
  \label{tbl:GSVDvsGEVDt}
  \centering
  \begin{tabular}{@{}cccccc@{}}
    \toprule
    $n$ & \texttt{D}~[$\times$] & \texttt{Z}~[$\times$] & $n$ & \texttt{D}~[$\times$] & \texttt{Z}~[$\times$]\\
    \midrule
    \hphantom{0}512 & \hphantom{0}6.57 & \hphantom{0}8.28 & 5632 & 49.92 & 26.63\\
               1024 & \hphantom{0}9.39 &            10.49 & 6144 & 49.25 & 25.30\\
               1536 &            10.76 &            11.75 & 6656 & 41.63 & 24.02\\
               2048 &            11.92 &            12.53 & 7168 & 42.38 & 23.60\\
               2560 &            13.88 &            15.05 & 7680 & 41.82 & 23.45\\
               3072 &            26.32 &            24.72 & 8192 & 45.45 & 26.95\\
               3584 &            27.76 &            22.83 & 8704 & 46.42 & 25.92\\
               4096 &            30.06 &            21.93 & 9216 & 44.71 & 25.07\\
               4608 &            34.12 &            21.35 & 9728 & 43.41 & 24.50\\
               5120 &            34.15 &            21.53 &  --  &  --   &  --  \\
    \bottomrule
  \end{tabular}
\end{table}
\begin{biogs}
  \emph{Vedran Novakovi\'{c}\/} received his Ph.D.\ in Mathematics in
  2017 from University of Zagreb, Croatia, where he started his career
  as a teaching assistant.  He also worked at STFC Daresbury
  Laboratory, UK, as a computational scientist and at Universidad
  Jaime I, Spain.  He is interested in parallel algorithms of
  numerical linear algebra, especially eigenvalue and singular value
  algorithms.\\
  \includegraphics[height=2ex,bb=40 52 296 308]{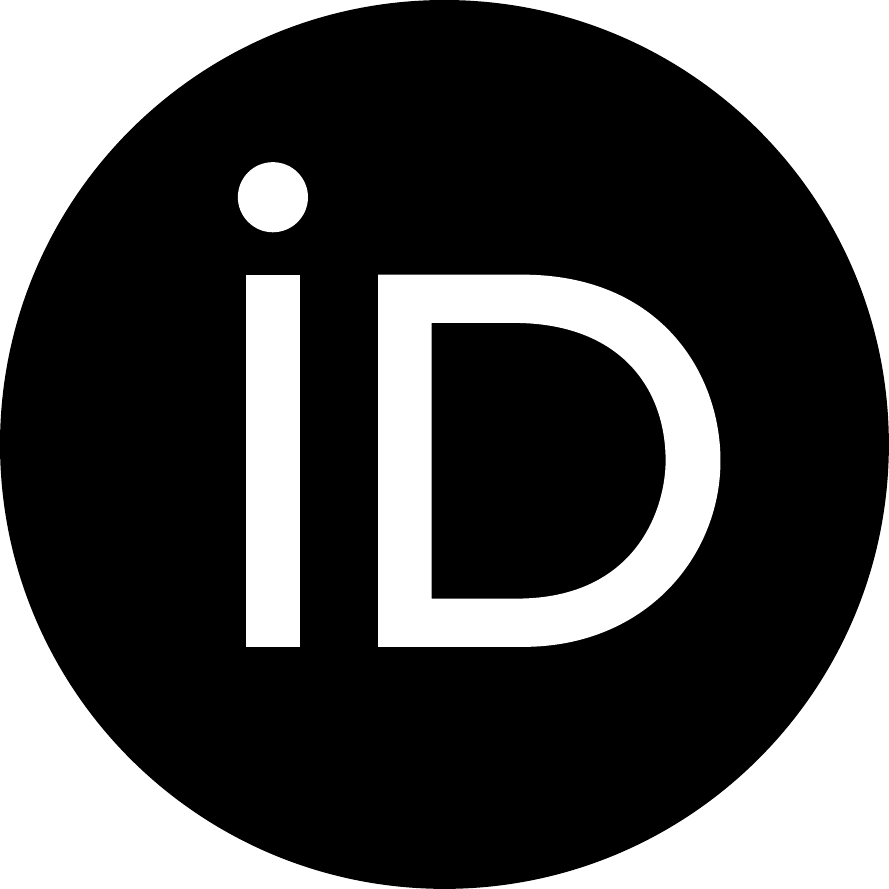}\,\url{https://orcid.org/0000-0003-2964-9674}

  \emph{Sanja Singer\/} is a tenured professor of Mathematics at the
  Faculty of Mechanical Engineering and Naval Architecture, University
  of Zagreb, Croatia.  She received her Ph.D.\ in Mathematics in 1997
  from the same University.  Her interests lie in accurate and
  high-performance algorithms of numerical linear algebra, especially
  matrix factorization algorithms, eigenvalue and singular value
  algorithms.\\
  \includegraphics[height=2ex,bb=40 52 296 308]{orcid.pdf}\,\url{https://orcid.org/0000-0002-4358-1840}
\end{biogs}
\end{document}